\documentclass[11pt,legal]{article}
\usepackage[utf8]{inputenc}
\usepackage{theorem,ifthen,algorithm,algorithmic}
\usepackage{amssymb,amsfonts,latexsym,dsfont}
\usepackage{tikz,pgflibraryplotmarks}
\usepackage{fullpage}
\usepackage{graphicx}
\usepackage{mathrsfs}
\usepackage{algorithm}
\usepackage{algorithmic}
\usepackage{amsmath}   
\usepackage{bm}
\usepackage{float}
\usepackage{caption}
\usepackage{hyperref}
\usepackage{epstopdf}

\newcommand{\bfb}{{\bm b}}

\newcommand{\bfe}{{\bm e}}

\newcommand{\bfs}{{\bf s}}
\newcommand{\bfx}{{\bm x}}
\newcommand{\bfy}{ {\bm y}}
\newcommand{\bfu}{{\bm u}}
\newcommand{\bfq}{{\bf q}}
\newcommand{\bfp}{{\bf p}}

\newcommand{\bfr}{{\bf r}}

\newcommand{\bfv}{{\bm v}}

\title{Iterative Methods at Lower Precision}

\author{Yizhou Chen, Xiaoyun Gong, Xiang Ji
\\Guest Editor: James Nagy\footnote{jnagy@emory.edu}}

\begin{document}
\maketitle

\begin{abstract}
\noindent Since numbers in the computer are represented with a fixed number of bits, loss of accuracy during calculation is unavoidable. At high precision where more bits (e.g. 64) are allocated to each number, round-off errors are typically small. On the other hand, calculating at lower precision, such as half (16 bits), has the advantage of being much faster. This research focuses on experimenting with arithmetic at different precision levels for large-scale inverse problems, which are represented by linear systems with ill-conditioned matrices. We modified the Conjugate Gradient Method for Least Squares (CGLS) and the Chebyshev Semi-Iterative Method (CS) with Tikhonov regularization to do arithmetic at lower precision using the MATLAB \textbf{chop} function, and we ran experiments on applications from image processing and compared their performance at different precision levels. We concluded that CGLS is a more stable algorithm, but overflows easily due to the computation of inner products, while CS is less likely to overflow but it has more erratic convergence behavior. When the noise level is high, CS outperforms CGLS by being able to run more iterations before overflow occurs; when the noise level is close to zero, CS appears to be more susceptible to accumulation of round-off errors.
\end{abstract}

\section{Introduction}

Most computer processors today use double-precision binary floating point arithmetic, which represents floating-point numbers with 64 bits. The convention follows from the IEEE-754 standard established in 1985, specifying the number of bits for the sign, exponent, and mantissa (fraction) for various floating-point formats including binary16 (half precision), binary32 (single precision), and binary64 (double precision). The difference between these formats is illustrated in the diagram below \cite{nvidiatechnicalblog_2021}:
\\
 \begin{center}
    \includegraphics[width=8cm]{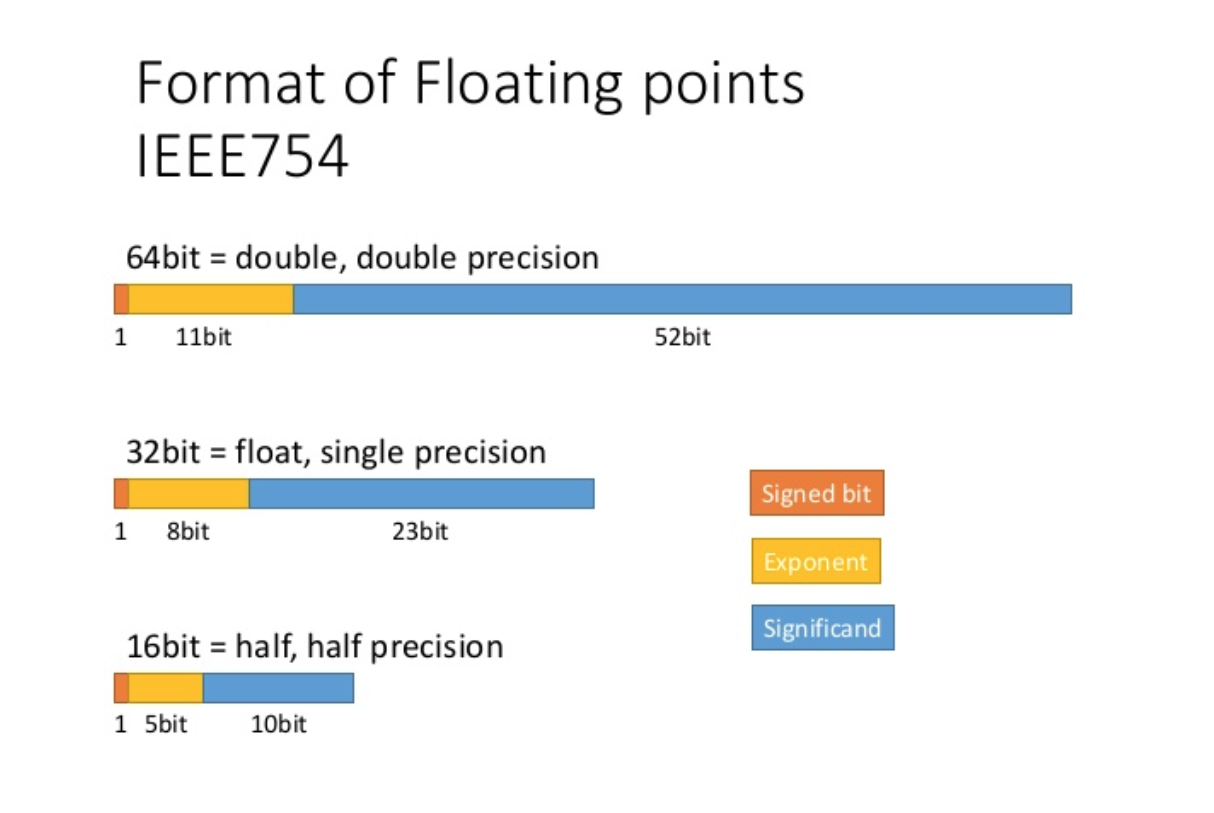}
 \end{center}
Taking a quarter of memory of their traditional 64-bit counterparts, computing with half precision numbers has attracted interest from researchers and manufacturers because of the potential for significantly reduced computational time. The format accelerates deep learning training, allows larger models, and improves gaming latency, providing players with better experiences \cite{san2021low}. However, aside from these benefits, half precision also brings one obvious disadvantage: increased inaccuracies in number representation. One goal of this research project is to investigate the impact of transferring from double precision operations to low precision arithmetic when solving large-scale inverse problems.\\
\\
Due to their wide variety of applications, we specifically utilized and modified codes for iterative methods in solving inverse problems to evaluate the performance of low precision algorithms. Inverse problems arise when using outside measurements to acquire information about internal or hidden data \cite{hansen2010discrete}. We operate on known outputs with some errors to compute the true input value. For example, X-ray computed tomography is a contactless imaging method for reconstructing target objects, most commonly, pathologies in the human body. Another example is image deblurring problems, which occur when the true picture is to be reconstructed from its blurry (and sometimes noisy) version \cite{Eps08}. These cases can be abstracted as linear systems $\bfb = A\bfx + \bfe$, where $A$ is a large-scale, typically ill-conditioned matrix and $\bfb$ is a vector output blended with noise, $\bfe$.\\
\\
A naive solution to this problem can be obtained by directly calculating a solution to $A\bfx = \bfb$. However, this naive solution is often corrupted by noise due to the ill-conditioning of matrix $A$. Regularization methods are often needed to balance signal and noise. One approach for regularization uses the singular value decomposition of $A$. However, when the matrix $A$ is large, as is often the case for many real-life applications (including our test cases), such direct methods may be difficult to implement, since decomposing the matrix can be very computationally costly. When the matrix is sparse, meaning it contains a lot of zero entries, computing matrix-vector products can be done very efficiently. Therefore, for such problems, implementing iterative methods, which only require matrix-vector products and vector operations, is a more efficient way to solve inverse problems than direct methods, with the former requiring $O(n^3)$ work for an $n \times n$ matrix $A$, while the latter involves significantly less than $O(n^2)$ work per iteration \cite{trefethen1997numerical}.

\section{Simulating Low Precision}
\subsection{Chop}

The MATLAB function \textbf{chop} written by Higham and Pranesh provides users a way to simulate low precision numbers and arithmetic. The function can simulate precisions, such as fp16, bfloat16, as well as custom formats where the user is able to choose the number of bits in the significand and the maximum value of the exponent. The function doesn't create the new data type; instead, it keeps numbers in double or single precision and makes sure that they have the right scale as in lower precision. Therefore, once the rounding meets the bit limit of the corresponding target precision, the remaining bits are set to zero. Computer operations are also carried out in double or single precision, and after that, the result is rounded to low precision \cite{higham2019simulating}. 
\\
\\
Fully simulating computations in low precision requires users to call \textbf{chop} after each arithmetic operation. For example, for the operation
$$a = x + y \times z,$$ the appropriate use of \textbf{chop} is:
$$ a = chop(x + chop(y \times z)).$$ This is burdensome yet unavoidable, so our modified version of iterative methods using \textbf{chop}, which only simulates low precision arithmetic, usually takes a long time to run, despite the fact that low precision arithmetic would run significantly faster with the corresponding hardware. However, for vector or matrix operations, there are ways to reduce the number of calls to \textbf{chop}, improving the efficiency of our codes. For example, for the vector inner product between $\bfx$ and $\bfy$, instead of using the line:
\begin{verbatim}
sum = 0;
for i = 1:vector_length
    sum = chop(sum + chop(x(i) * y(i)));
end
\end{verbatim}
we use an element-wise operation at the beginning:
 \begin{verbatim}
sum = 0;
z = chop(x.*y)
for i = 1:vector_length
    sum = chop(sum + z(i));
end
\end{verbatim}  
such that we successfully reduce $2n$ calls of \textbf{chop} to only $n+1$ calls, accelerating our running process.  
\subsection{Blocking}
\begin{figure}
\includegraphics[width=15cm]{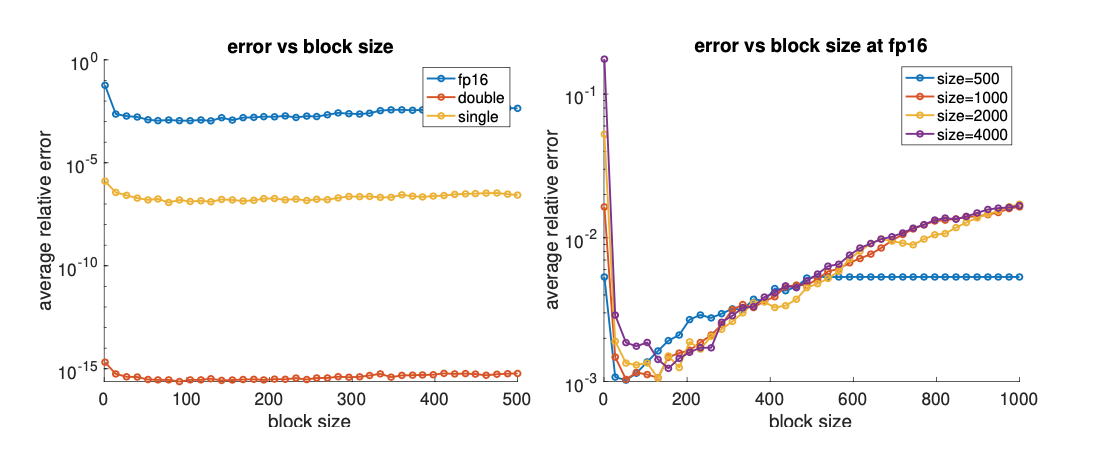}
\caption{Relationship between errors and block size.}\label{fig:error_block}
\end{figure}

Since we are using floating-point arithmetic, inaccuracy is inevitable in computations. However, blocking can be used to reduce the error bound. The method breaks a large number of operations into several smaller pieces, computes them independently, and sums them up. Consider the inner product between two vectors:
$$
\bfx = [x_1,x_2,x_3,x_4,x_5,x_6], \;\; 
\bfy = [y_1,y_2,y_3,y_4,y_5,y_6]
$$
Instead of calculating it directly, we could break it into:
$$
\bfx_1 = [x_1,x_2,x_3],\; \; \bfx_2 = [x_4,x_5,x_6]
$$
$$
\bfy_1 = [y_1,y_2,y_3],\; \; \bfy_2 = [y_4,y_5,y_6]
$$
Then we calculate $\bfx_1^T\bfy_1$ and $\bfx_2^T\bfy_2$, and add the sum. The result can have less errors than the direct calculation. The error bound for inner products $\bfx^T\bfy$ is \cite{higham2002accuracy}:
$$
|\bfx^T\bfy-fl(\bfx^T\bfy)| \leq \gamma_n|\bfx|^T|\bfy|
$$
where $\displaystyle \gamma_n = \frac{nu}{1-nu}$ and $u$ is the unit round-off of floating point arithmetic (e.g. $9.77*10^{-4}$ for half precision, $1.19*10^{-7}$ for single precision, and $2.22*10^{-16}$ for double precision).
Nonetheless, with blocking, the error bound is reduced to \cite{higham2002accuracy}:
$$
|s_n-\hat{s}_n| \leq \gamma_{(\log_2 n)+1}|\bfx|^T|\bfy|.
$$
 $s_n$ is the summation of the results for each block using the true value of $\bfx$ and $\bfy$, whereas $\hat s_n$ is the summation of each block's result using the floating number representation of $\bfx$ and $\bfy$.
The detailed proof is in \cite[Chapter 3]{higham2002accuracy}.
We plotted the error graph for double, single, and half precision against block size when computing the inner product of two random vectors in Figure \ref{fig:error_block}. We did inner products for each precision $20$ times and calculated the average error, which is computed as the difference between the result from our chopped version of inner products and MATLAB's default double precision computation.
\\

\noindent The graph on the left-hand-side includes all three precision levels, and we can see that the error for half precision is the largest, since it has the least bits. Focusing on half precision only, the graph on the right shows that errors decrease sharply as blocking is introduced. However, the errors increase as the block size increases, indicating that an appropriate block size needs to be carefully chosen. The reason behind this increase is that as the block size grows larger, it has the similar effect as doing no blocking; the whole vector or matrix is put into the first block, not divided into smaller sections. In our modified codes, we chose the block size to be $256$, appropriate for the $4096 \times 4096$ matrices we used in the test cases.

\section{Conjugate Gradient Method for Least Squares}

\subsection{Method Overview}
The Conjugate Gradient (CG) Method was introduced by Hestenes and Stiefel ~\cite{hestenes1952methods}. It is an iterative method for solving the linear system $A\bfx=\bfb$, where $A$ is a symmetric and positive definite matrix. The CG method can be viewed as an optimization problem of minimizing a convex quadratic function: $$\phi(\bfx)=\frac{1}{2}\bfx^{T}A\bfx-\bfx^{T}\bfb.$$ The gradient is zero at minimum point, meaning $\nabla \phi (\bfx) = A\bfx-\bfb=0$, which is exactly the linear system we are trying to solve.\\
\\
The CG method is a Krylov subspace method, which means its approximated solutions lie in Krylov subspaces. In each iteration, $\bfx$ is allowed to explore subspaces of increasing dimensions. An interesting property of CG is that during each iteration, $\bfx$ is updated to the point within the subspaces where the $A$-norm of the error is minimized \cite{trefethen1997numerical}. If we assume zero round-off errors, CG is guaranteed to converge within a limited number of steps. Specifically, if $A$ is an $n \times n$ matrix, the method will find the solution in no greater than $n$ steps. Besides, the updated $\bfx$ is closer to the true solution than the previous $\bfx$ in each iteration.\\
\\
The CGLS algorithm is the least squares version of the CG method, applied to the normal equations $A^{T}A\bfx = A^{T}\bfb$.\\
\\
One potential problem with this method for low precision is that the calculation of inner products can easily result in overflow, as illustrated in the experiment section below.
\\
\\The following algorithm describes CGLS \cite{bjorck1996numerical}.

\begin{algorithm}
\caption{Conjugate Gradient Method Least Squares}\label{alg:cap}
\begin{algorithmic}
\STATE Let $\bfx^{(0)}$ be an initial approximation, set $\bfr^{(0)} = \bfb - A\bfx^{(0)}, \bfp^{(0)} = \bfs^{(0)} = A^T\bfr^{(0)}, \psi_0 = ||\bfs^{(0)}||_2^2$
\WHILE{$\psi_k > tol$} 
    \STATE $\bfq^{(k)} = A\bfp^{(k)}$,
    \STATE $\alpha_k = \psi_k/||\bfq^{(k)}||_2^2$,
    \STATE $\bfx^{(k+1)} = \bfx^{(k)} + \alpha_k\bfp^{(k)}$,
    \STATE $\bfr^{(k+1)} = \bfr^{(k)} + \alpha_k\bfq^{(k)}$,
    \STATE $\bfs^{(k+1)} = A^T\bfr^{(k+1)}$,
    \STATE $\psi_{k+1} = ||\bfs^{(k+1)}||_2^2$,
    \STATE $\beta_k = \psi_{k+1}/\psi_k$,
    \STATE $\bfp^{(k+1)} = \bfs^{(k+1)} + \beta_k\bfp^{(k)}$.
\ENDWHILE
\end{algorithmic}
\end{algorithm}

\subsection{Experiment}
With the \textbf{chop} function, we modified the original CGLS method in the IRtools package, which offers various iterative methods for large-scale, ill-posed inverse problems and a set of test problems for these iterative methods \cite{gazzola2019ir}. Then we used two large-scale, ill-posed inverse problems from the same package: image deblurring and tomography reconstruction. The first is to reconstruct an approximation of the true image from the observed blurry version, whereas the second one is to reconstruct an image from measured projections, which can be obtained, for example, from X-ray beams. We investigated these two test problems in different sizes, floating-point precision levels, and noise levels. In each of these examples, the problem is modeled as $\bfb = A\bfx + \bfe$, where $A$ and $\bfb$ are given.

\subsubsection{Image Deblurring}
\begin{figure}[!htb]
\minipage{0.32\textwidth}
  \includegraphics[width=\linewidth]{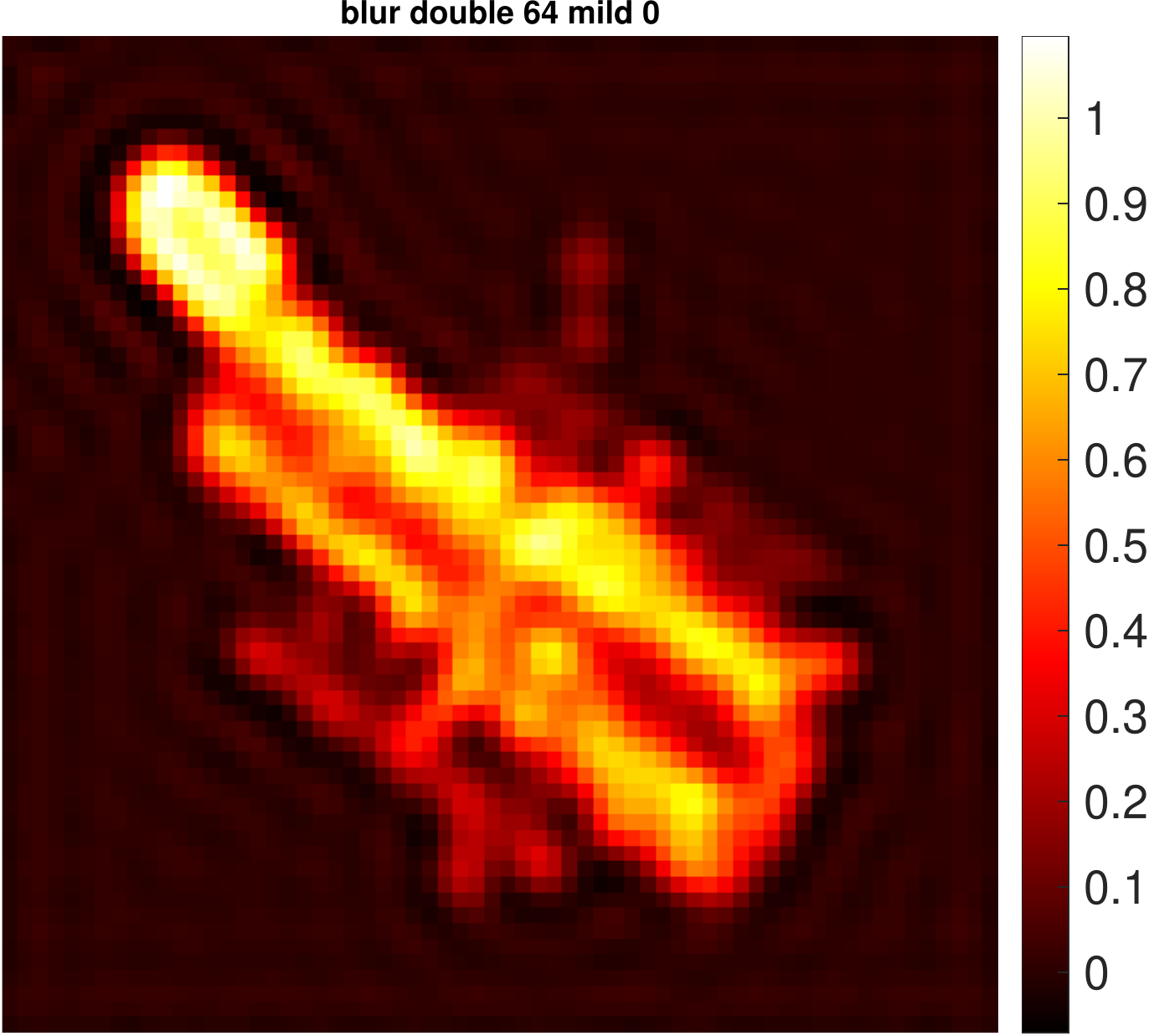}
  \caption{Double precision, size 64, zero noise.}\label{fig:Bdouble_precision_size64_zeroNoise}
\endminipage\hfill
\minipage{0.32\textwidth}
  \includegraphics[width=\linewidth]{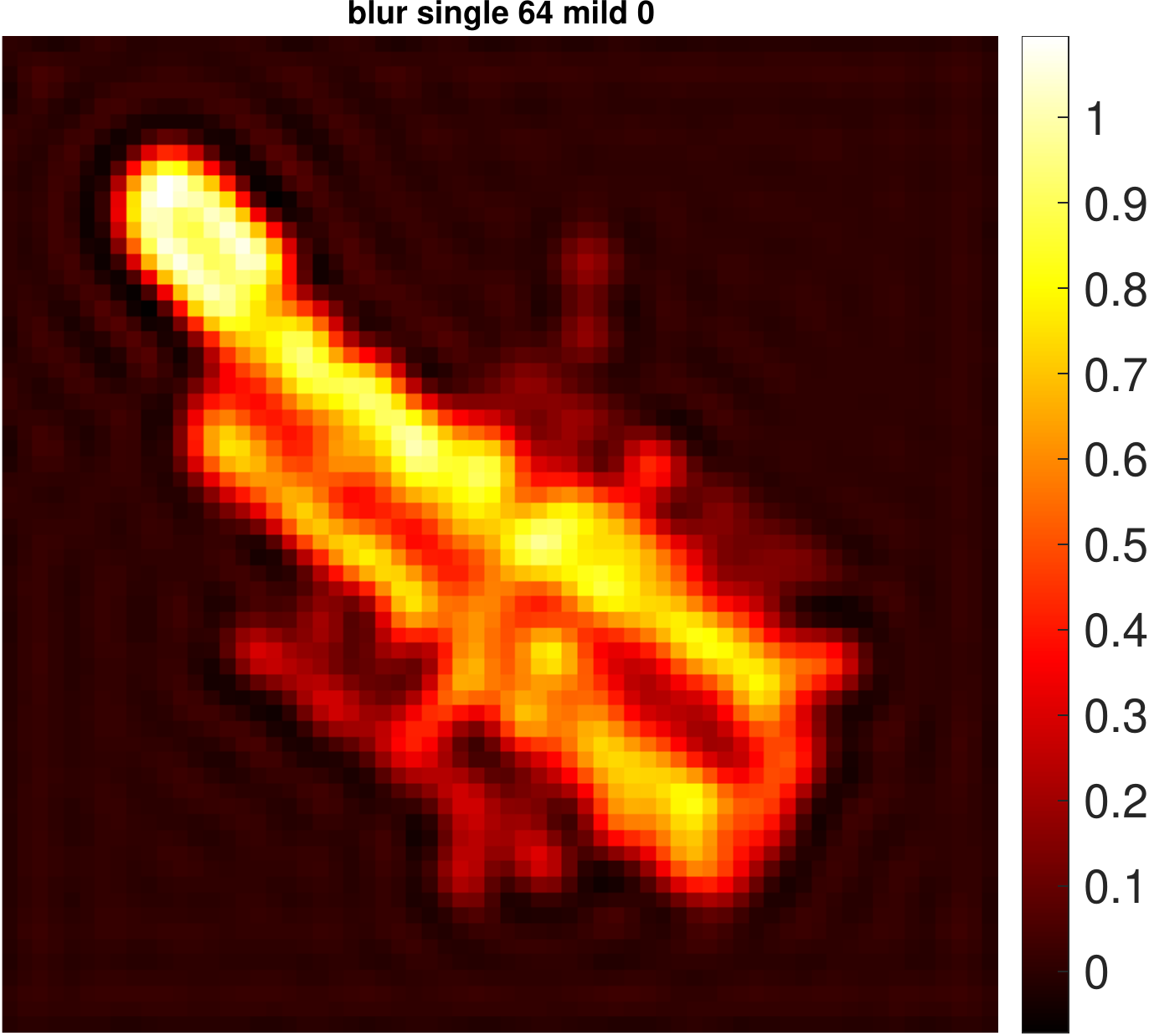}
  \caption{Single precision, size 64, zero noise.}\label{fig:Bsingle_precision_size64_zeroNoise}
\endminipage\hfill
\minipage{0.32\textwidth}%
  \includegraphics[width=\linewidth]{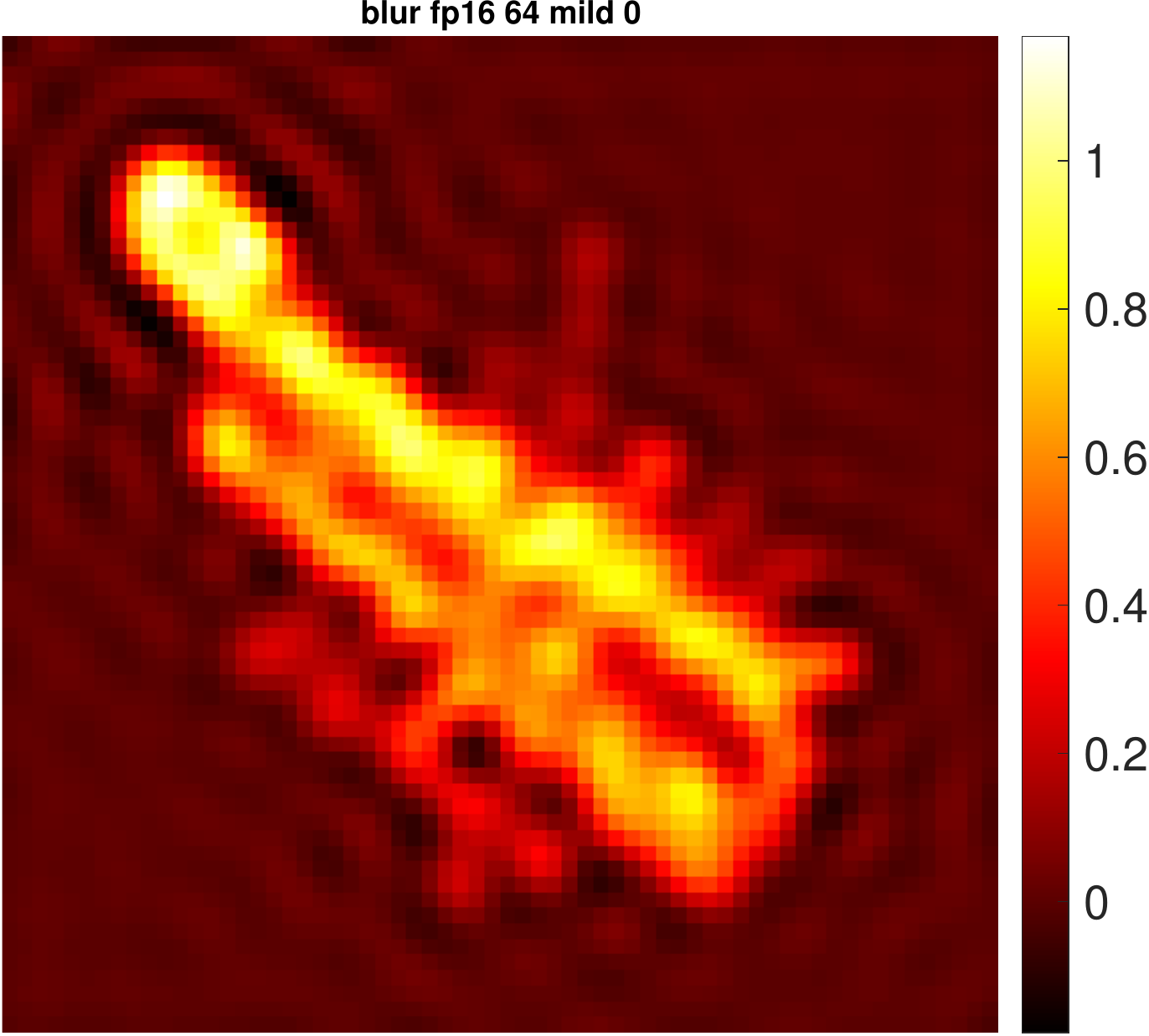}
  \caption{Half precision, size 64, zero noise.}\label{fig:Bhalf_precision_size64_zeroNoise}
\endminipage
\end{figure}

\begin{figure}
\begin{center}
\includegraphics[width=8cm]{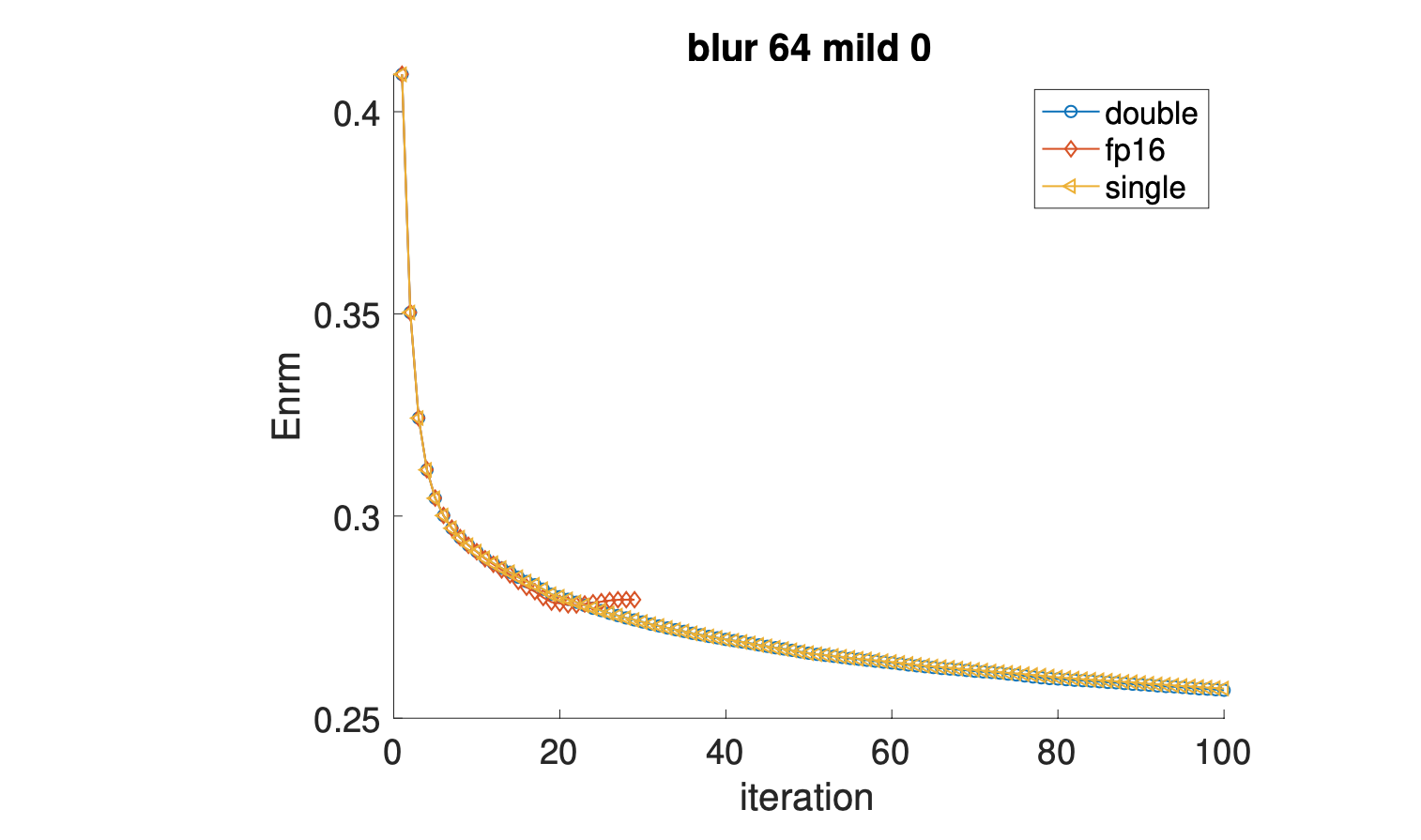}
\caption{The error norm of a size 64 problem with mild blurring of different precisions.}\label{fig:Benrms_size64_zeroNoise}
\end{center}
\end{figure}

Here we generate an image deblurring test problem using the \textbf{PRblur} function in the IRtools package. The true image is a picture of the Hubble space telescope, and the observed data is corrupted by Gaussian blurring function; for details, see \cite{gazzola2019ir}. We first added no noise to $\bfb$ (i.e. $\bfe = 0$) to see how our modified CGLS method works. In Figures \ref{fig:Bdouble_precision_size64_zeroNoise}, \ref{fig:Bsingle_precision_size64_zeroNoise}, and \ref{fig:Bhalf_precision_size64_zeroNoise}, we displayed the computed approximation of $\bfx$ at the final iteration.
\\
\\The single-precision result is very similar to that computed in double-precision, but the half-precision result is more blurry, and the background contains more artifacts. We also plotted the error graphs for all three formats in Figure \ref{fig:Benrms_size64_zeroNoise}.
\\

\begin{figure}[!htb]
\minipage{0.32\textwidth}
  \includegraphics[width=\linewidth]{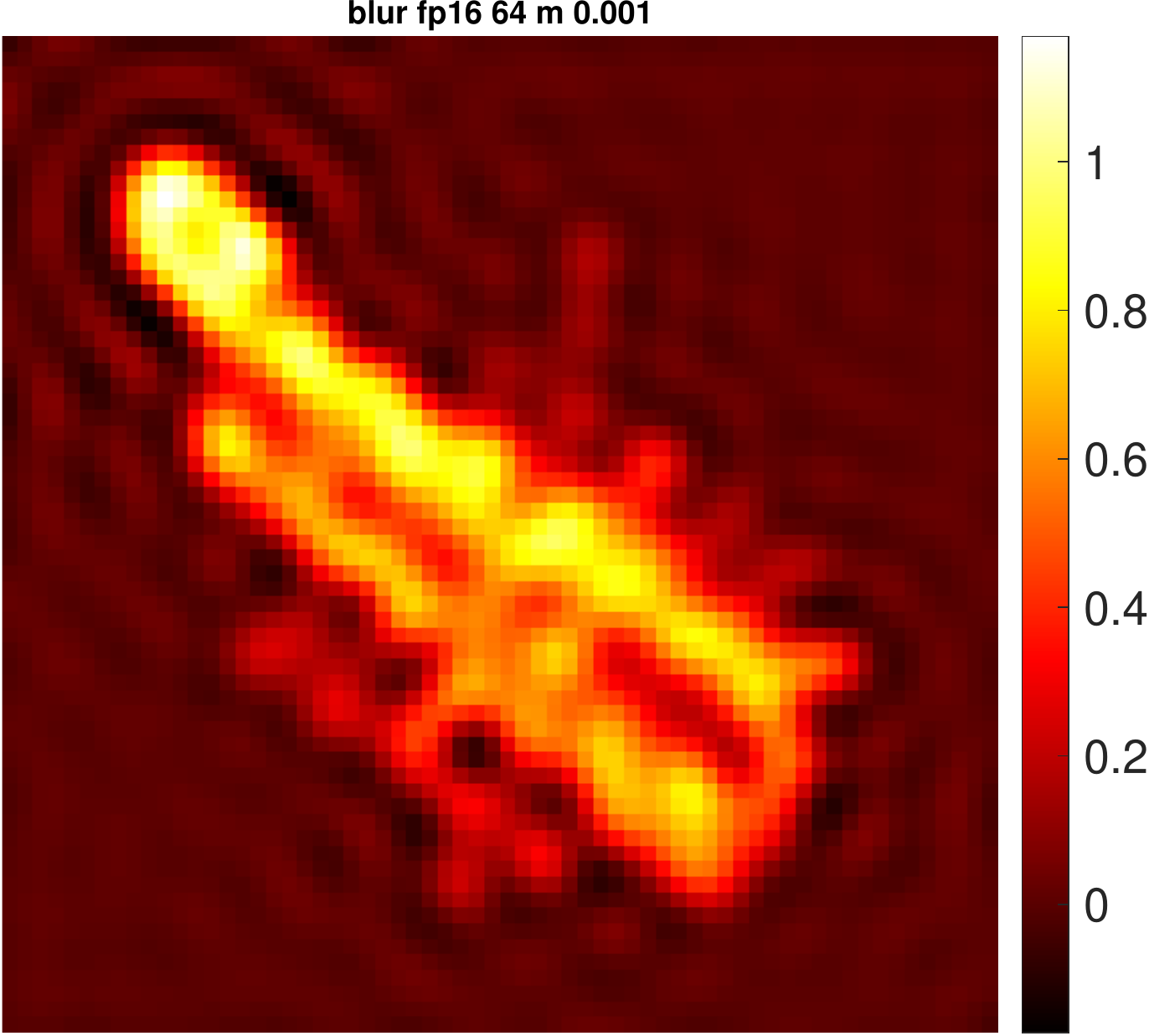}
  \caption{Half precision, size 64, 0.1\% noise.}\label{fig:Bhalf__size64_0.1PNoise}
\endminipage\hfill
\minipage{0.32\textwidth}
  \includegraphics[width=\linewidth]{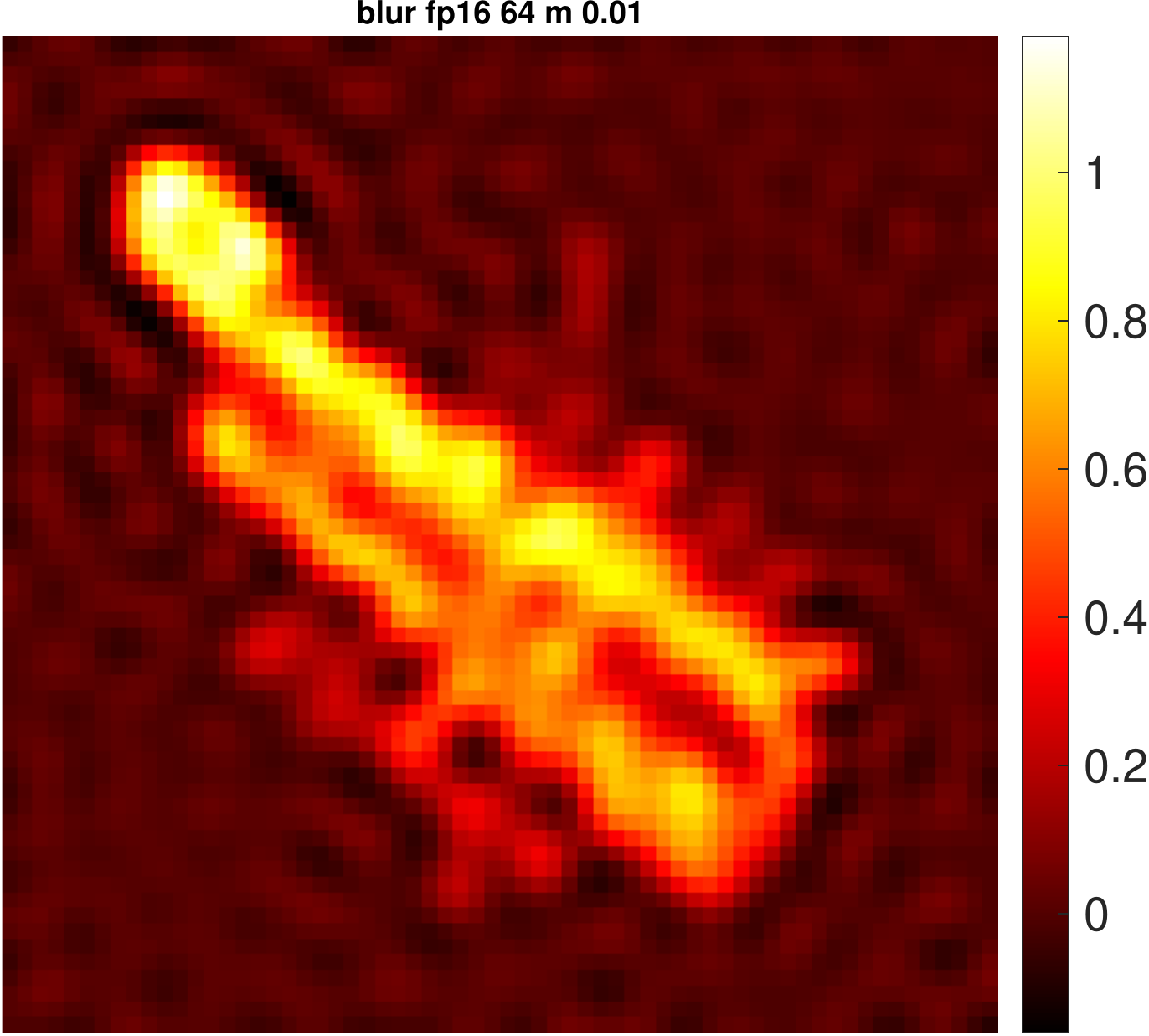}
  \caption{Half precision, size 64, 1\% noise.}\label{fig:Bhalf_size64_1PNoise}
\endminipage\hfill
\minipage{0.32\textwidth}%
  \includegraphics[width=\linewidth]{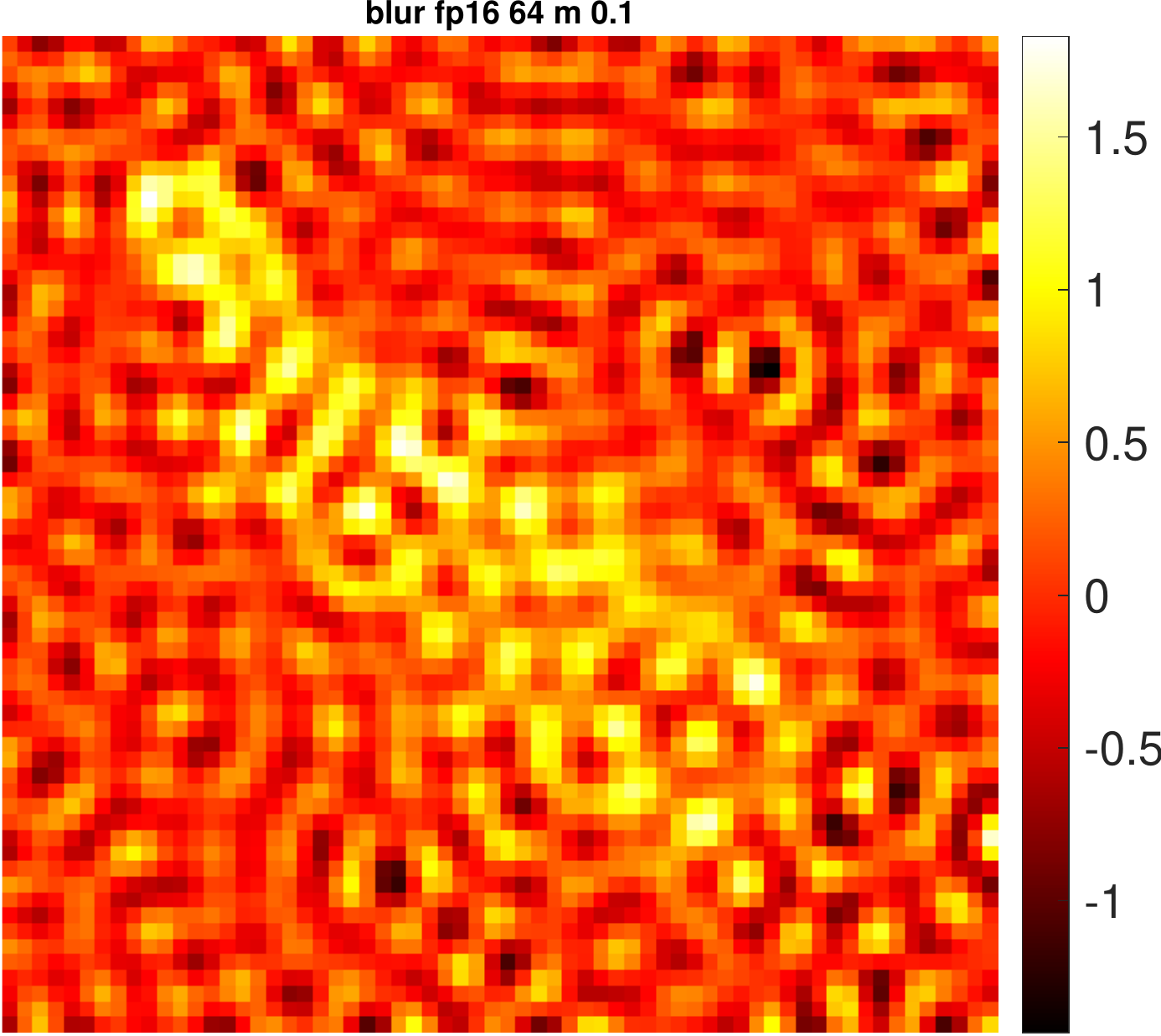}
  \caption{Half precision, size 64, 10\% noise.}\label{fig:Bhalf_size64_10PNoise}
\endminipage
\end{figure}

\noindent All three error norms overlap from the beginning until around $20^{th}$ iteration, where the half-precision errors begins to deviate from those in single and double precision. The difference is due to the round-off errors of half precision, which add up and take over. Moreover, the error norms of half precision terminates at $28^{th}$ iteration, because overflow of inner products causes NaNs (Not a Number) to be computed during the iteration.  \\
\\
In real life, the noise-free image, $\bfb$, is not often obtainable but most likely it is blended with noise. Therefore, it is imperative to make sure that the chopped CGLS method works effectively with noisy data as well. We added noise at different levels to $\bfb$ and displayed the computed approximations of $\bfx$ in the last iteration in half precision Figures \ref{fig:Bhalf__size64_0.1PNoise}, \ref{fig:Bhalf_size64_1PNoise}, and \ref{fig:Bhalf_size64_10PNoise}.\\

\begin{figure}[!htb]
\minipage{0.32\textwidth}
  \includegraphics[width=\linewidth]{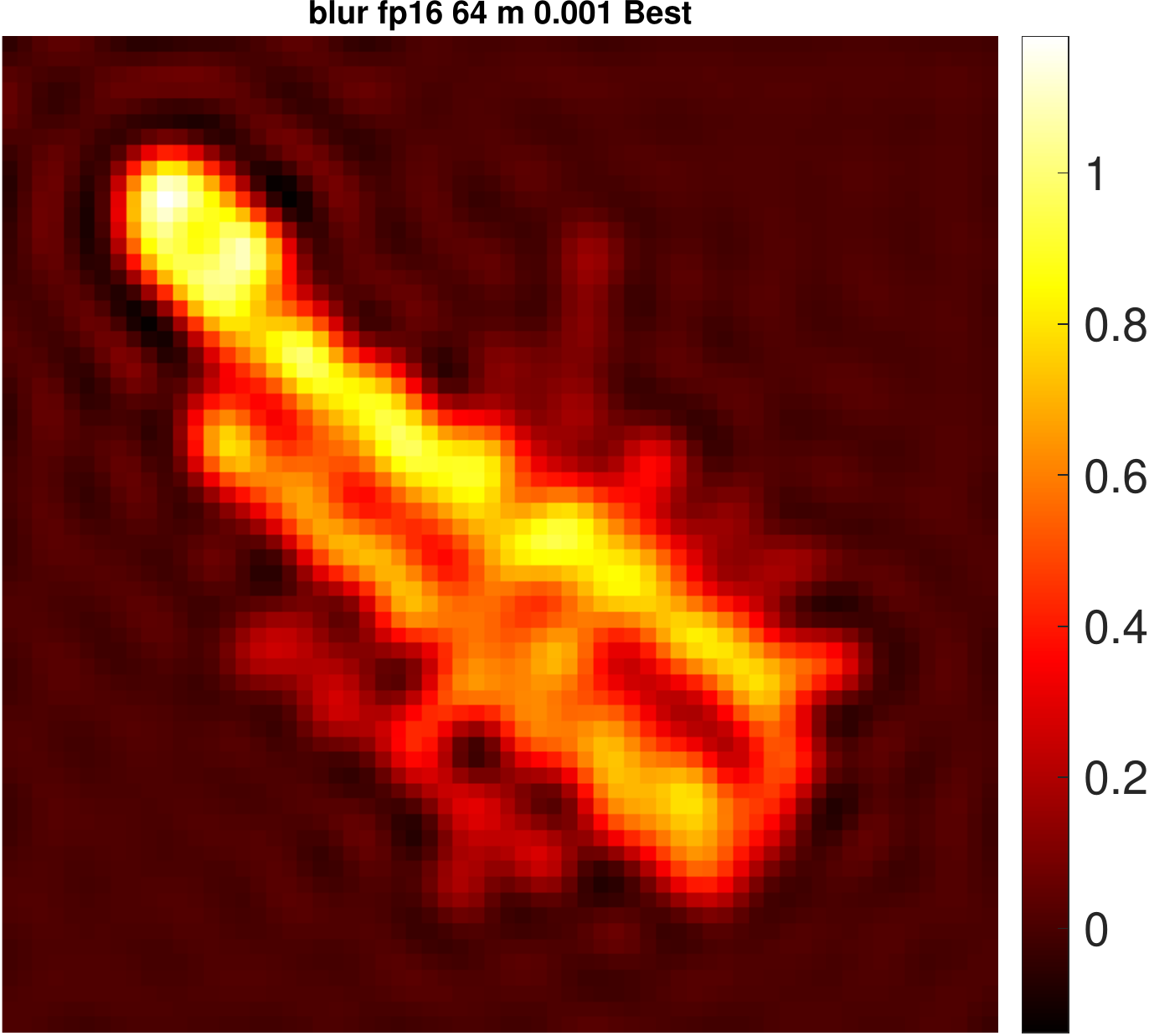}
  \caption{Half precision, problem size 64 with mild blurring and 0.1\% noise at best iteration.}\label{fig:cgls_fp16_64_m_0.001_Best}
\endminipage\hfill
\minipage{0.32\textwidth}
  \includegraphics[width=\linewidth]{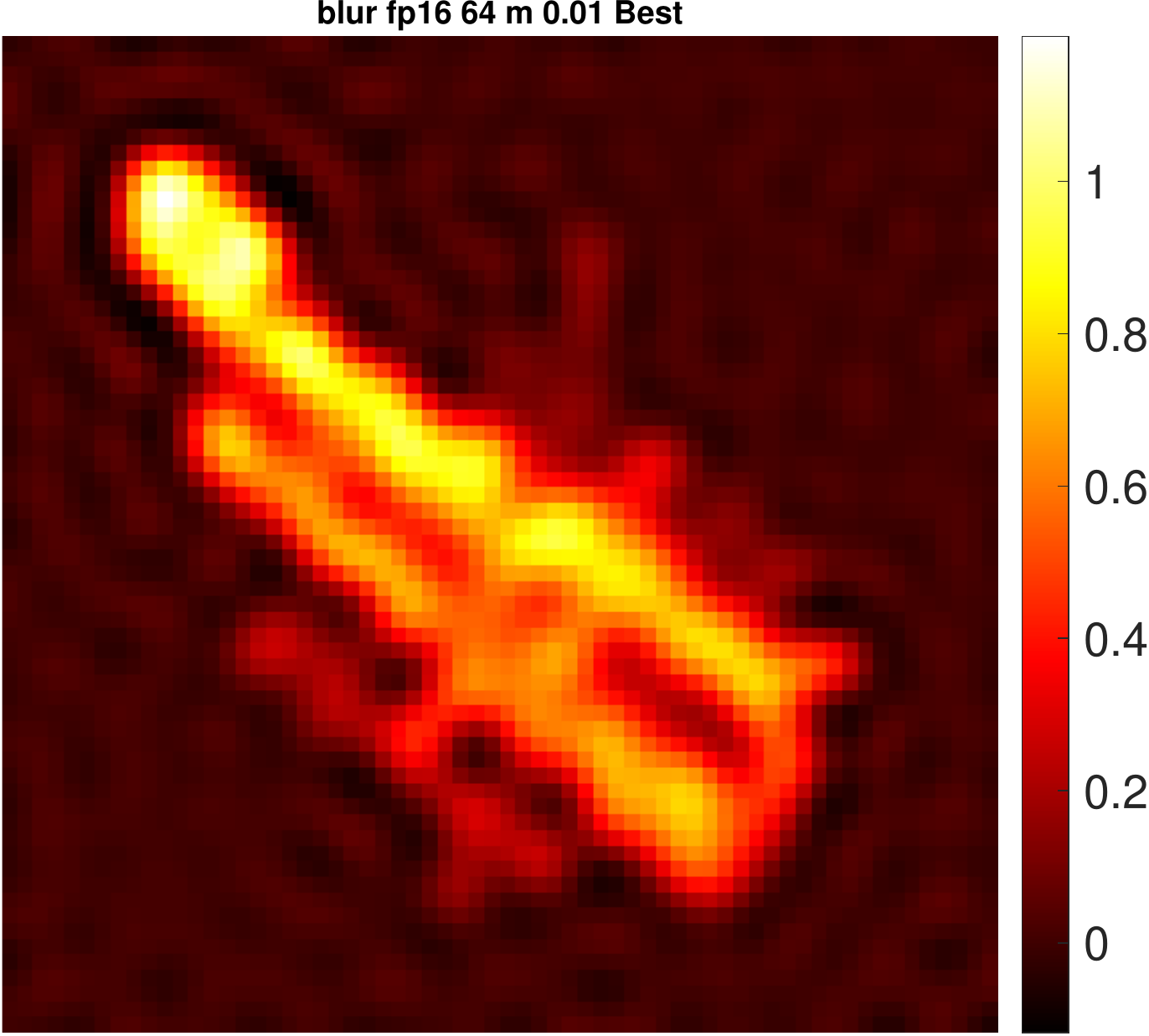}
  \caption{Half precision, problem size 64 with mild blurring and 1\% noise at best iteration.}\label{fig:cgls_fp16_64_m_0.01_Best}
\endminipage\hfill
\minipage{0.32\textwidth}%
  \includegraphics[width=\linewidth]{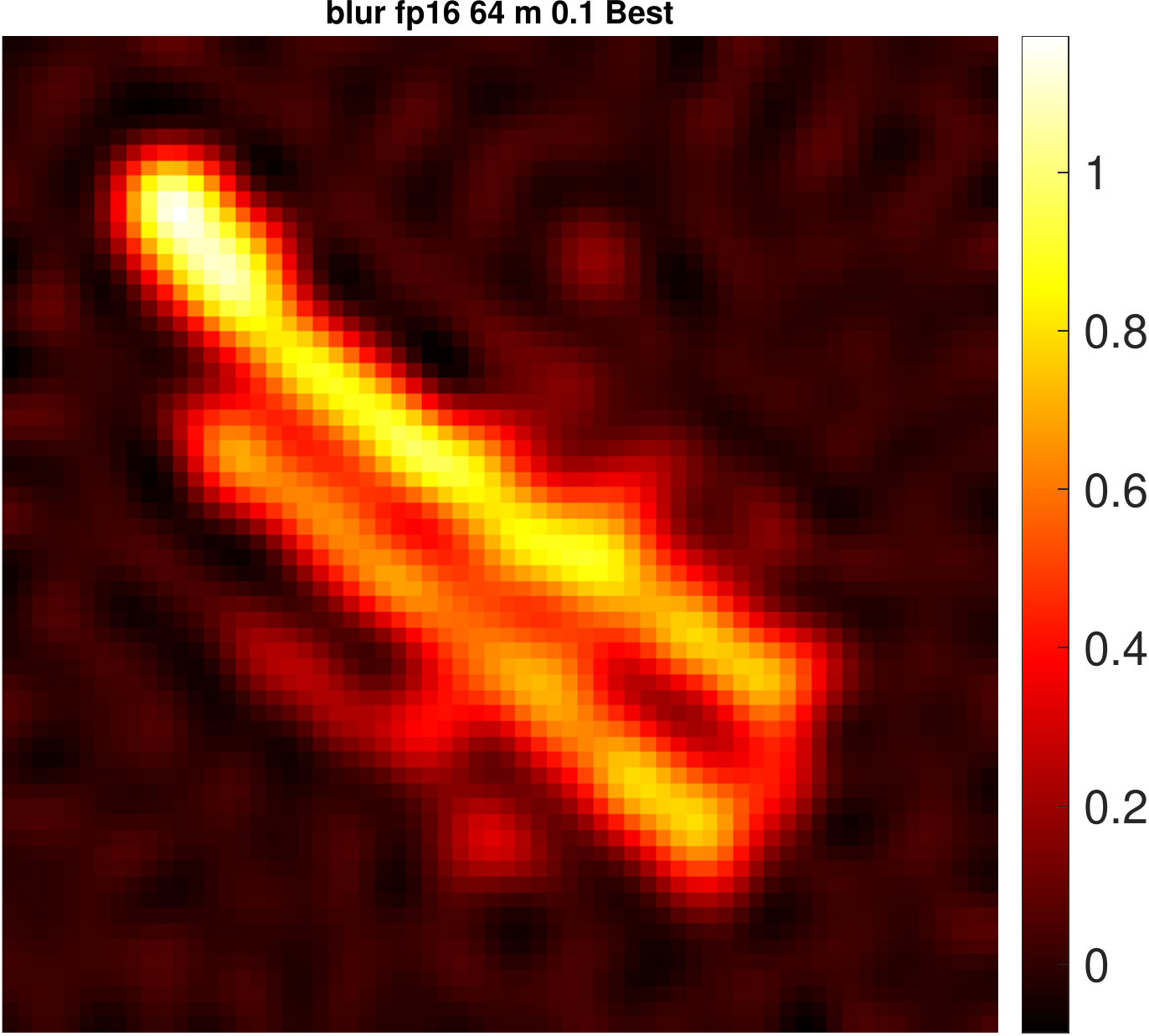}
  \caption{Half precision, problem size 64 with mild blurring and 10\% noise at best iteration.}\label{fig:cgls_fp16_64_m_0.1_Best}
\endminipage
\end{figure}

\noindent However, the last iteration doesn't necessarily mean it's the iteration with the best performance, so the results generated using $\bfx$ from the best iteration (i.e. the iteration with smallest error norms) are also shown in Figures \ref{fig:cgls_fp16_64_m_0.001_Best}, \ref{fig:cgls_fp16_64_m_0.01_Best}, and \ref{fig:cgls_fp16_64_m_0.1_Best}.
\noindent For results in the last iteration, while images obtained when using half precision are more distorted than their single counterparts, they follow the same trend in terms of the impact of noise. $0.1\%$ noise behaves very similarly to the cases with no noise; however, $1\%$ noise begins to dominate, and the background has substantially more artifacts, while the object is still identifiable. Eventually, the computed result is horribly corrupted by artifacts in the 10\% noise case; the picture contains very little meaningful information.\\
\begin{figure}[!htb]
\minipage{0.32\textwidth}
  \includegraphics[width=\linewidth]{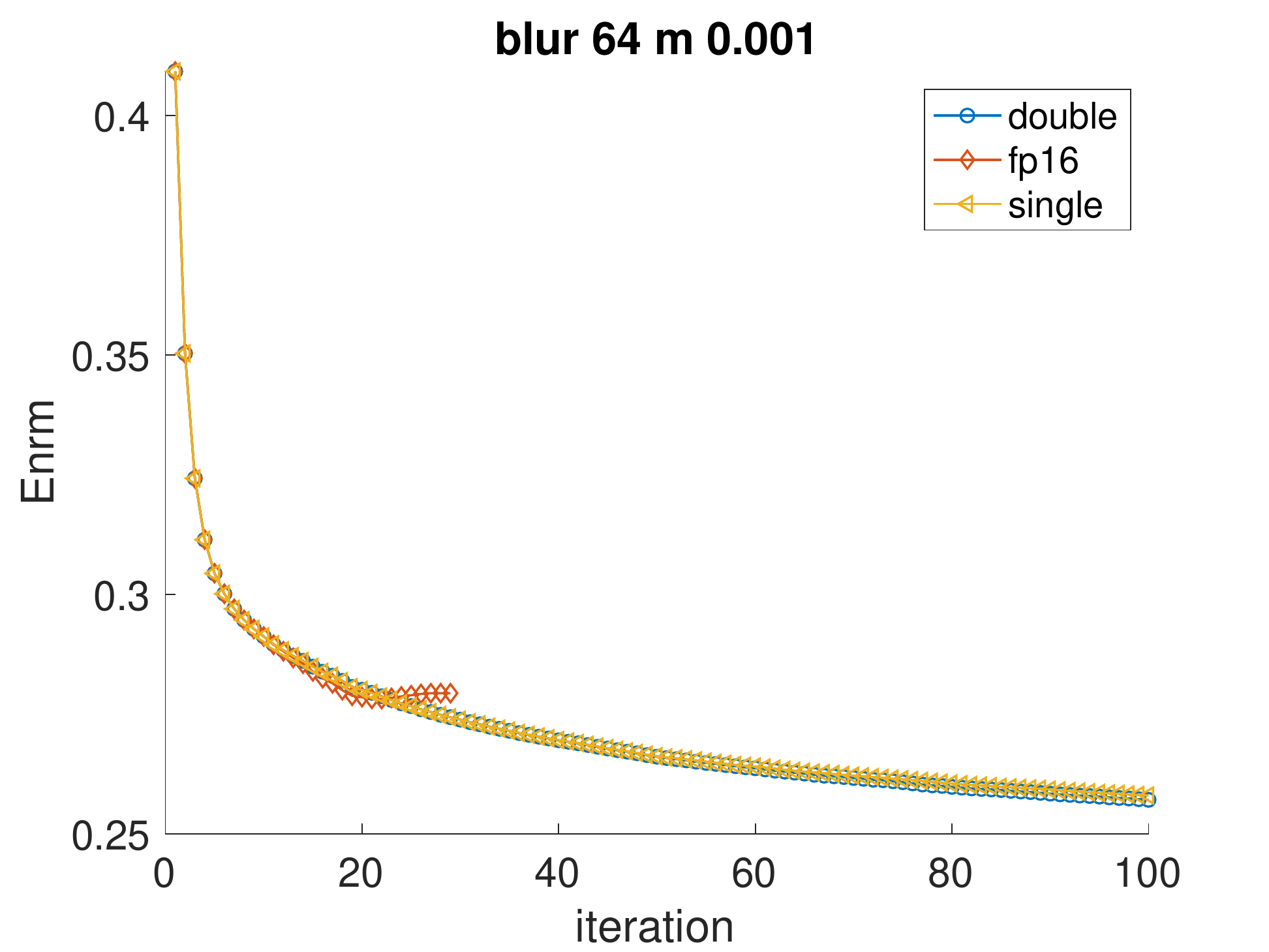}
  \caption{Error norm for problem size 64 with mild blurring and 0.1\% noise.}\label{fig:BEnrm_size64_m_0.1PNoise}
\endminipage\hfill
\minipage{0.32\textwidth}
  \includegraphics[width=\linewidth]{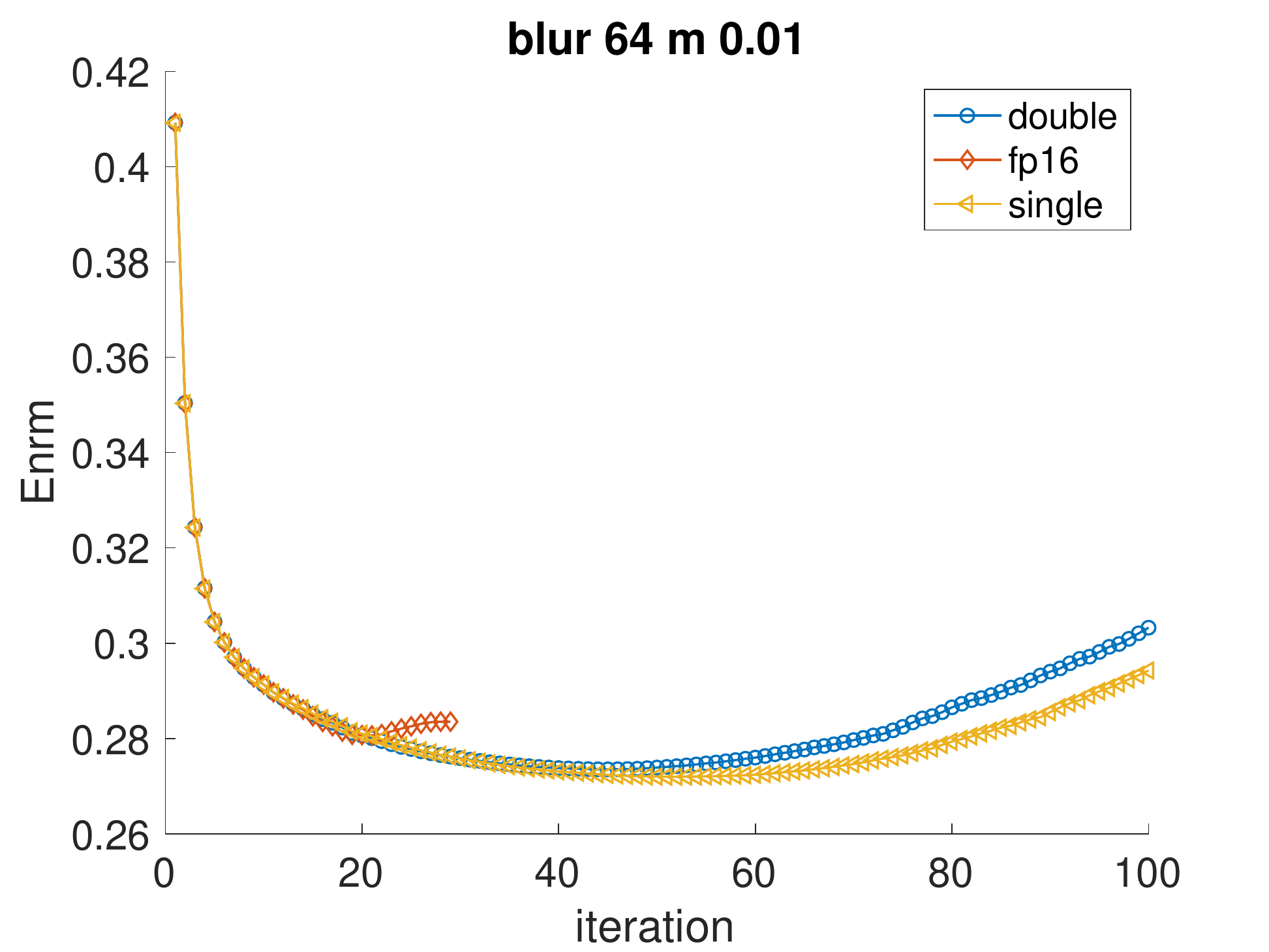}
  \caption{Error norm for problem size 64 with mild blurring and 1\% noise.}\label{fig:BEnrm_size64_m_1PNoise}
\endminipage\hfill
\minipage{0.32\textwidth}%
  \includegraphics[width=\linewidth]{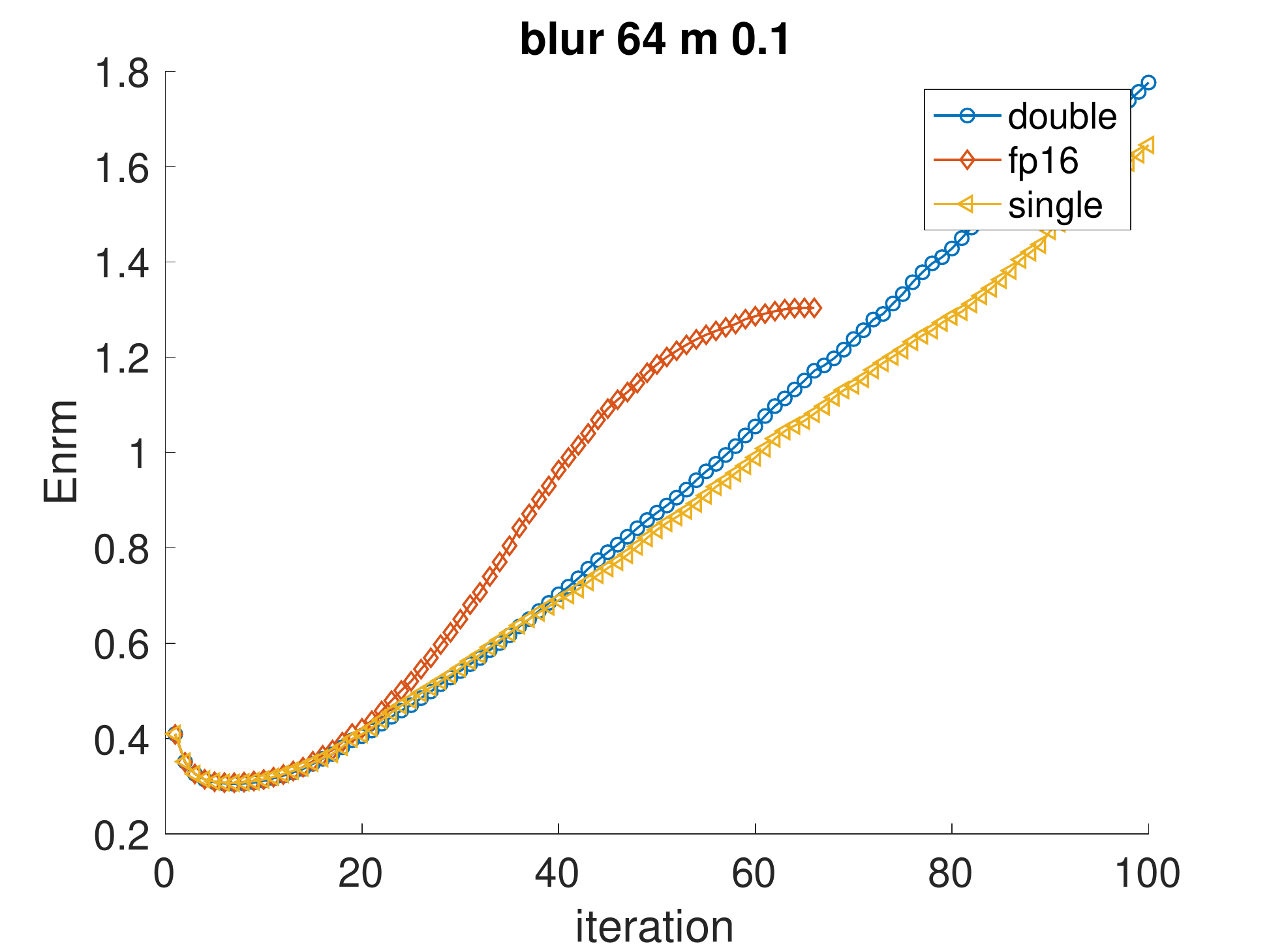}
  \caption{Error norm for problem size 64 with mild blurring and 10\% noise.}\label{fig:BEnrm_size64_m_10PNoise}
\endminipage
\end{figure}
\\An interesting phenomenon emerges if we examine the error norm for each iteration in Figures \ref{fig:BEnrm_size64_m_0.1PNoise}, \ref{fig:BEnrm_size64_m_1PNoise}, and \ref{fig:BEnrm_size64_m_10PNoise}. While the norm decreases as the iteration continues for the $0.1\%$ noise cases, those with $1\%$ and $10\%$ noise present a different trend: the error norm starts to increase again after certain iterations. The reason for this reversal is that the noise in $\bfb$ is inverted during each iteration along with the true $\bfb$; the noisy data eventually outweigh the actual solution, and the error norm increases. This phenomenon is known as semi-convergence in the inverse problems literature \cite{hansen2010discrete}.\\
\\So far all experimental results are acquired using our chopped version of the CGLS code without regularization. Section $4$ discusses another iterative method called Chebyshev semi-iterative method, which has the advantage of not requiring any inner products. However, the Chebyshev method does not have the same convergence behavior as CGLS, and it is therefore necessary to incorporate another regularization technique known as Tikhonov regularization. In order to make comparisons between these two methods, we later added Tikhonov regularization in our modified CGLS. The results and graphs will be displayed in Section $5$.

\subsubsection{Tomography Reconstruction}

\begin{figure}[!htb]
\minipage{0.32\textwidth}
  \includegraphics[width=\linewidth]{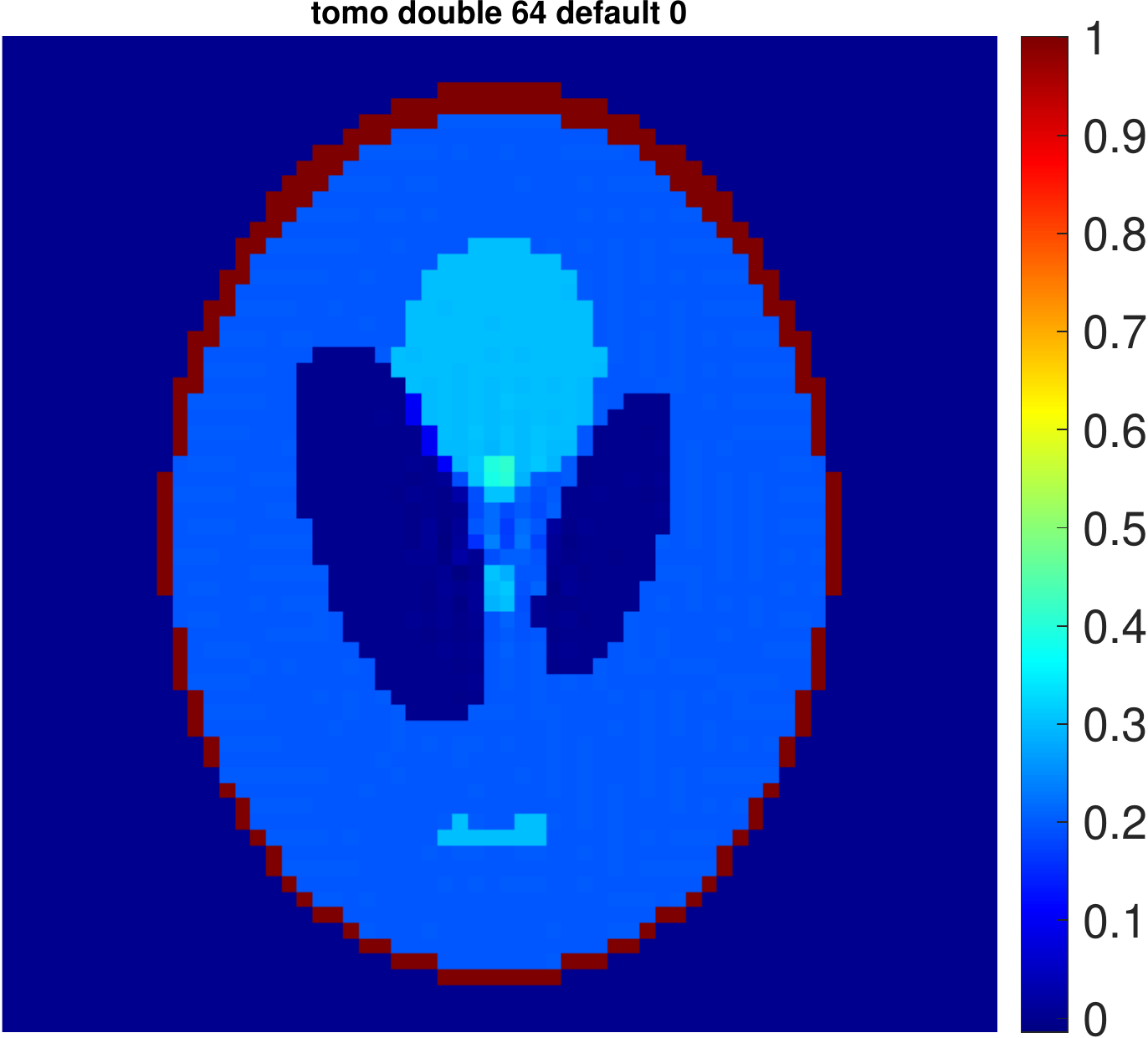}
  \caption{Double precision, size 64, zero noise.}\label{fig:Tdouble_size64_zeroNoise}
\endminipage\hfill
\minipage{0.32\textwidth}
  \includegraphics[width=\linewidth]{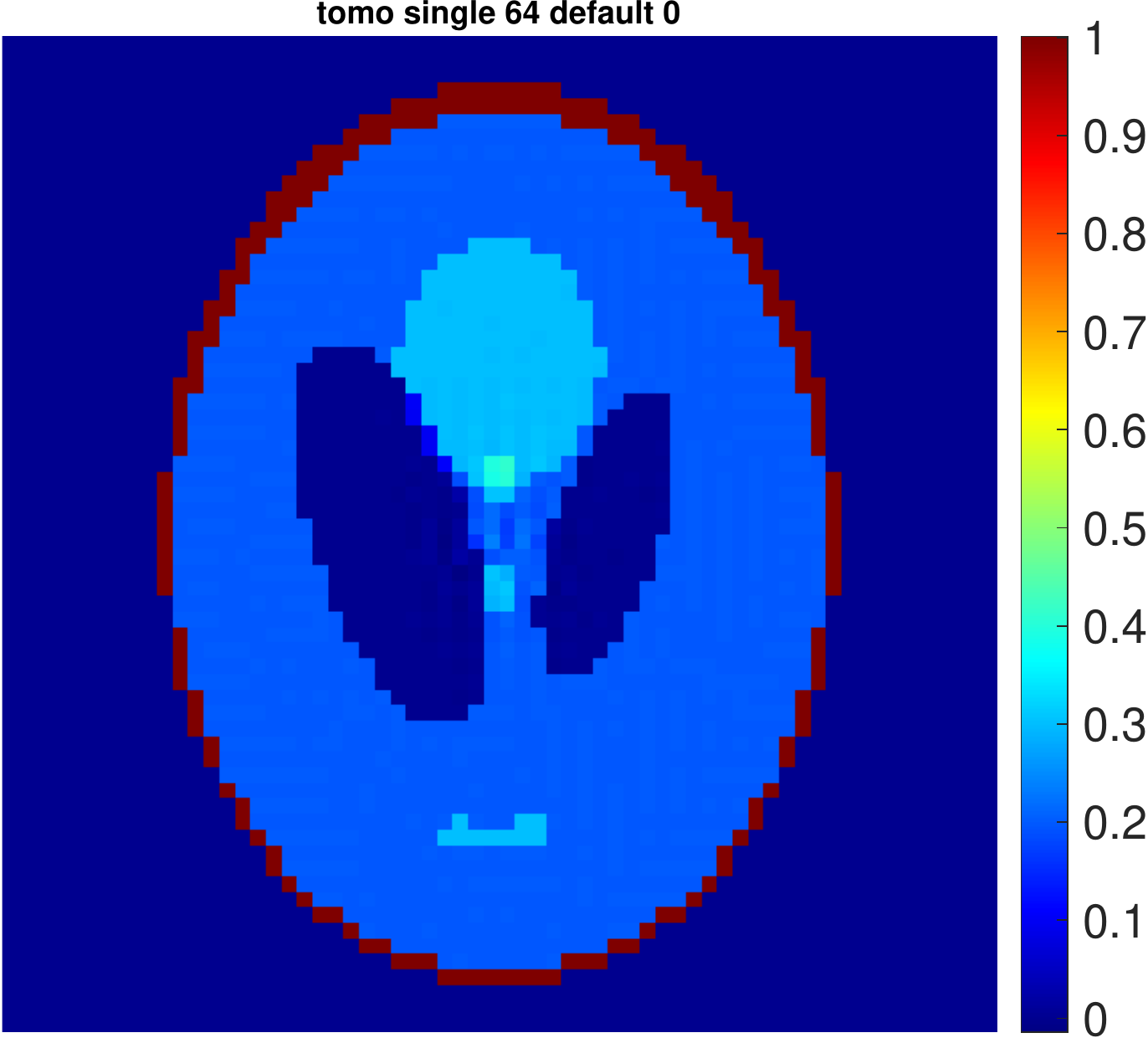}
  \caption{Single precision, size 64, zero noise.}\label{fig:Tsingle_size64_zeroNoise}
\endminipage\hfill
\minipage{0.32\textwidth}%
  \includegraphics[width=\linewidth]{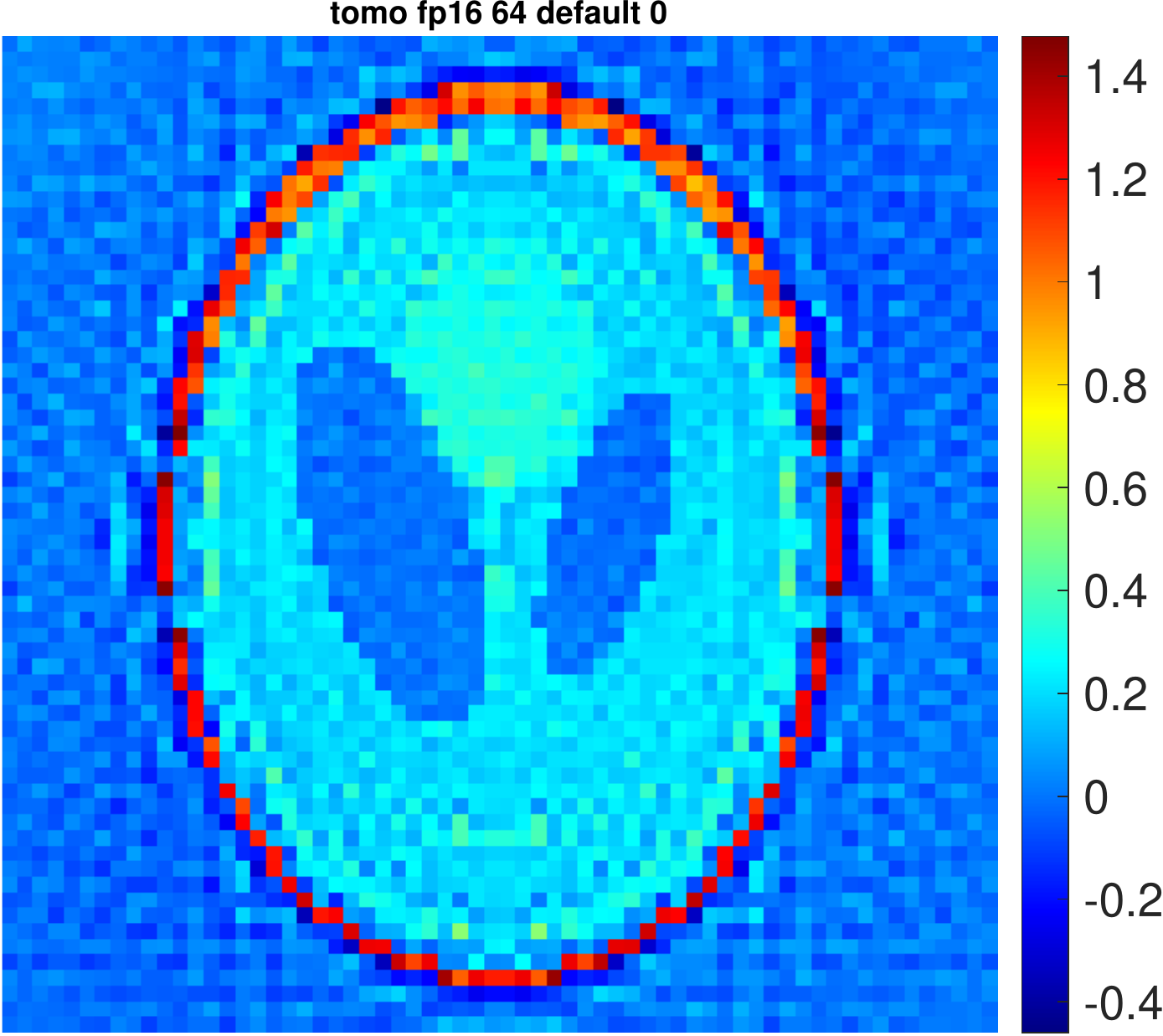}
  \caption{Half precision, size 64, zero noise (after rescaling).}\label{fig:Thalf_precision_size 64_zero noise}
\endminipage
\end{figure}

\begin{figure}[!htb]
\minipage{0.32\textwidth}
  \includegraphics[width=\linewidth]{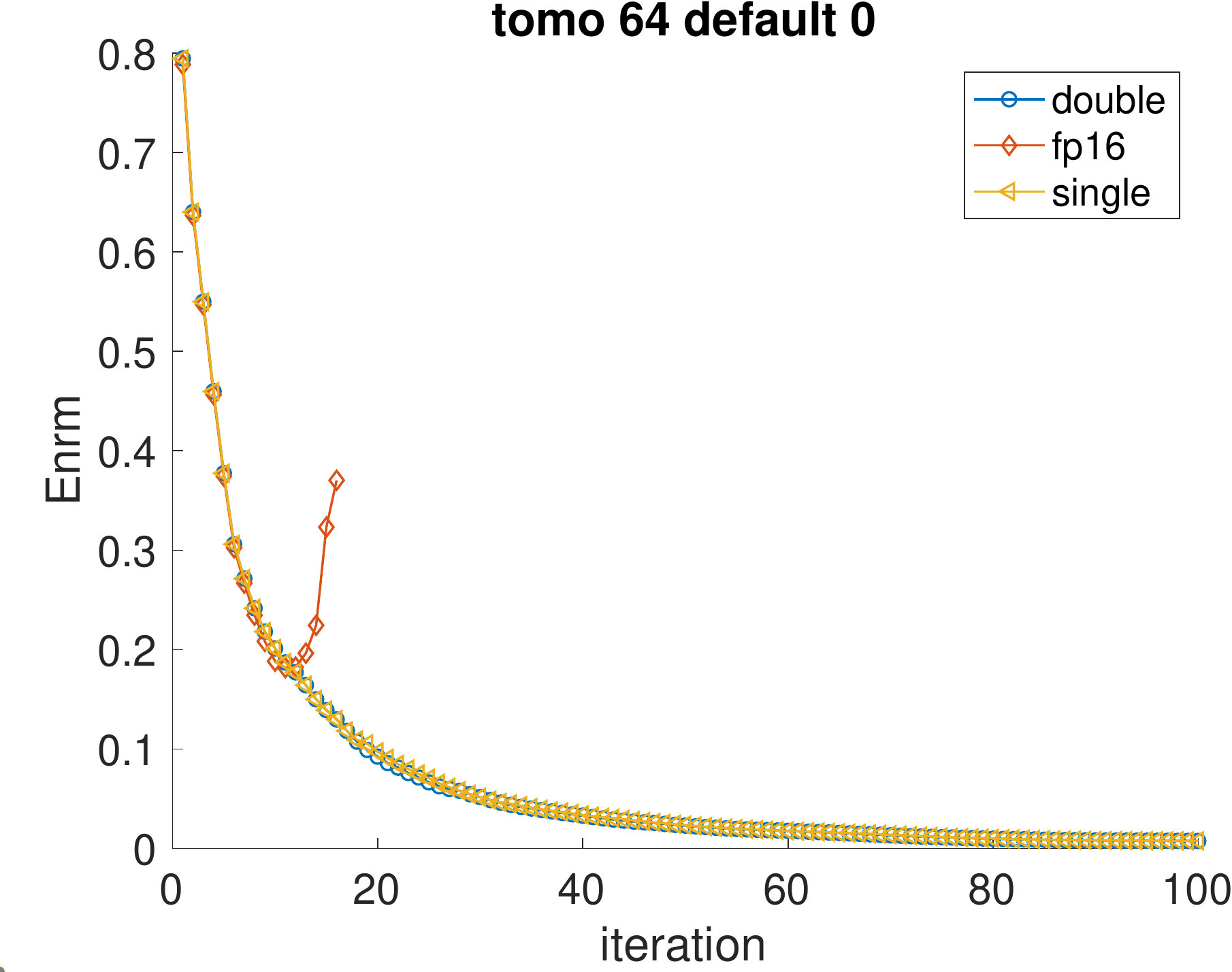}
  \caption{Error norm for size 64 problem with zero noise.}\label{fig:TEnrm_size64_zeroNoise}
\endminipage\hfill
\minipage{0.32\textwidth}
  \includegraphics[width=\linewidth]{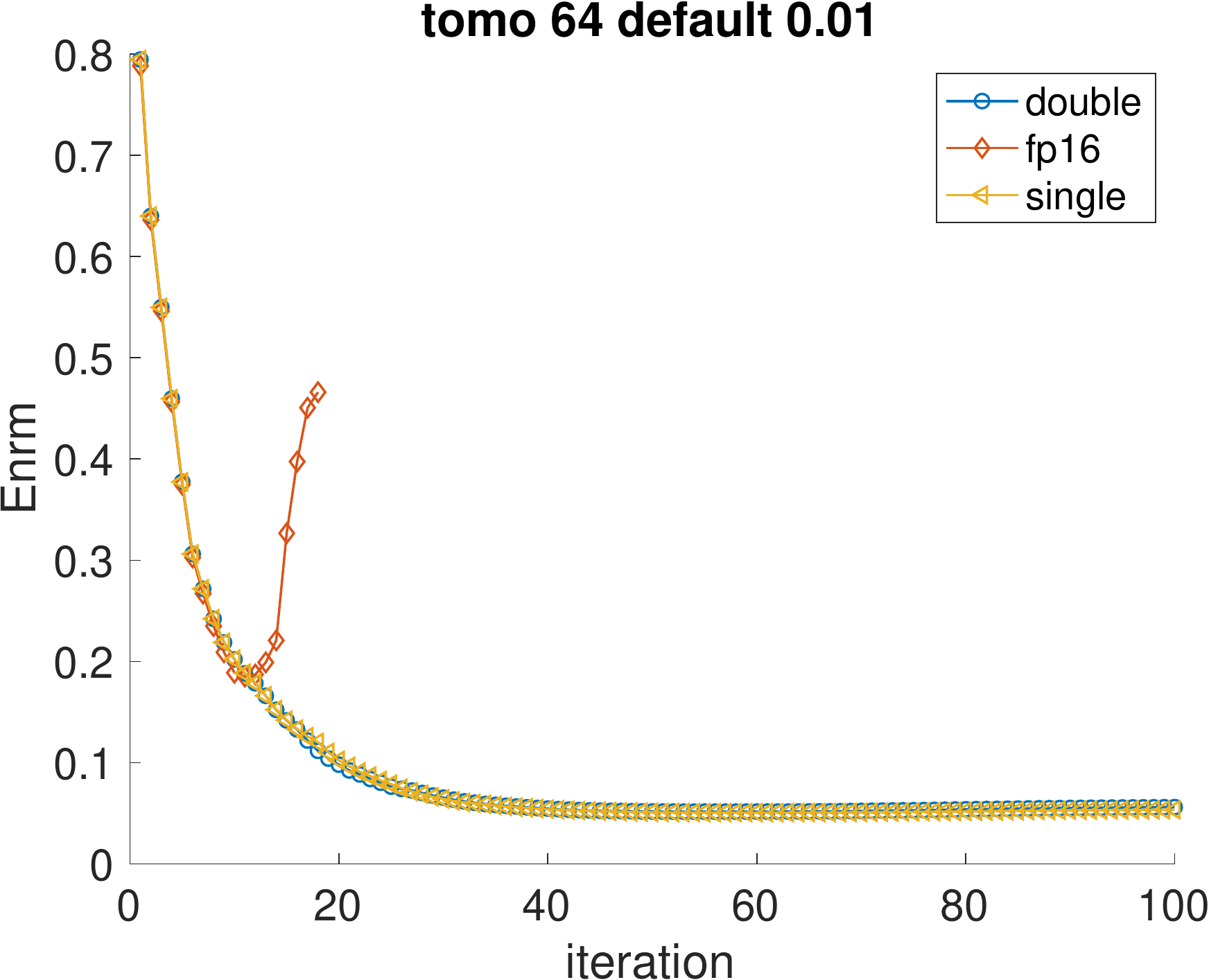}
  \caption{Error norm for size 64 problem with 1\% noise.}\label{fig:TEnrm_size64_1PNoise}
\endminipage\hfill
\minipage{0.32\textwidth}%
  \includegraphics[width=\linewidth]{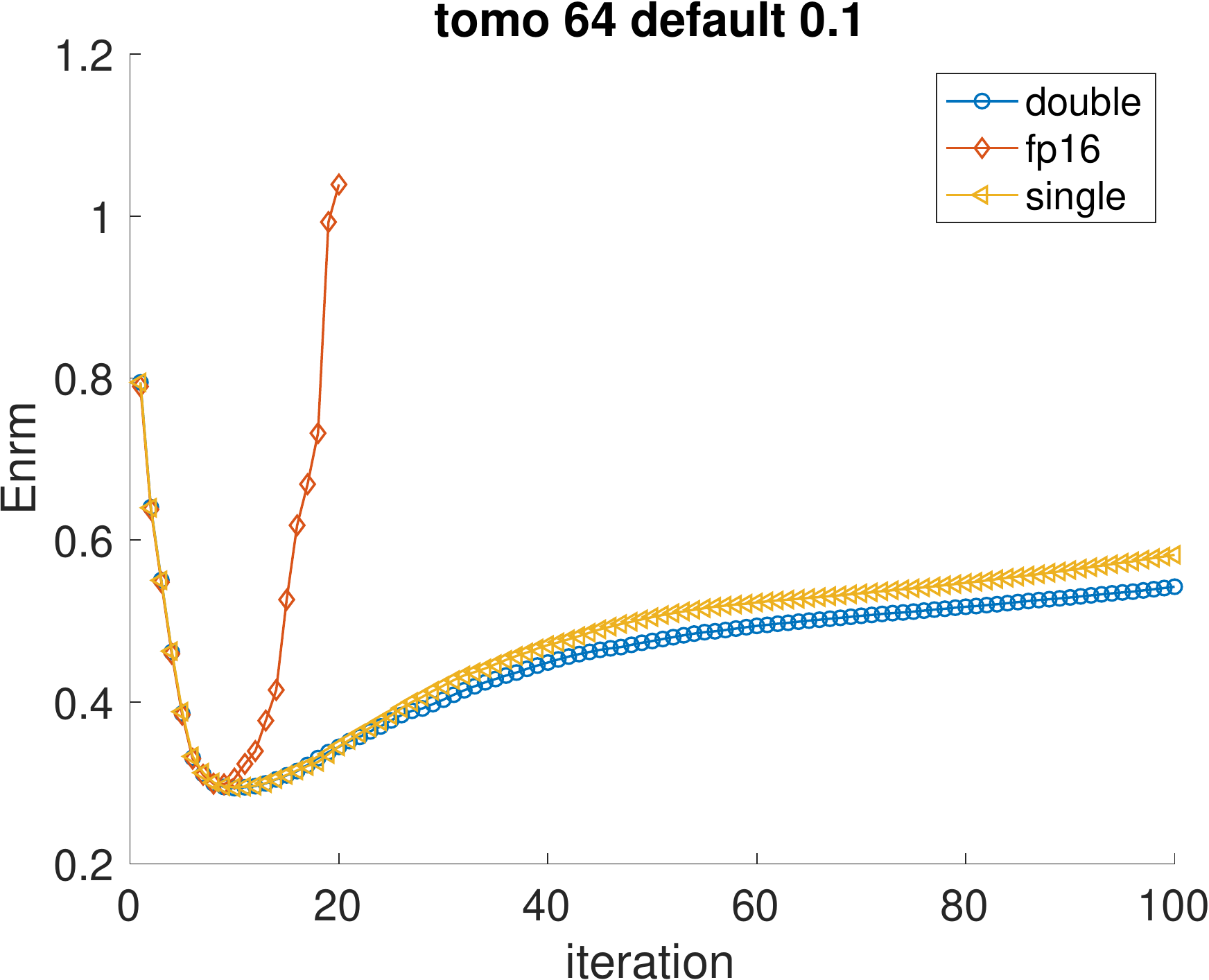}
  \caption{Error norm for size 64 problem with 10\% noise.}\label{fig:TEnrm_size64_10PNoise}
\endminipage
\end{figure}

Tomography reconstruction is another type of inverse problem that produces images from X-ray projection data of various angles. We used it as a test problem since the IRtools software package includes a simulation in the \textbf{PRtomo} function. The modified CGLS method works well for double and single precision and generates nice reconstructions as shown in Figures \ref{fig:Tdouble_size64_zeroNoise} and \ref{fig:Tsingle_size64_zeroNoise}, but issues arise for half precision because entries of the solution $\bfx$ are all NaNs. Since half precision formats allocate only $5$ bits for the exponent, overflow easily occurs when the inner product is calculated. NaNs are the results from dividing infinity by infinity in the CGLS method.\\
\\One of the solutions to the overflow is to rescale the matrices. After scaling both $A$ and $\bfb$ by $0.01$, we obtain Figure \ref{fig:Thalf_precision_size 64_zero noise}. Although the image is still not as clear as those from double and single precision, the inner products did not overflow and a meaningful picture was obtained.\\
\\
\noindent The noise level affects the quality of the reconstruction as well. As Figures \ref{fig:Tdouble precision, size 64, zero noise}, \ref{fig:Tsingle precision, size 64, zero noise}, \ref{fig:Thalf precision, size 64, zero noise} illustrate, while random noise does not completely destroy the reconstruction and take over the original information as is the case for image deblurring, the artifacts caused by the noise make it difficult to see small objects in the image.
\begin{figure}[!htb]
\minipage{0.32\textwidth}
  \includegraphics[width=\linewidth]{tomo_fp16_64_default_0}
  \caption{Half precision, size 64, zero noise.}\label{fig:Tdouble precision, size 64, zero noise}
\endminipage\hfill
\minipage{0.32\textwidth}
  \includegraphics[width=\linewidth]{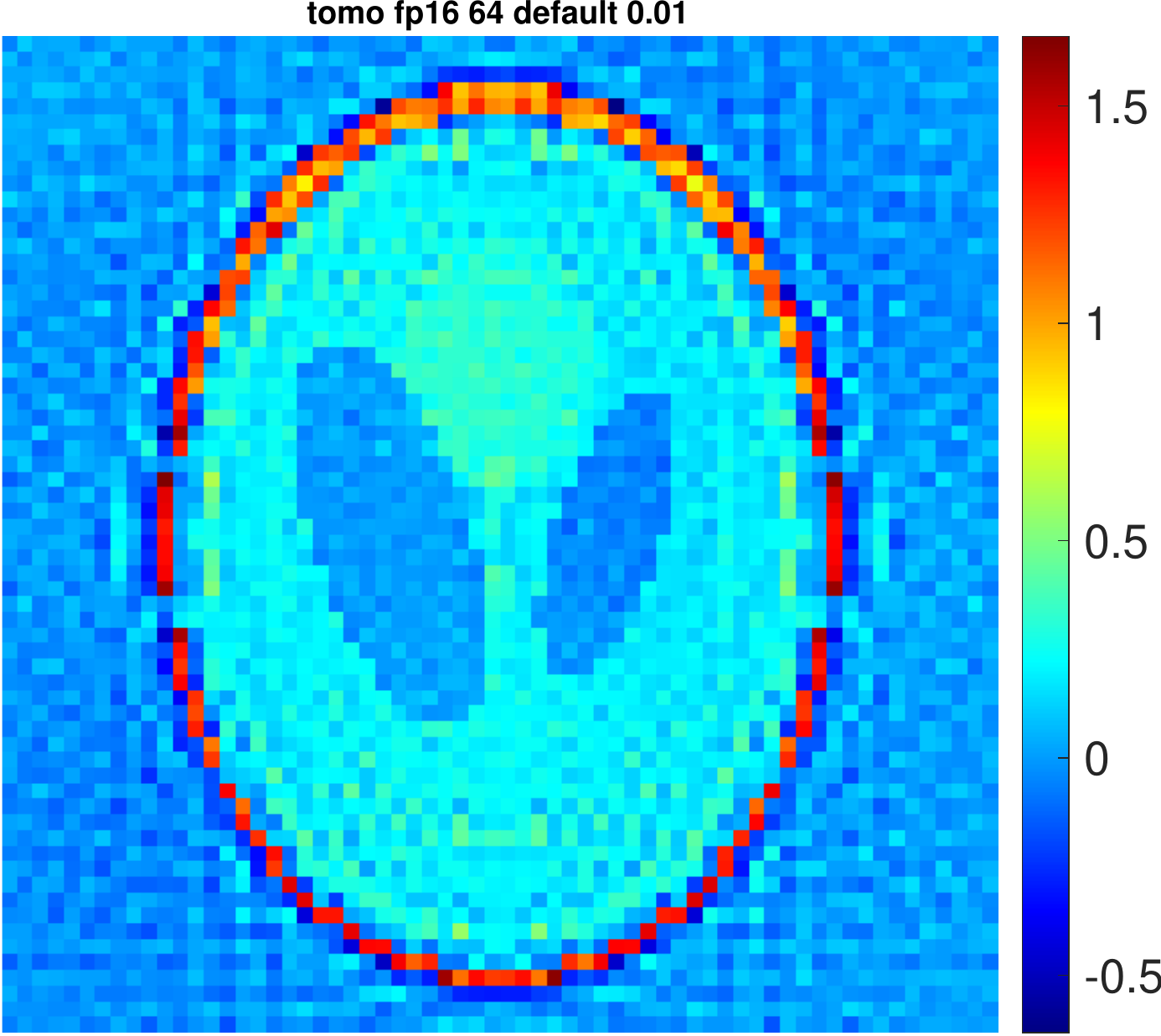}
  \caption{Half precision, size 64, 1\% noise.}\label{fig:Tsingle precision, size 64, zero noise}
\endminipage\hfill
\minipage{0.32\textwidth}%
  \includegraphics[width=\linewidth]{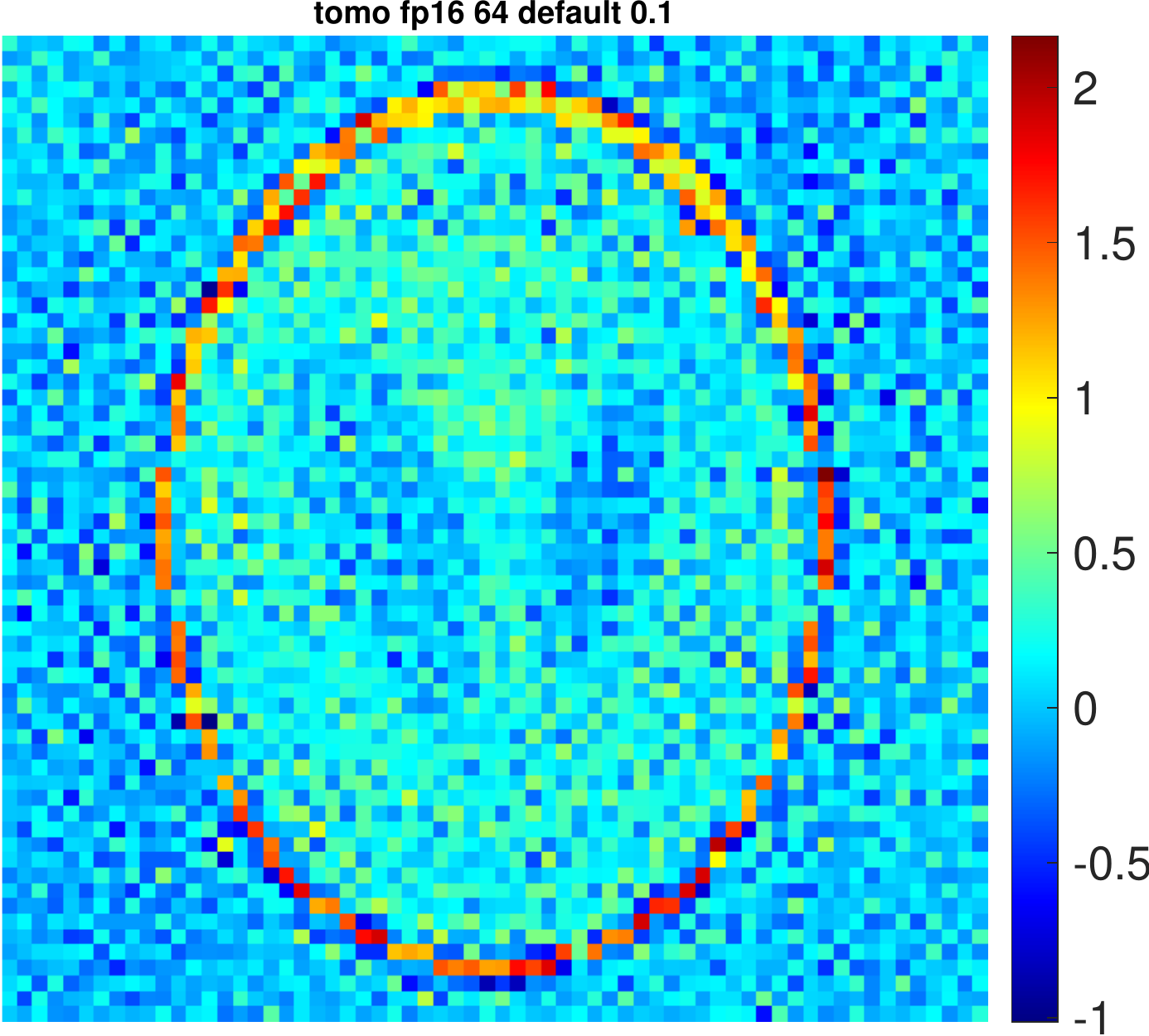}
  \caption{Half precision, size 64, 10\% noise (after rescaling).}\label{fig:Thalf precision, size 64, zero noise}
\endminipage
\end{figure}
\\
\noindent As was the case with the deblurring problem, the error norm begins to increase at some point of the iteration due to the accumulation of the inverted noise, as shown in Figures \ref{fig:TEnrm_size64_zeroNoise}, \ref{fig:TEnrm_size64_1PNoise}, and \ref{fig:TEnrm_size64_10PNoise}. At half precision, the norm grows after the $11^{th}$ iteration even without noise in Figure \ref{fig:TEnrm_size64_zeroNoise}. The reason this time is not the attached noise but the truncation errors.

\section{Chebyshev Semi-Iterative Method}

\subsection{Method Overview}
\label{sec:chebyshev}
The Chebyshev Semi-Iterative (CS) Method requires no inner product computation, so it is safer in terms of overflow issues when at low precision. The trade-off here is that the CS method needs prior knowledge of the range of singular values of the matrix $A$~\cite{wang2015chebyshev}.\\
As with CG, the CS method is generally derived for symmetric positive definite linear systems, but there are variations for least squares problems, which are applied to the normal equations $A^{T}A \bfx = A^{-1}\bfb$.
\noindent The following algorithm\footnote{We are using Algorithm 3 from \cite{meng2014lsrn}, but we remark that there is a typographical error for the $\alpha$ parameter for the case $k=1$. We show the correct formula in our paper.} describes the implementation of the CS method for least squares problems that we used in our work ~\cite{meng2014lsrn}~\cite{gutknecht2002chebyshev}:

\begin{algorithm}
\caption{Chebyshev semi-iterative method}\label{alg:cap}
\begin{algorithmic}
\STATE Given $A \in \mathbb{R}^{m \times n}$, $\bfb \in \mathbb{R}^m$, and a tolerance $\epsilon > 0$, choose $0 < \sigma_L < \sigma_U$ such that all nonzero singular values of A in [$\sigma_L$, $\sigma_U$], and let $d = \frac{{\sigma_U}^2 + {\sigma_L}^2}{2}$ and $c = \frac{{\sigma_U}^2 - {\sigma_L}^2}{2}$ 
\STATE Let $\bfx = 0$, $\bfv = 0$, and $\bfr = \bfb$.
\FOR {$k = 0,1,\ldots, \lceil{\log \epsilon - \log 2}/{\log \frac{\sigma_U - \sigma_L}{\sigma_U + \sigma_L}}\rceil$} 
\STATE $\beta \leftarrow
\begin{cases}
    0, & \text{if } k = 0\\
    \frac{1}{2}(c/d)^2, & \text{if } k = 1\\
    (\alpha c/2)^2, & \text{otherwise, } 
\end{cases}
       \alpha \leftarrow
\begin{cases}
    1/d, & \text{if } k = 0\\
    1/(d - c^2/(2d)), & \text{if } k = 1 \\ 
    1/(d-\alpha c^2/4), & \text{otherwise }
\end{cases} $
    
\STATE $\bfv \gets \beta \bfv + A^T\bfr$.
\STATE $\bfx \gets \bfx + \alpha \bfv$.
\STATE $\bfr \gets \bfr - \alpha A\bfv$.
\ENDFOR
\end{algorithmic}
\end{algorithm}

\subsection{Experiment}
\begin{figure}[!htb]
\minipage{0.32\textwidth}
  \includegraphics[width=\linewidth]{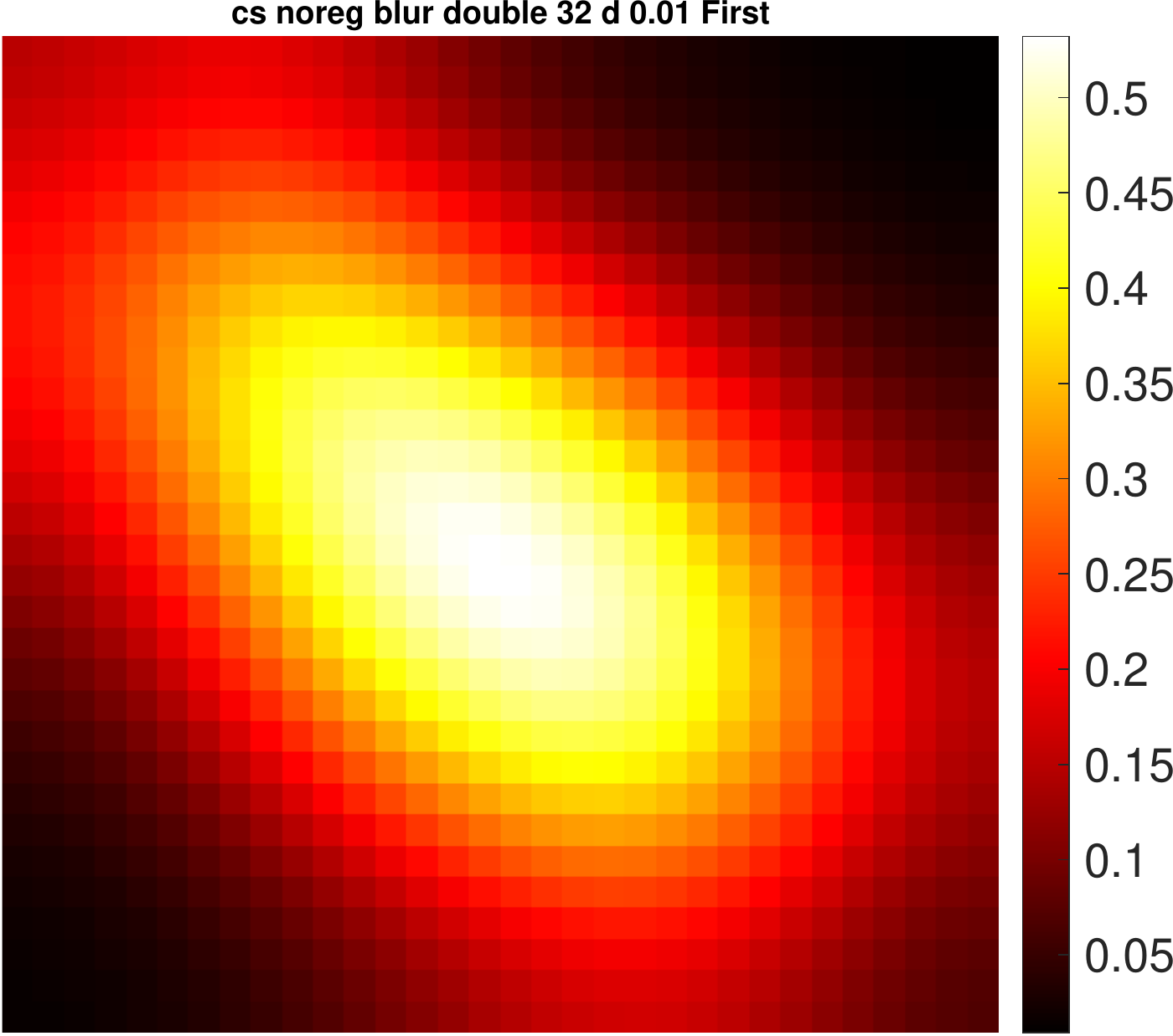}
  \caption{Reconstruction at first iteration.}\label{fig:cs_noreg_double_32_1PNoise_it1}
\endminipage\hfill
\minipage{0.32\textwidth}
  \includegraphics[width=\linewidth]{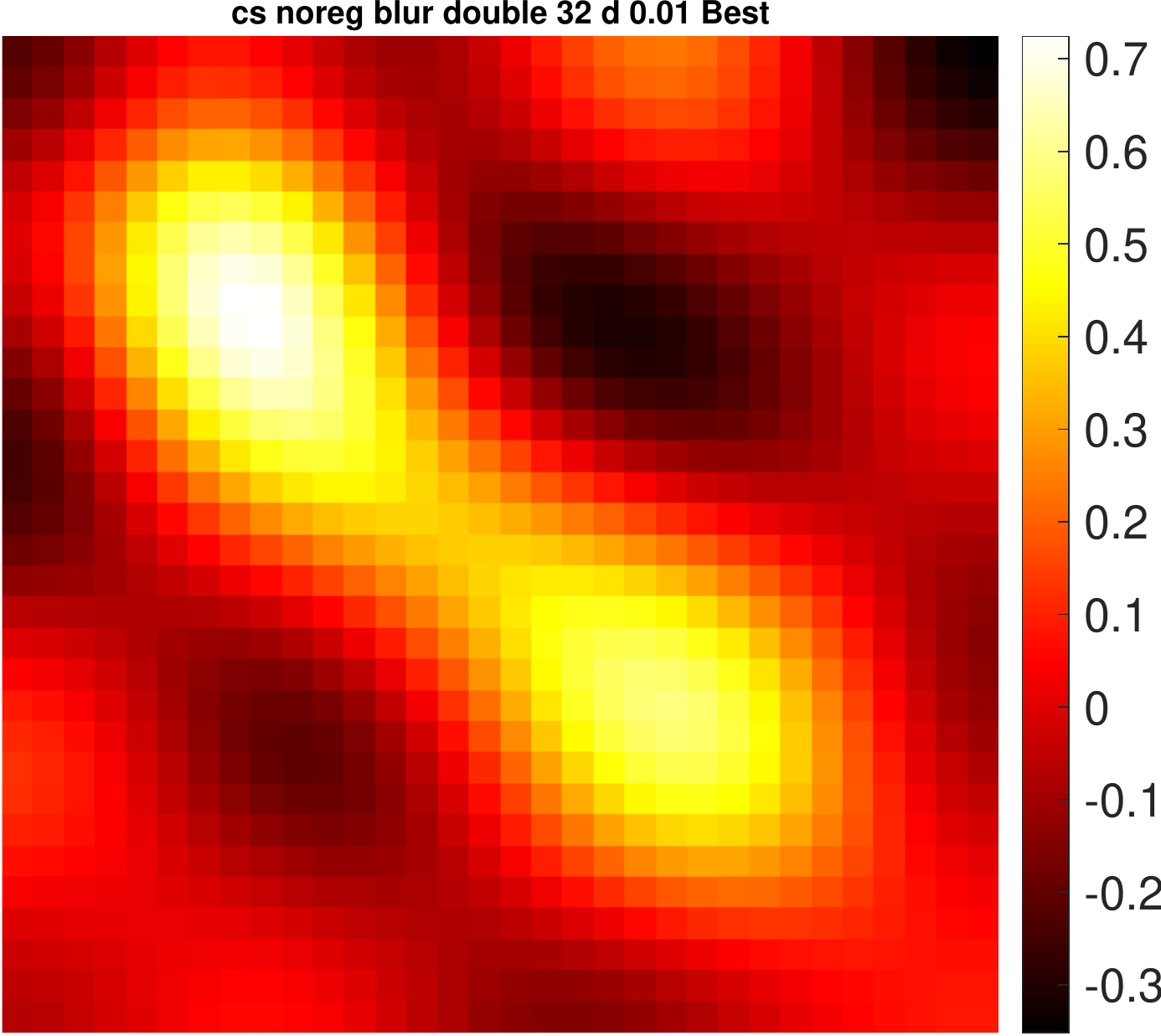}
  \caption{Reconstruction at the best iteration.}\label{fig:cs_noreg_double_32_1PNoise_it5}
\endminipage\hfill
\minipage{0.32\textwidth}%
  \includegraphics[width=\linewidth]{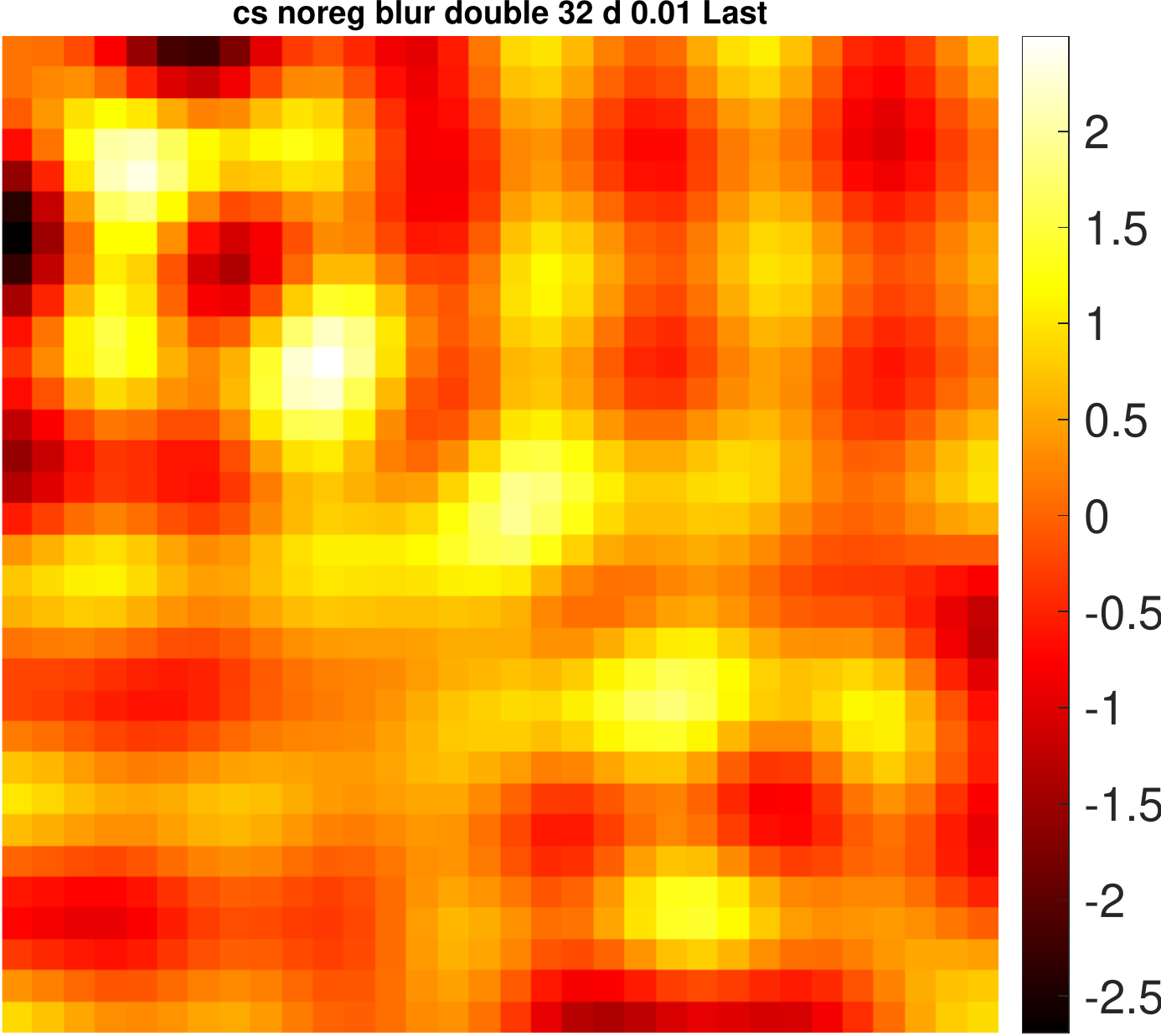}
  \caption{Reconstruction at the last iteration.}\label{fig:cs_noreg_double_32_1PNoise_it500}
\endminipage
\end{figure}

Again, we used image deblurring and tomography reconstruction as two test problems for the CS method in low precision. At first we ran the image deblurring test problems of size 32 with default blur and $1\%$ noise in double precision. We got estimations for the bound of singular values of matrix $A$ with the built-in MATLAB function \textbf{SVDS}. We set the maximum number of iterations to be $500$, and plotted the graph at the first, $5^{th}$ (where the error norm is smallest), and last iteration, as shown in Figures \ref{fig:cs_noreg_double_32_1PNoise_it1}, \ref{fig:cs_noreg_double_32_1PNoise_it5}, and \ref{fig:cs_noreg_double_32_1PNoise_it500}. \\

\noindent The reconstruction is doing poorly even at the best iteration (Figure \ref{fig:cs_noreg_double_32_1PNoise_it5}). At the last iteration, it is even worse due to over-fitting. After a closer look, we noticed that the number of iterations recommended by the algorithm is at order of magnitude of $10^{12}$. If we increase the size of the problem to $64\times 64$, the number overflows to -Inf, which results from the tiny estimation of the lower bound of matrix $A$'s singular values. The estimation is so close to zero that we suspect this is actually a result of round-off errors on a zero entry. Therefore, we implemented Tikhonov regularization for two reasons: a) to avoid over-fitting; and b) to obtain a more valid lower bound of singular values by increasing them to a larger value. Specifically, the singular values for the matrix $A$ after Tikhonov regularization with regularization parameter $\lambda$ would be $\sqrt{\sigma_{i}^2+\lambda^2}$. The technique is discussed in more details in the next section \ref{sec:tikhonov}.

\subsubsection{Tikhonov Regularization}
\label{sec:tikhonov}

Tikhonov Regularization includes a regularization term to the original least squares problem:\\
$$\min_{x} \{||A\bfx-\bfb||_2^2+\lambda^2||\bfx||_2^2\},$$
where the regularization parameter $\lambda$ balances the residual term $||A\bfx-\bfb||_2^2$ and the regularization term $||\bfx||_2^2$. 
We can rewrite the Tikhonov problem as a least squares problem
$$\min_{\bfx} \begin{Vmatrix}\begin{pmatrix}
A\\
\lambda I
\end{pmatrix}\bfx-
\begin{pmatrix}
\bfb\\
0
\end{pmatrix}\end{Vmatrix}_2,$$
and the solution to this least squares problem is
$$\bfx_{\lambda}=(A^{T}A+\lambda^2I)^{-1}A^{T}\bfb.$$
To see why Tikhonov regularization is effective, observe that if we substitute the singular value decomposition, $\displaystyle A = U\Sigma V^{T}$, into the expression for $\bfx_{\lambda}$, and expand the matrix multiplication column-wise, we would get
$$\bfx_{\lambda}=\sum_{i=1}^{n}\phi_i^{[\lambda]}\frac{\bfu_i^{T}\bfb}{\sigma_i}\bfv_i,$$
where the filter factors $\phi_i^{[\lambda]}$ are
$$\phi_i^{[\lambda]} = \frac{\sigma_i^{2}}{\sigma_i^{2}+\lambda^{2}}.$$
Notice that the filter factor $\phi_i^{[\lambda]}$ is approximately equal to $0$ for smaller singular values and approximately equal to $1$ for larger singular values. It therefore acts like a filter by decreasing the effects of magnifying noise in $\bfb$ when divided by tiny singular values.\\
With regularization, we are running the CS method with the matrix $\begin{pmatrix}
A\\
\lambda I
\end{pmatrix}$ instead of $A$.\\
Given $A = U\Sigma V^{T}$, we have 
$$A^{T}A = V\Sigma^T \Sigma V^{T}$$
and\\
$$\begin{pmatrix}
A^{T} & \lambda I\\
\end{pmatrix}
\begin{pmatrix}
A\\
\lambda I
\end{pmatrix} = A^{T}A+\lambda^2I = V\Sigma^T \Sigma V^{T}+V\lambda^2 V^{T}=V(\Sigma^T \Sigma+\lambda^2) V^{T}.$$
Therefore, we know the singular values of $\begin{pmatrix}
A\\
\lambda I
\end{pmatrix}$ are $\sqrt{\sigma_{i}^2+\lambda^2}$, where $\sigma_{i}$ are singular values of $A$. When $\sigma_{i}$ is small, $\sqrt{\sigma_{i}^2+\lambda^2} \approx \lambda$, and we can directly use $\lambda$ as the lower bound of singular values for the matrix after regularization. With Tikhonov regularization, the singular values are increased to a more reasonable lower bound.

\subsubsection{Image Deblurring}

\begin{figure}[!htb]
\minipage{0.32\textwidth}
  \includegraphics[width=\linewidth]{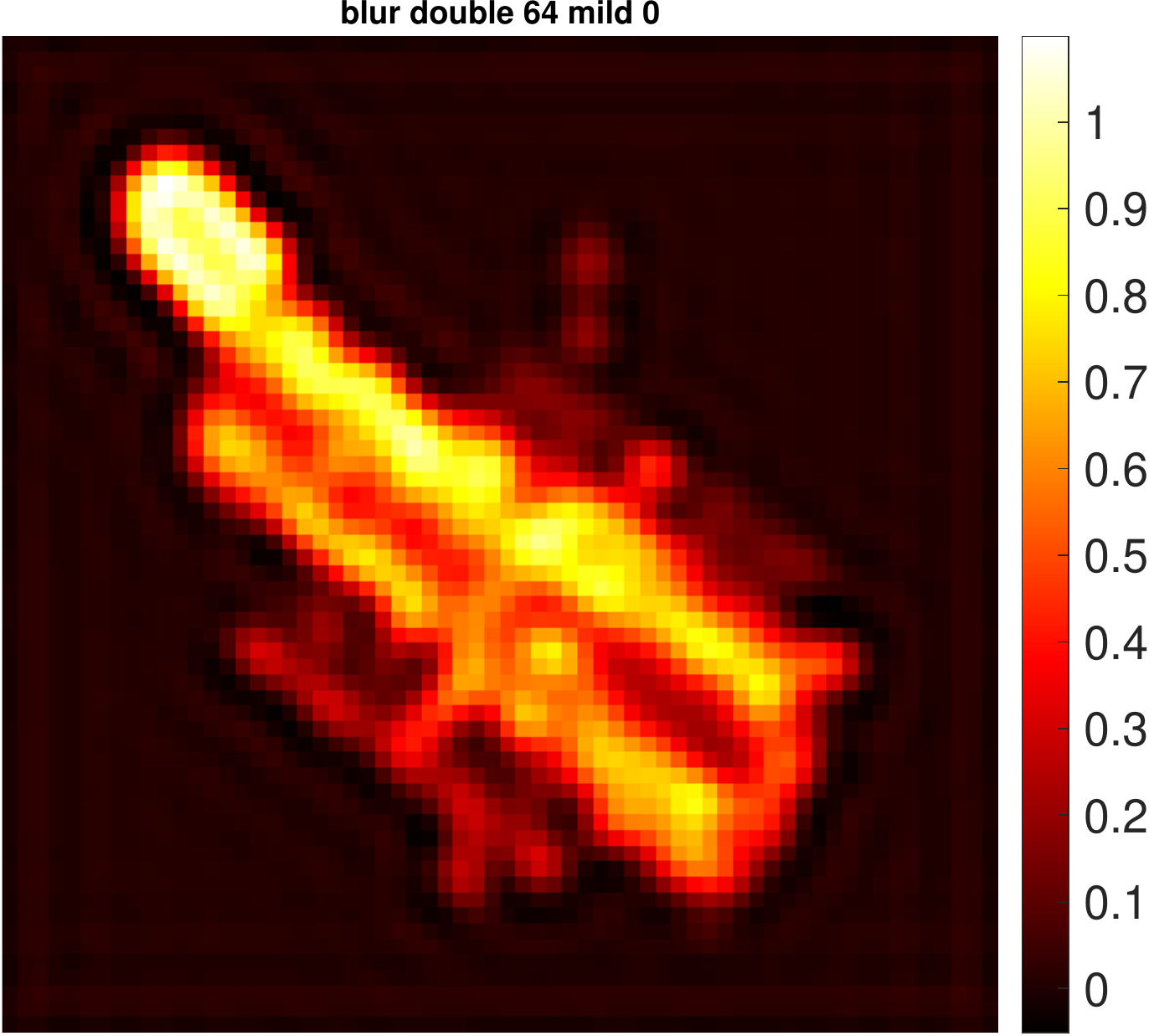}
  \caption{Image reconstruction in double precision of size 64 problem with mild blurring and no noise.}\label{fig:Bcs_m_double_64_0noise}
\endminipage\hfill
\minipage{0.32\textwidth}
  \includegraphics[width=\linewidth]{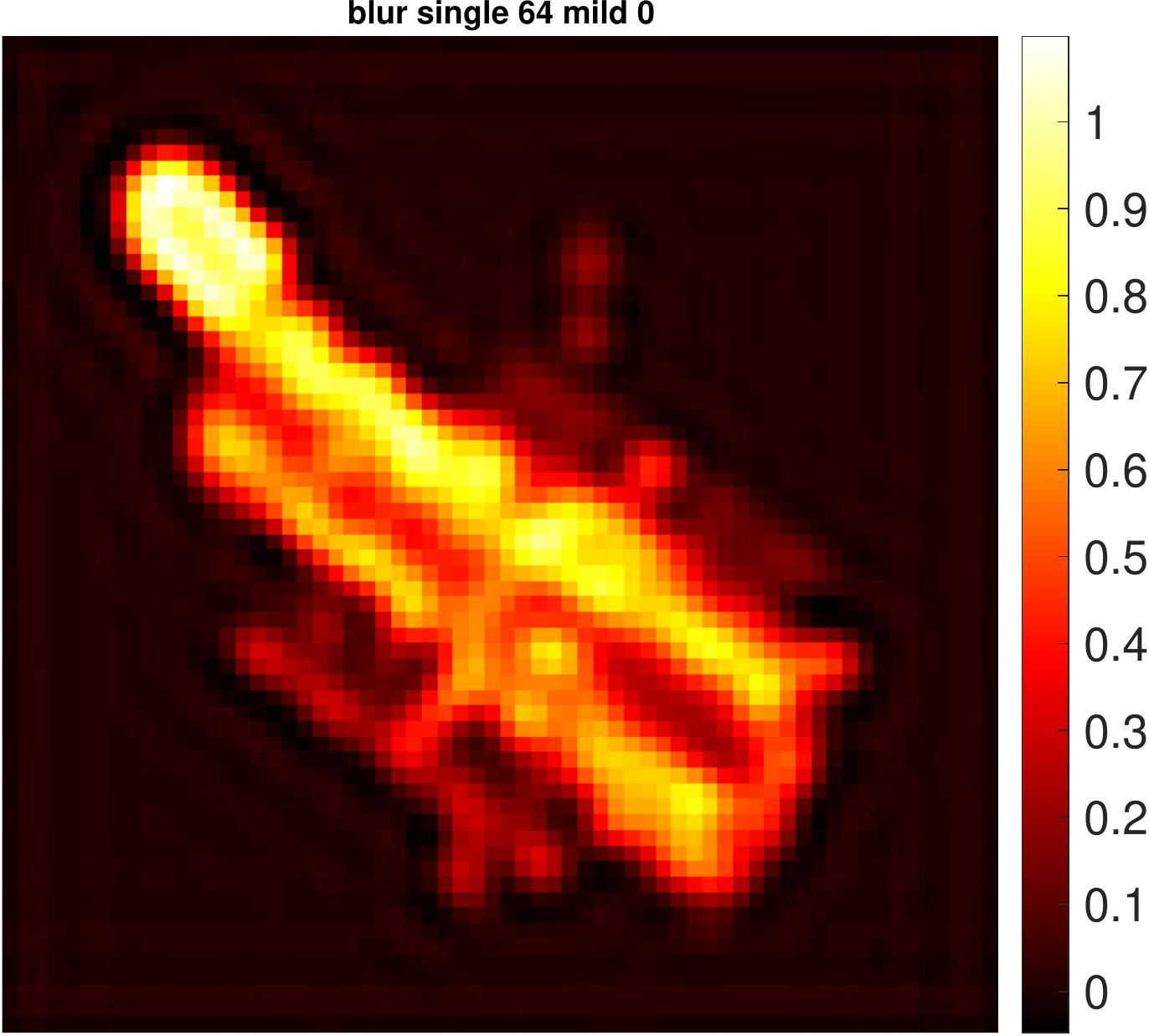}
  \caption{Image reconstruction in single precision of size 64 problem with mild blurring and no noise.}\label{fig:Bcs_m_single_64_0noise}
\endminipage\hfill
\minipage{0.32\textwidth}%
  \includegraphics[width=\linewidth]{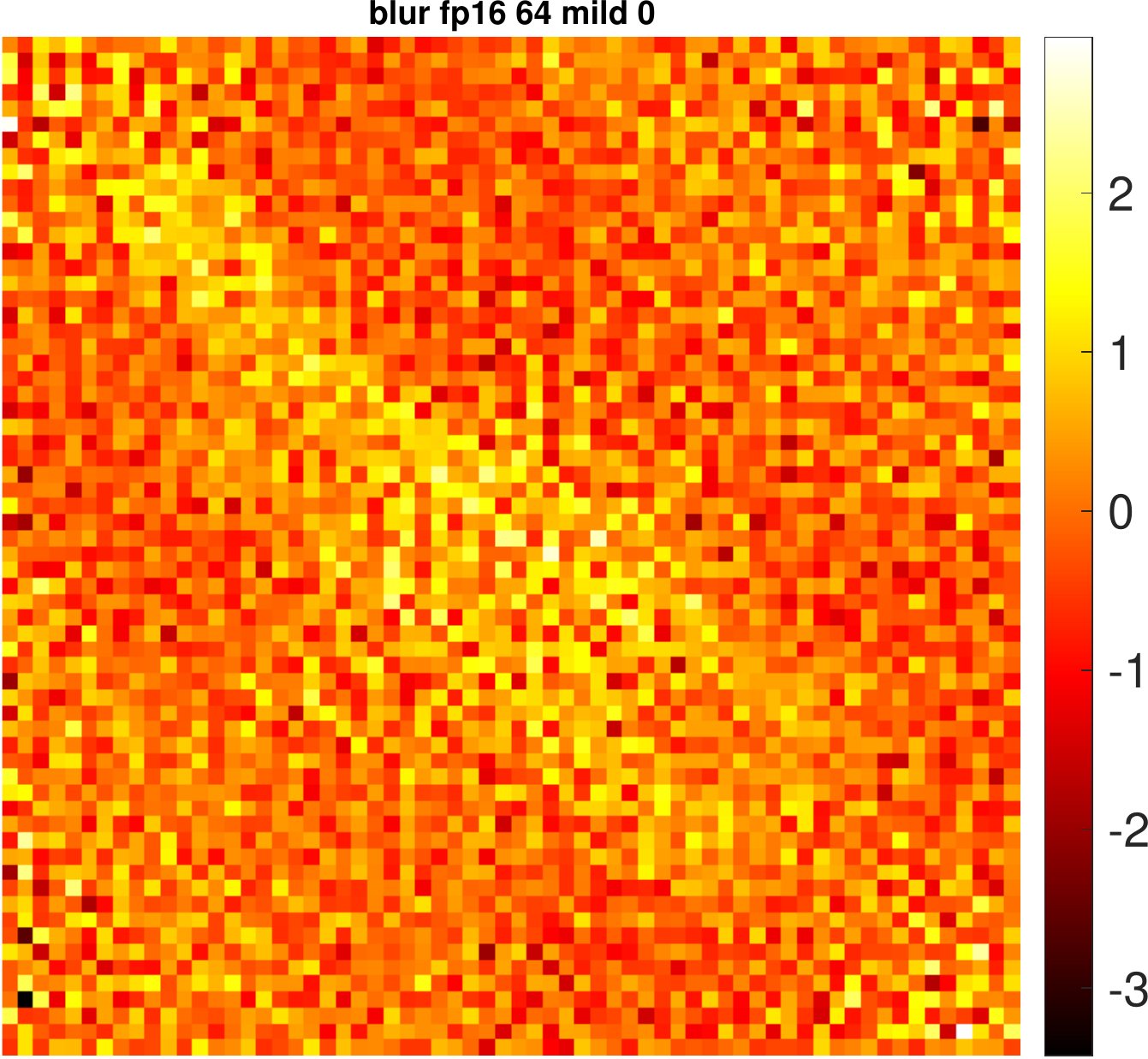}
  \caption{Image reconstruction in half precision of size 64 problem with mild blurring and no noise.}\label{fig:Bcs_m_half_64_0noise}
\endminipage
\end{figure}

\begin{figure}[!htb]
\minipage{0.32\textwidth}
  \includegraphics[width=\linewidth]{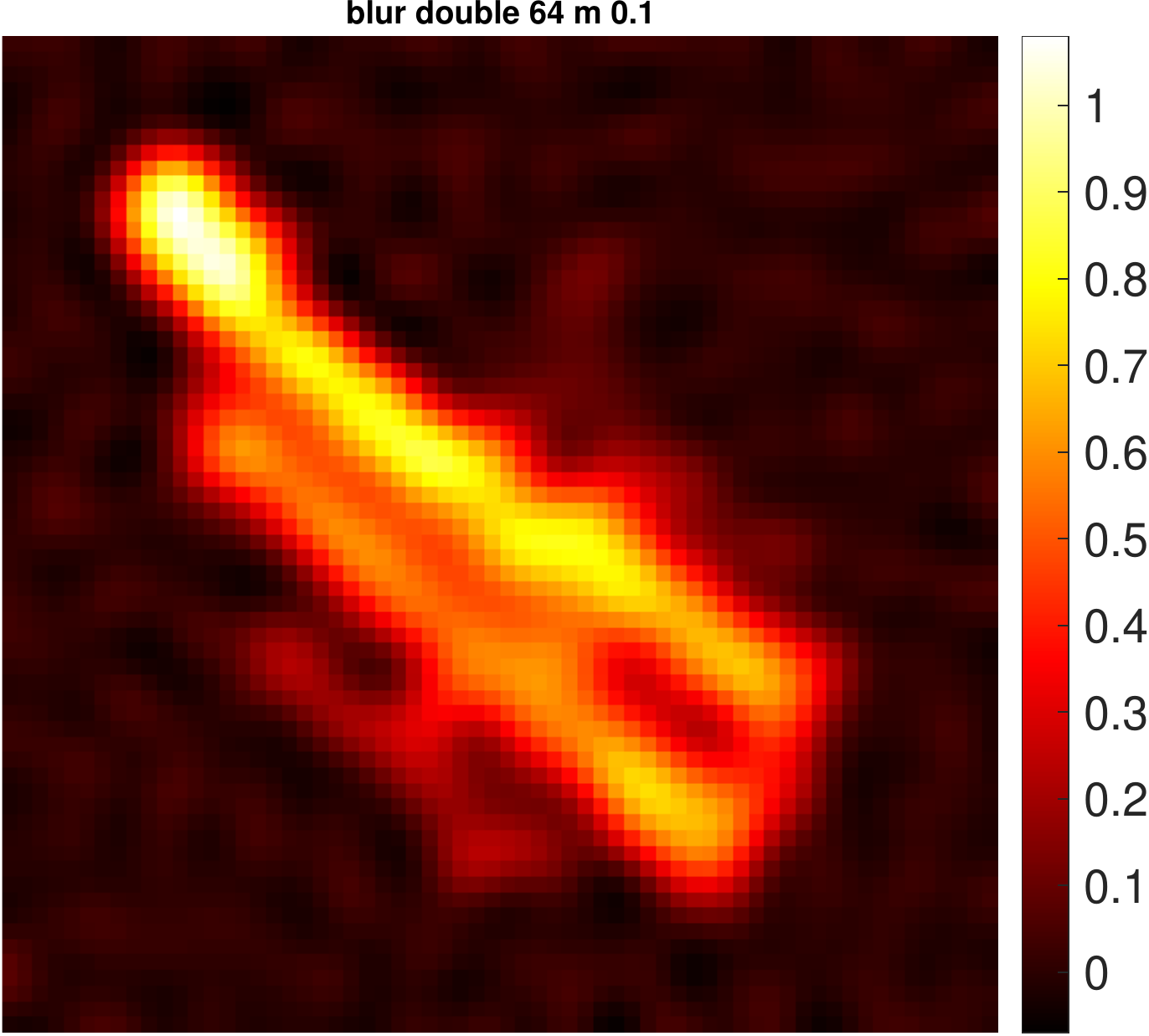}
  \caption{Image reconstruction in double precision of size 64 problem with mild blurring and 10\% noise.}\label{fig:Bcs_m_double_64_10PNoise}
\endminipage\hfill
\minipage{0.32\textwidth}
  \includegraphics[width=\linewidth]{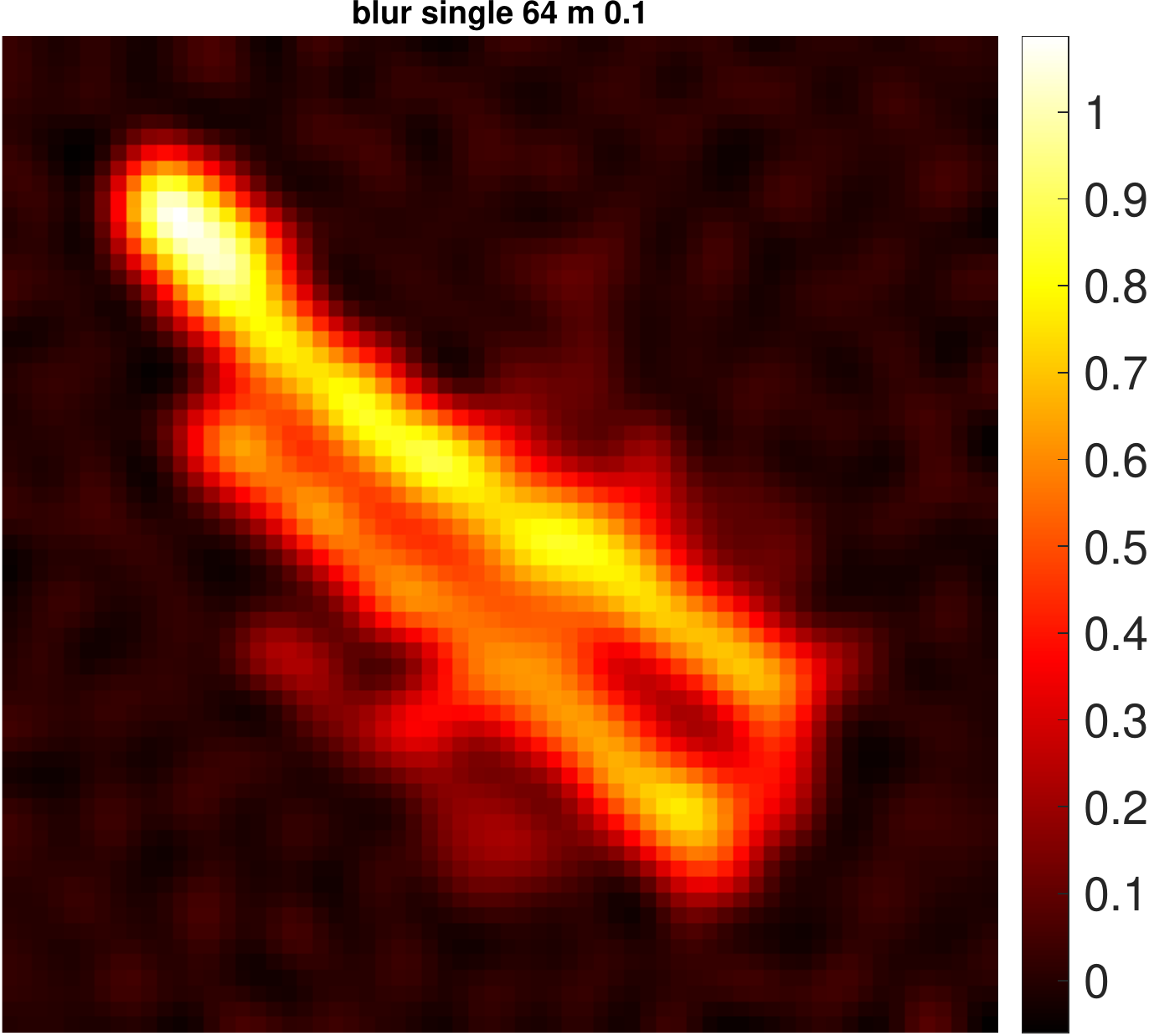}
  \caption{Image reconstruction in single precision of size 64 problem with mild blurring and 10\% noise.}\label{fig:Bcs_m_single_64_10PNoise}
\endminipage\hfill
\minipage{0.32\textwidth}%
  \includegraphics[width=\linewidth]{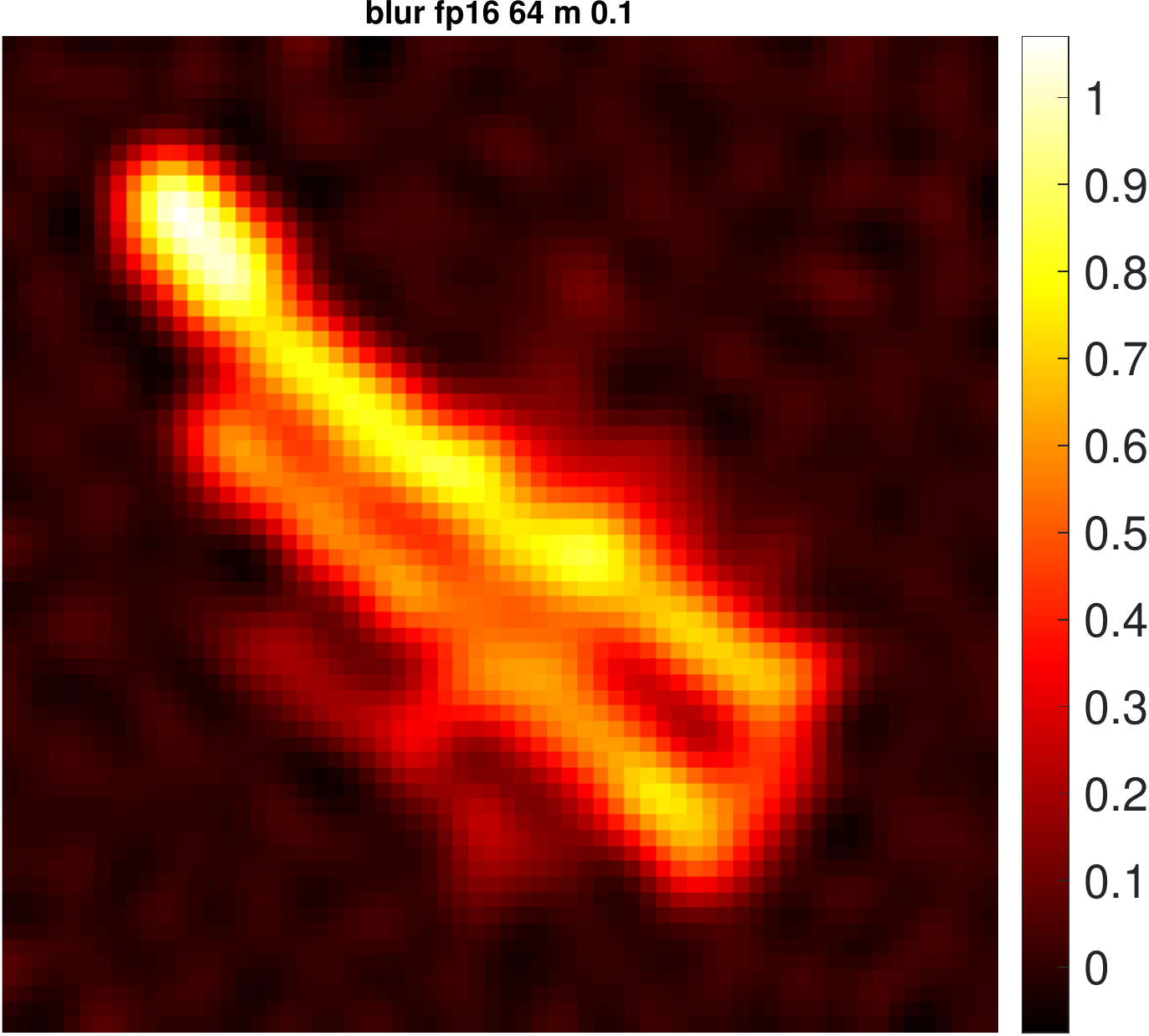}
  \caption{Image reconstruction in half precision of size 64 problem with mild blurring and 10\% noise.}\label{fig:Bcs_m_half_64_10PNoise}
\endminipage
\end{figure}

\begin{figure}[!htb]
\minipage{0.32\textwidth}
  \includegraphics[width=\linewidth]{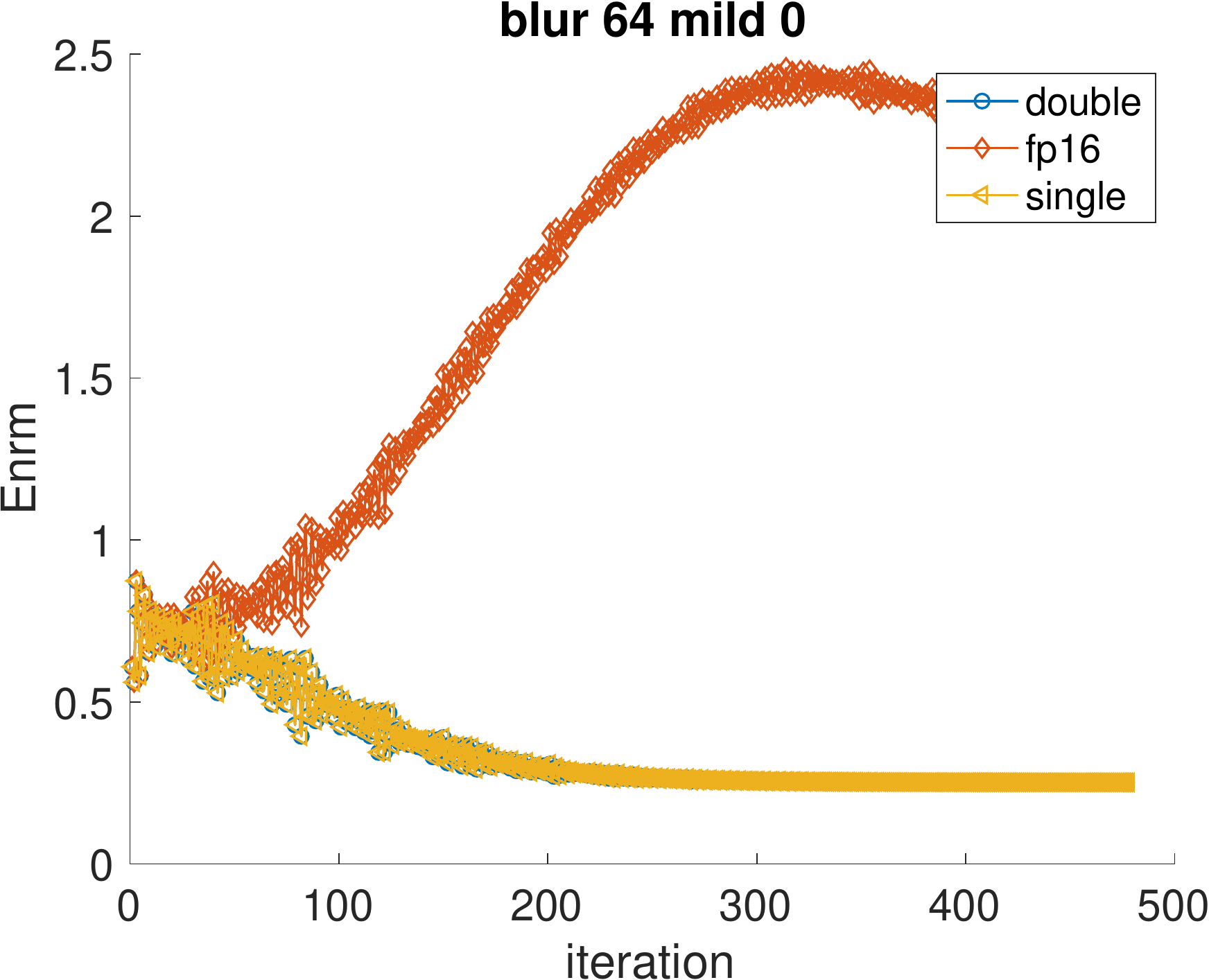}
\caption{Error norm of a size 64 problem with mild blurring of different precisions with 0 noise.}\label{fig:BEnrm_cs_size64_0noise_m}
\endminipage\hfill
\minipage{0.32\textwidth}
  \includegraphics[width=\linewidth]{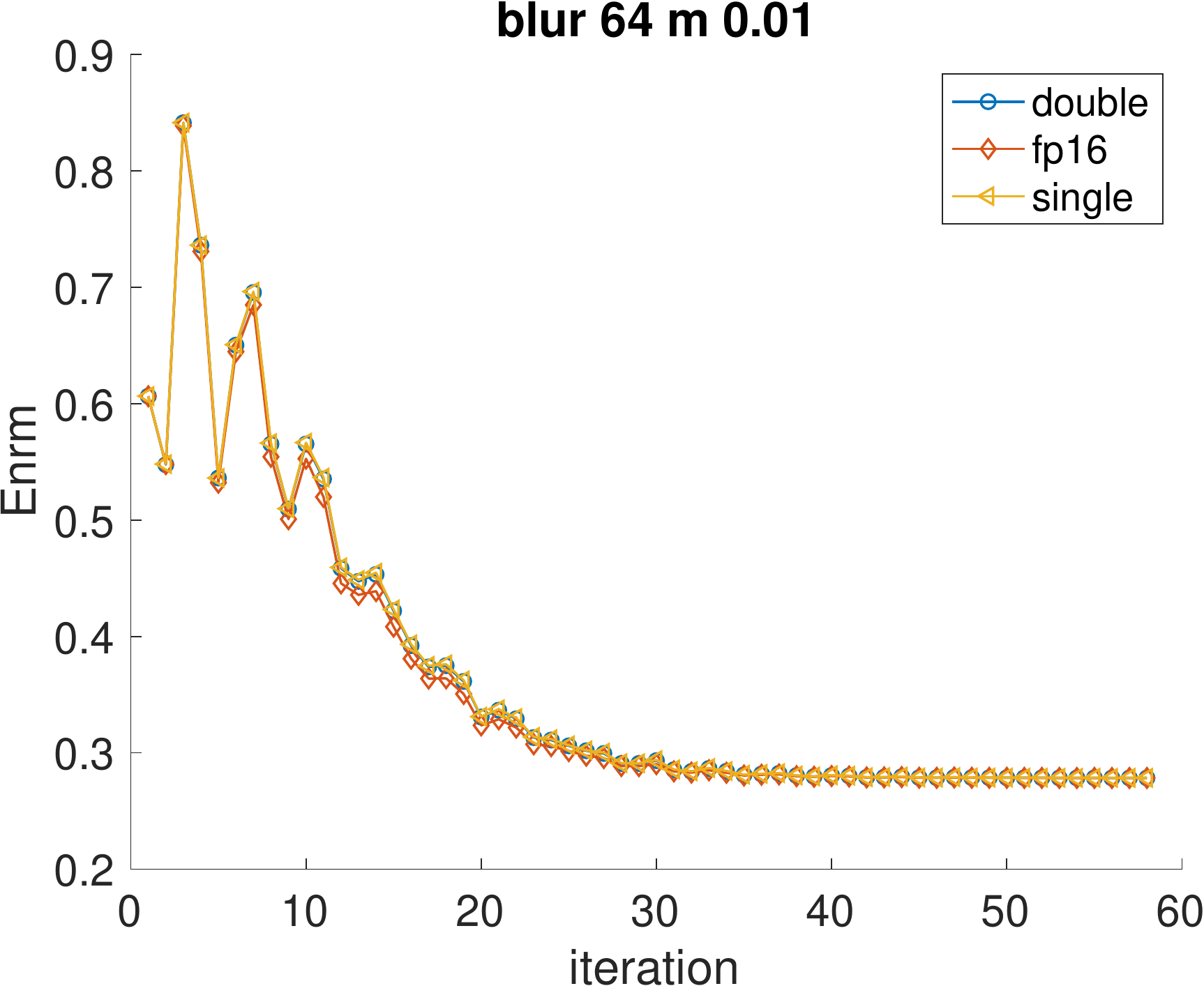}
  \caption{Error norm of a size 64 problem with mild blurring of different precisions with 1\% noise.}\label{fig:BEnrm_cs_size64_0.01noise_m}
\endminipage\hfill
\minipage{0.32\textwidth}%
  \includegraphics[width=\linewidth]{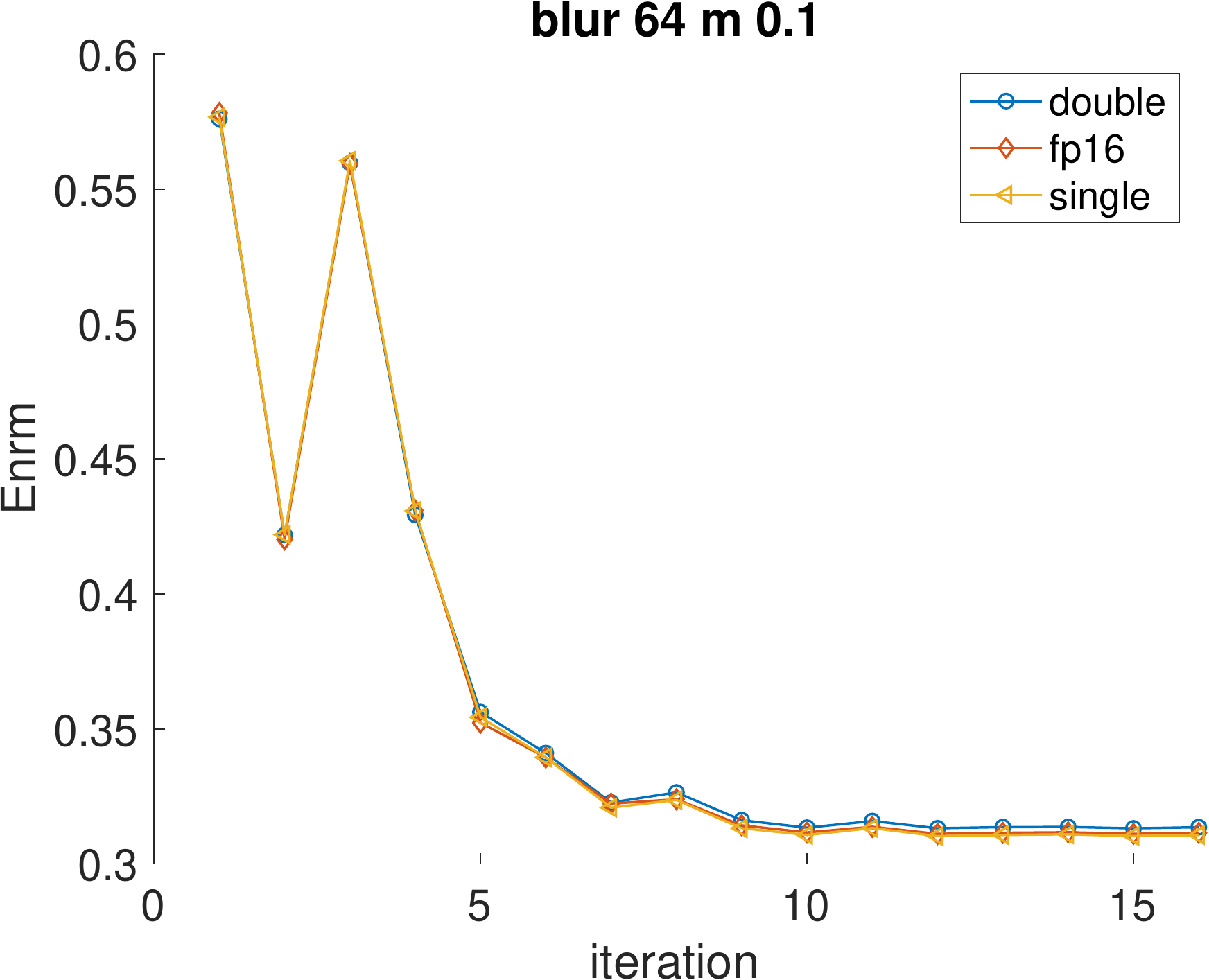}
 \caption{Error norm of a size 64 problem with mild blurring of different precisions with 10\% noise.}\label{fig:Benrms_cs_size64_10PNoise_m}
\endminipage
\end{figure}

\noindent With regularization, we successfully ran $478$ iterations in all three precisions, without the occurrence of NaNs. We plotted the results in Figures \ref{fig:Bcs_m_double_64_0noise}, \ref{fig:Bcs_m_single_64_0noise}, and \ref{fig:Bcs_m_half_64_0noise}. For half precision, the result is not clear as expected, and the image seems to be dominated by round-off errors. We plotted the error norm in Figure \ref{fig:BEnrm_cs_size64_0noise_m}.\\
\\
\noindent The error norms do not converge to a point like those generated by the CGLS algorithm. Instead, they oscillate at the beginning, and then converge for double and single precision. However, for half precision, the error norm increases rapidly, indicating that the round-off errors accumulate and take over. We suspect this is because the regularization parameter is too small, so we may need to develop better ways to determine a suitable regularization parameter for half precision. \\
\\We then added noise to $\bfb$ and displayed the resulting computed reconstructions in Figures \ref{fig:Bcs_m_double_64_10PNoise}, \ref{fig:Bcs_m_single_64_10PNoise}, and \ref{fig:Bcs_m_half_64_10PNoise}.
Here we only showed results with $10\%$ noise; results for other noise levels are consistent with those shown in previous sections. Surprisingly, the result in half precision looks similar to the one in double or single precision, and it is much better than the image without noise. Moreover, the clear result is acquired by only $16$ iterations, illustrating the efficiency of the CS method. 

\noindent By observing the error norm plots in Figure \ref{fig:Benrms_cs_size64_10PNoise_m}, we also see that there seems to be very little difference in the convergence behavior of CS for the different precision levels when the noise level is high. The more noise we have, the closer are the results between double and half precision.

\subsubsection{Tomography Reconstruction}

\begin{figure}[!htb]
\minipage{0.32\textwidth}
  \includegraphics[width=\linewidth]{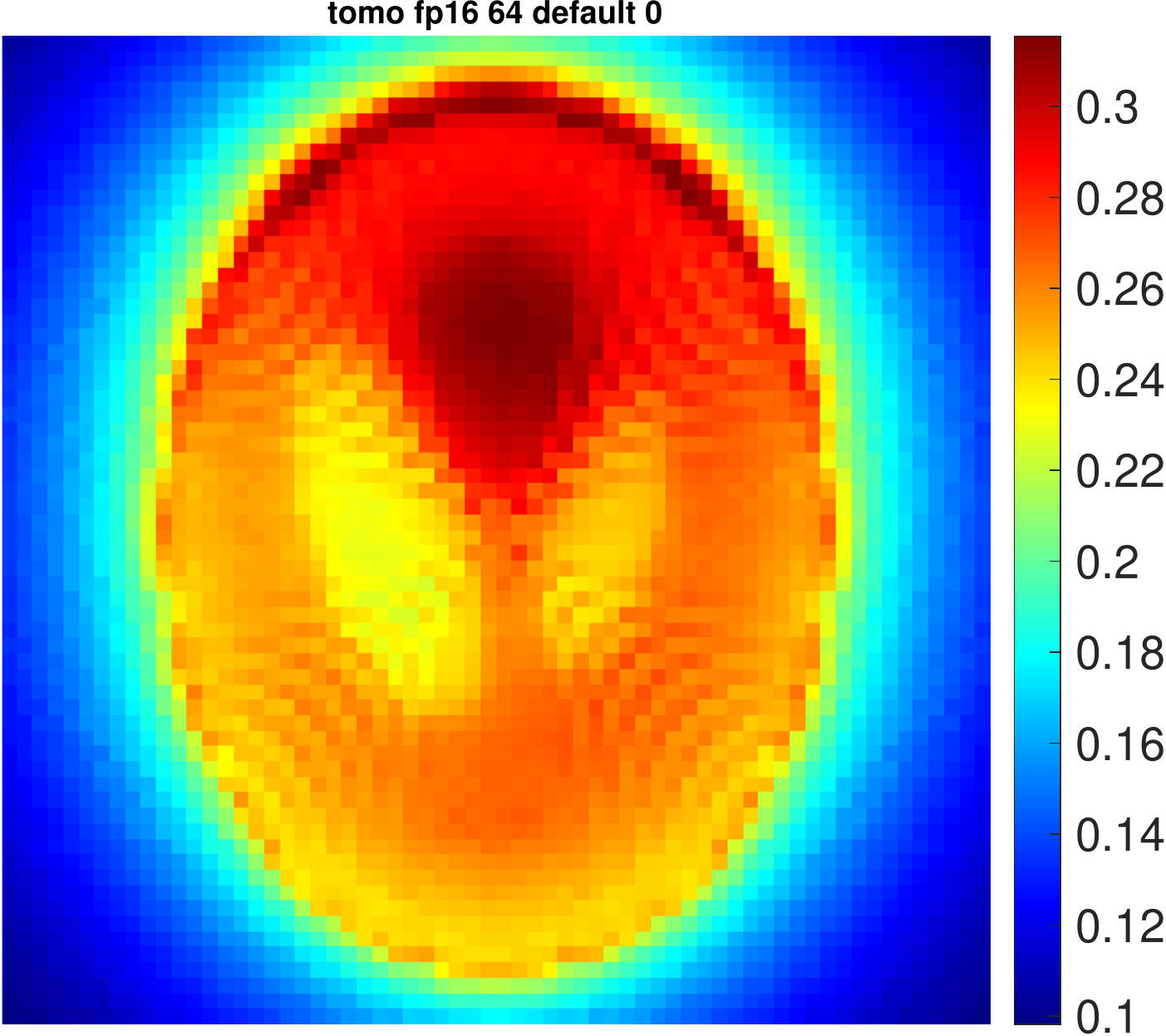}
  \caption{Reconstruction at first iteration, half precision, zero noise (after rescaling).}\label{fig:Tcs_half_0noise_it1_64}
\endminipage\hfill
\minipage{0.32\textwidth}
  \includegraphics[width=\linewidth]{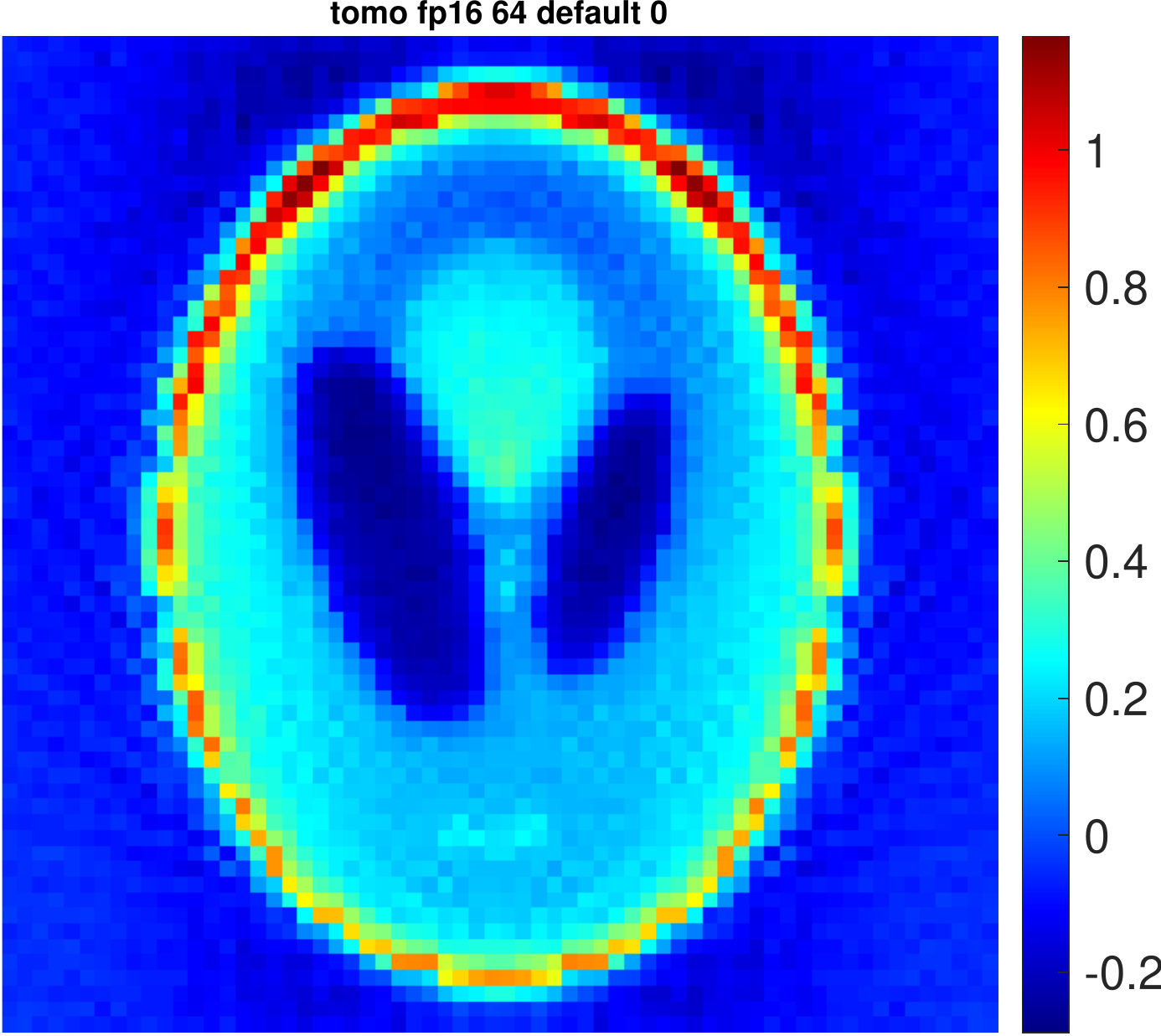}
  \caption{Reconstruction at best iteration, half precision, zero noise (after rescaling).}\label{fig:Tcs_half_0noise_it5_64}
\endminipage\hfill
\minipage{0.32\textwidth}%
  \includegraphics[width=\linewidth]{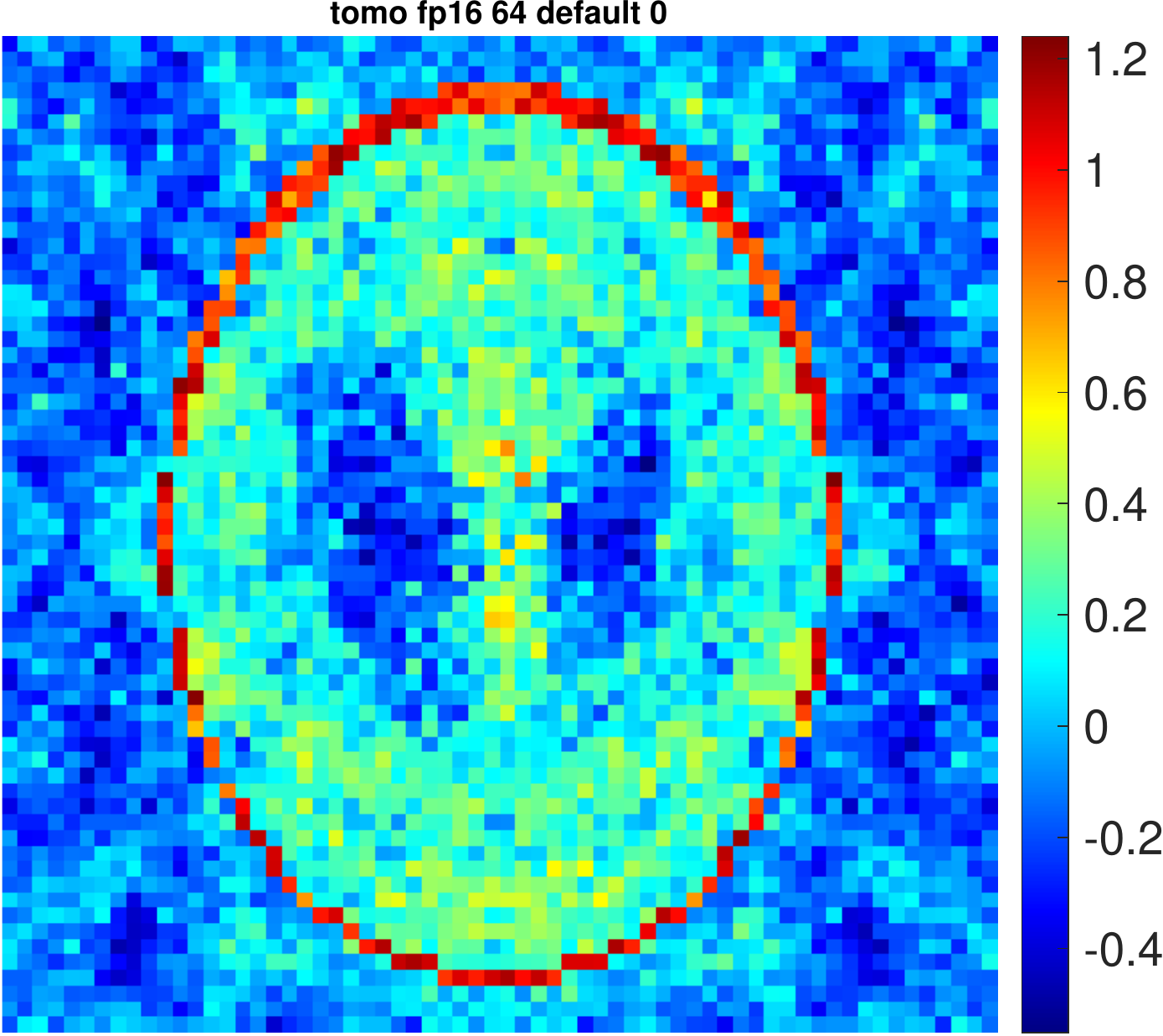}
  \caption{Reconstruction at $250^{th}$ iteration, half precision, zero noise (after rescaling).}\label{fig:Tcs_half_0noise_it250_64}
\endminipage
\end{figure}

\begin{figure}[!htb]
\minipage{0.32\textwidth}
  \includegraphics[width=\linewidth]{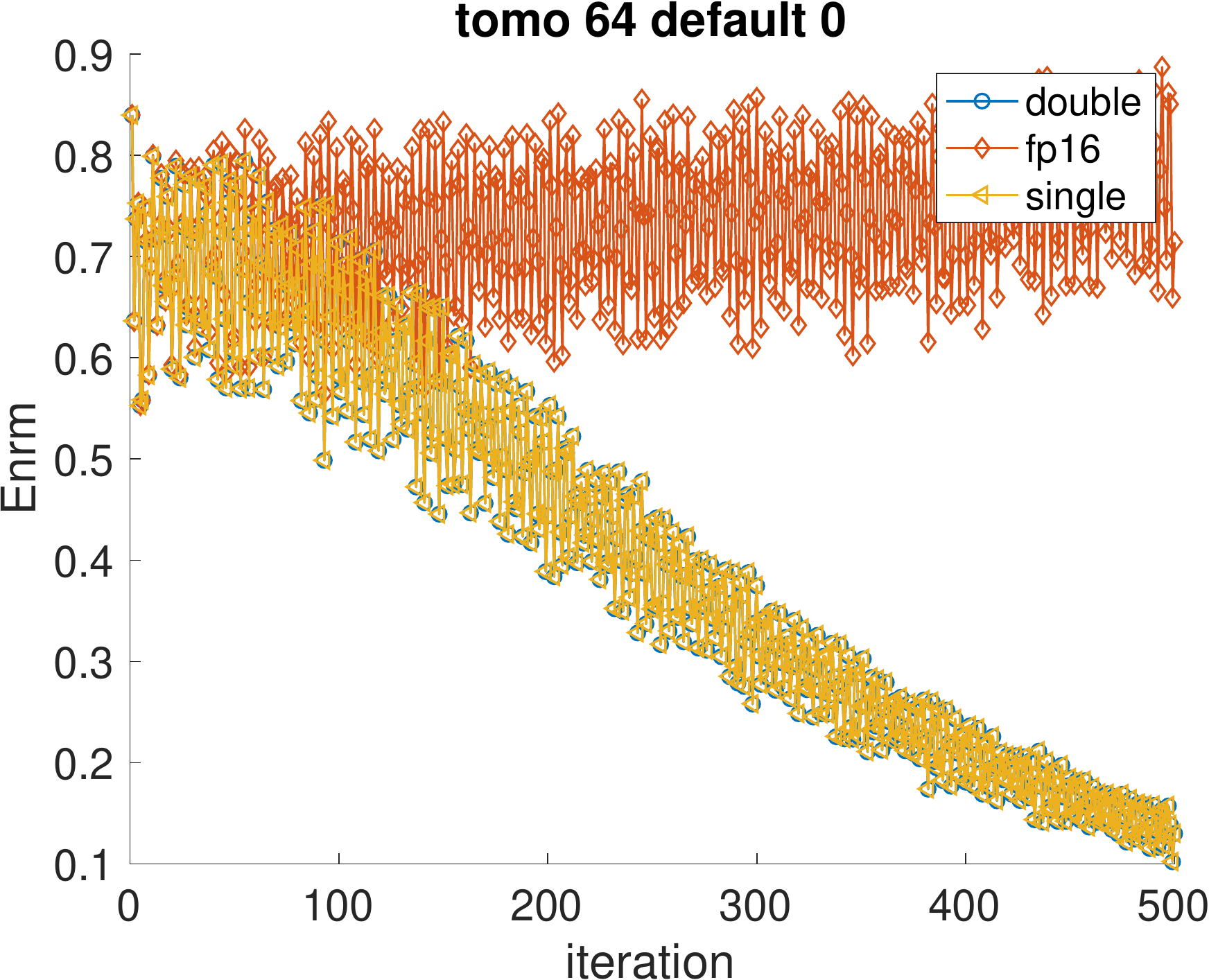}
  \caption{Error norm of a size 64 problem at different precisions with zero noise.}\label{fig:size64Tomo0Noise}
\endminipage\hfill
\minipage{0.32\textwidth}
  \includegraphics[width=\linewidth]{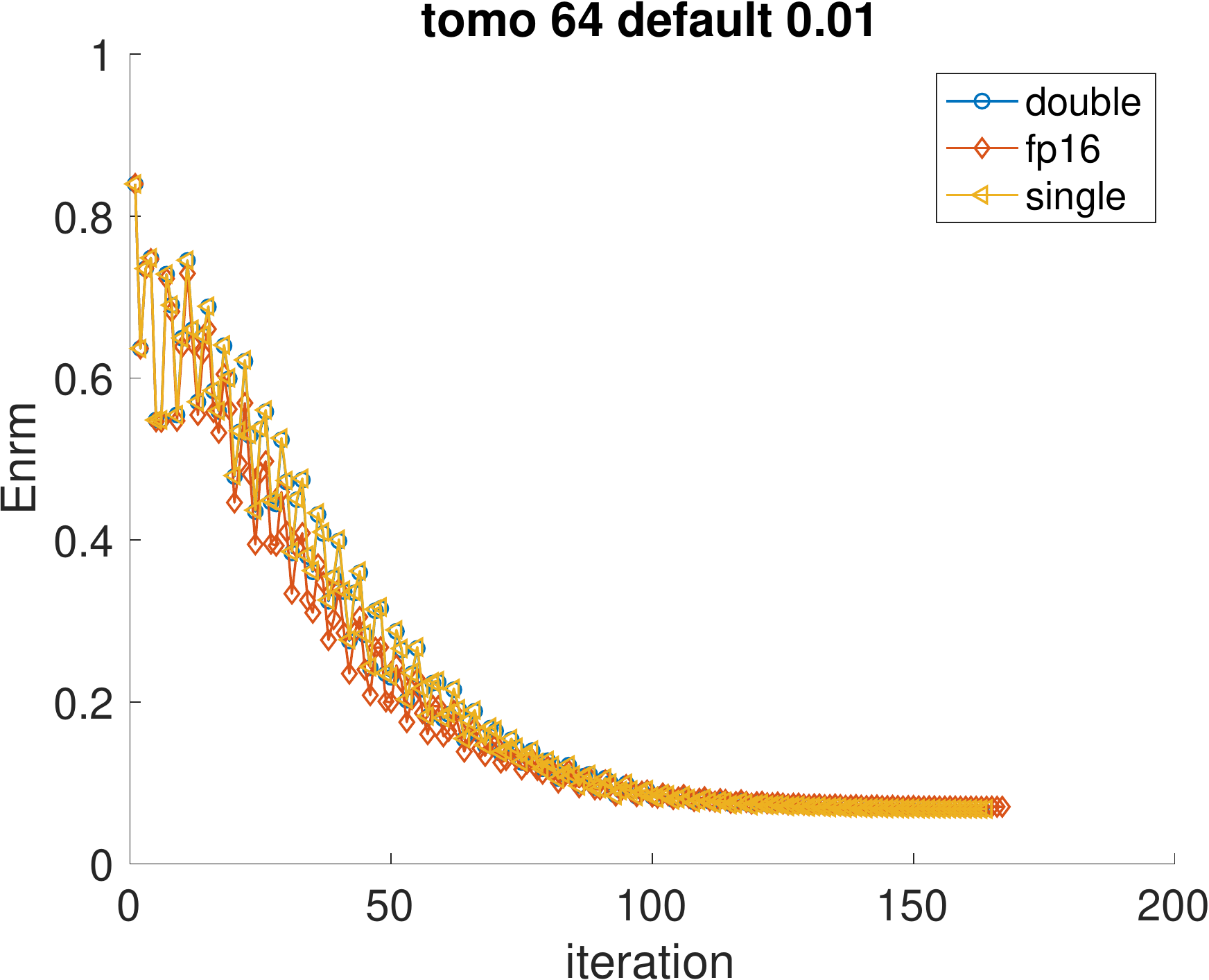}
  \caption{Error norm of a size 64 problem at different precisions with 1\% noise.}\label{fig:size64Tomo1Noise}
\endminipage\hfill
\minipage{0.32\textwidth}%
  \includegraphics[width=\linewidth]{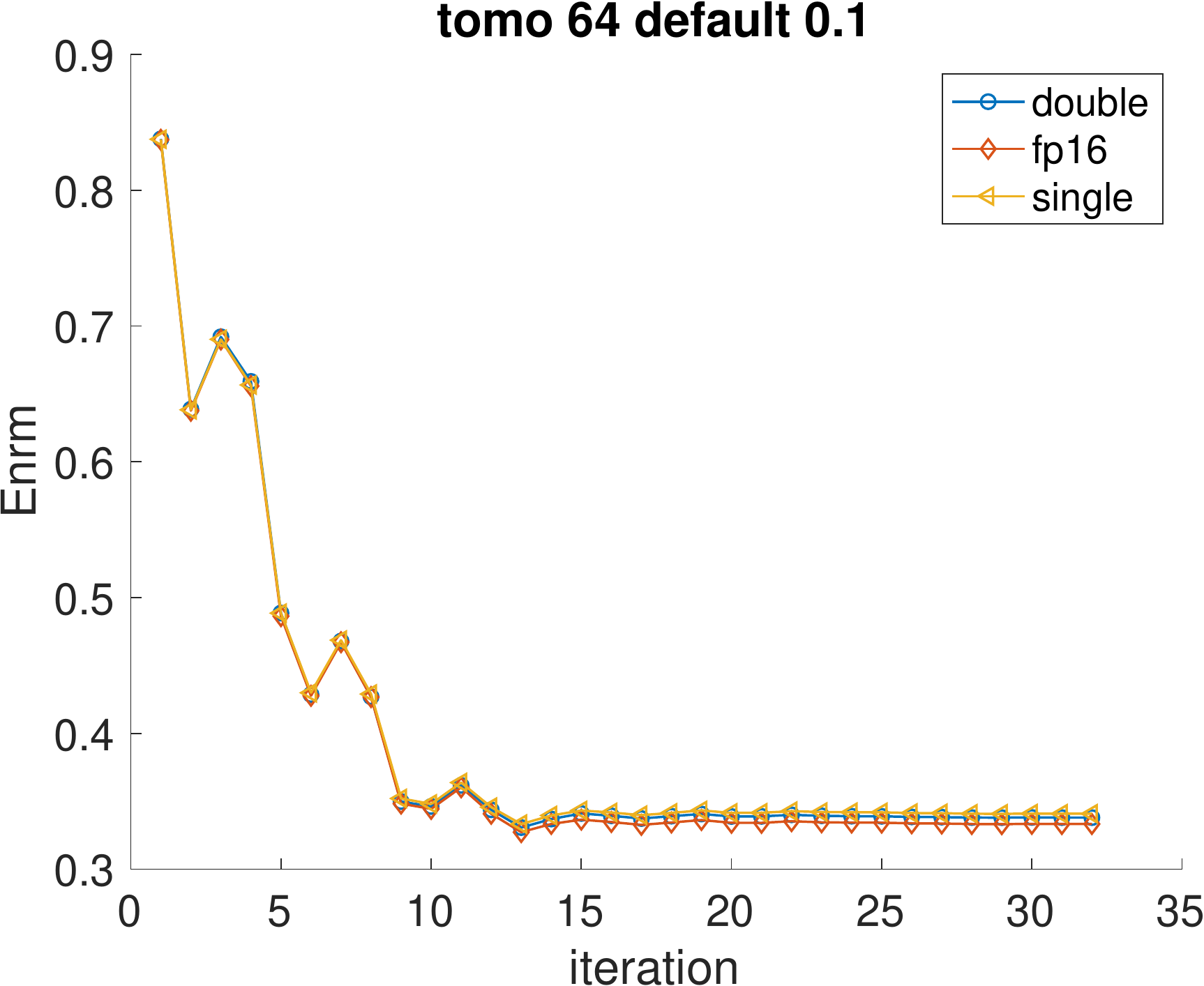}
\caption{Error norm of a size 64 problem at different precisions with 10\% noise.}\label{fig:size64Tomo10Noise}
\endminipage
\end{figure}

Although CS avoids the computation of inner products, overflow still occurs at half precision during the calculation of matrix-vector multiplications for the tomography reconstruction problem. NaNs start to occur at the second iteration. Figure \ref{fig:Tcs_half_0noise_it1_64} shows the resulting image at the end of the first iteration at half precision for a problem with $10\%$ noise, which is already reasonably good with visible shapes and boundaries.\\
\begin{figure}
\begin{center}
\includegraphics[width=7cm]{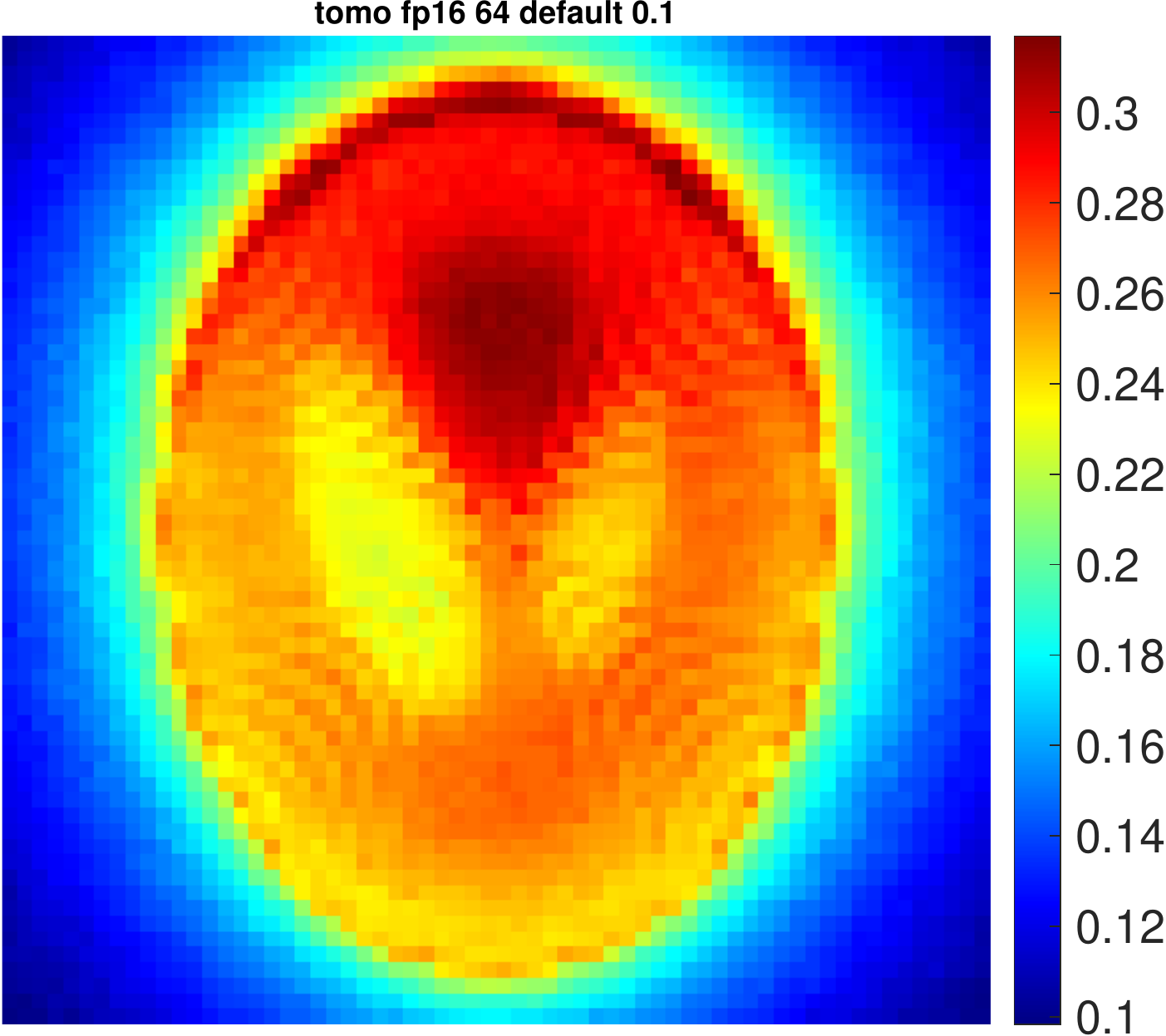}
\caption{Half precision, size 64, 10\% noise, without rescale.}\label{fig:Tcs_half_10Pnoise_64}
\end{center}
\end{figure}
\\Due to overflow, we rescaled the problem as for CGLS so that CS could run more iterations for half precision. The algorithm successfully runs to the end without occurrences of overflow. In the figures below, we show the results at the last iteration for problems with different precision and noise levels.\\
\begin{figure}[!htb]
\minipage{0.32\textwidth}
  \includegraphics[width=\linewidth]{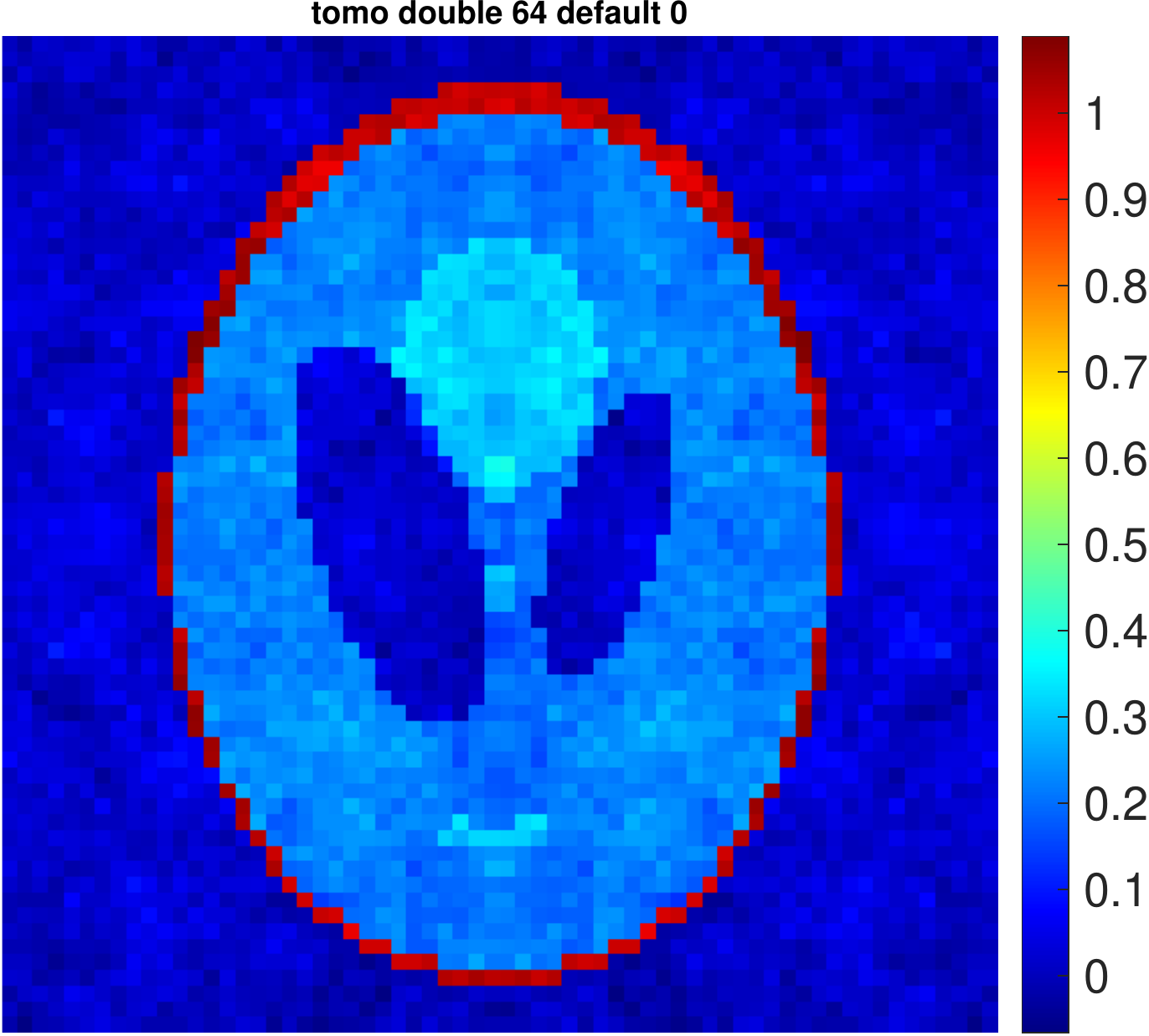}
  \caption{Double precision, size 64, zero noise.}\label{fig:Tcs_double_0noise_64}
  \includegraphics[width=\linewidth]{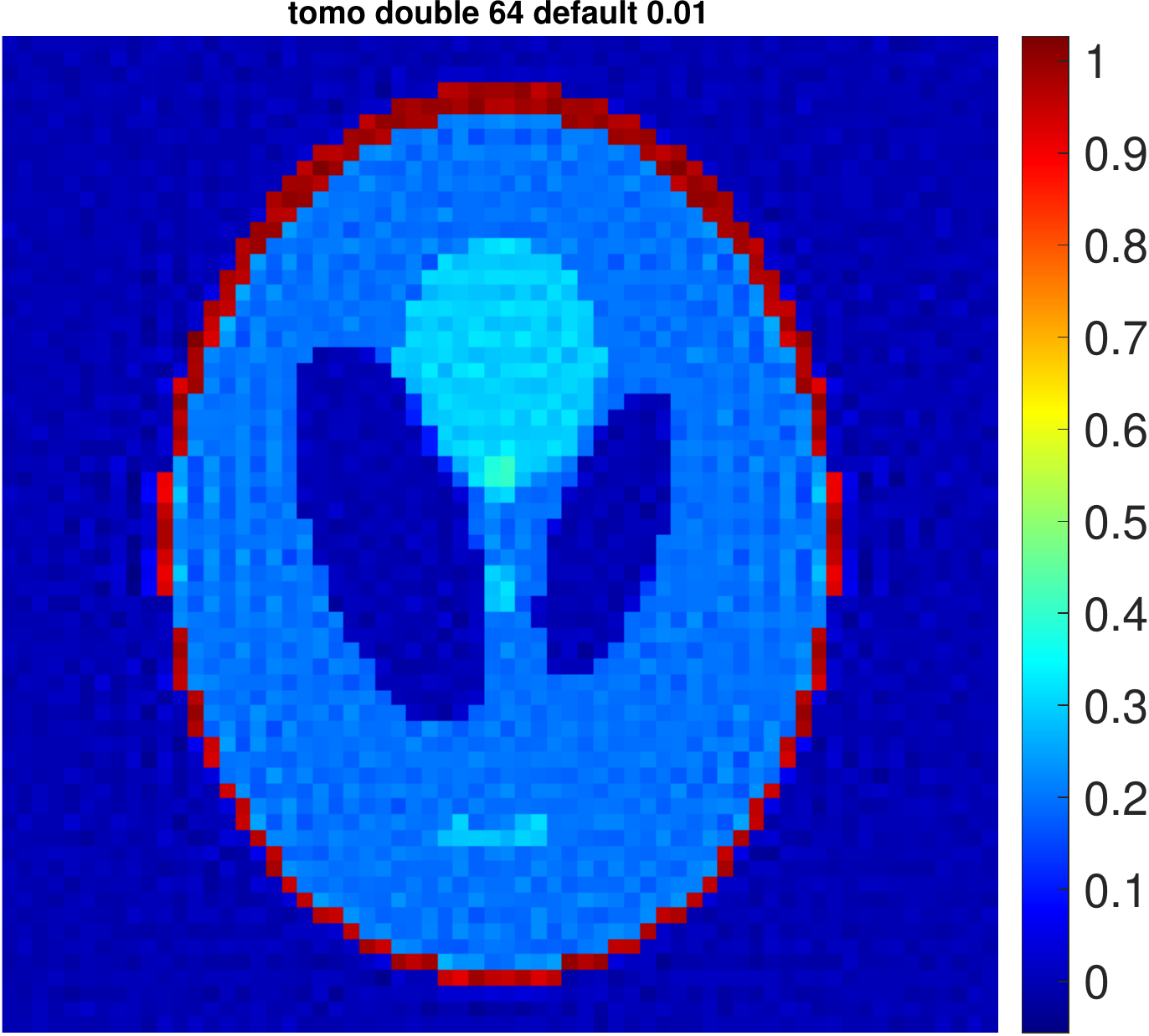}
  \caption{Double precision, size 64, 1\% noise.}\label{fig:Tcs_double_1pnoise_64}
  \includegraphics[width=\linewidth]{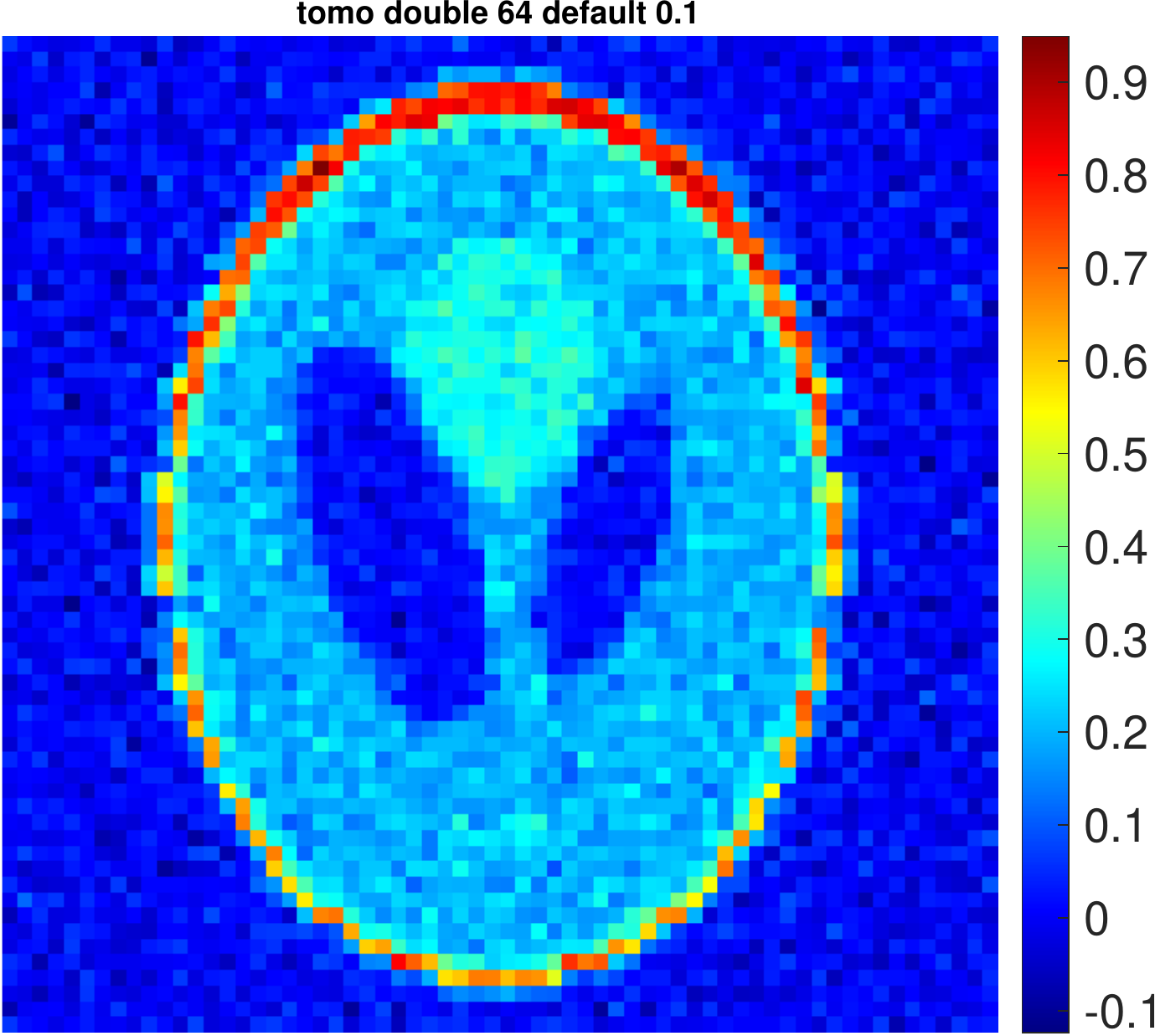}
  \caption{Double precision, size 64, 10\% noise.}\label{fig:Tcs_double_10pnoise_64}
\endminipage\hfill
\minipage{0.32\textwidth}
  \includegraphics[width=\linewidth]{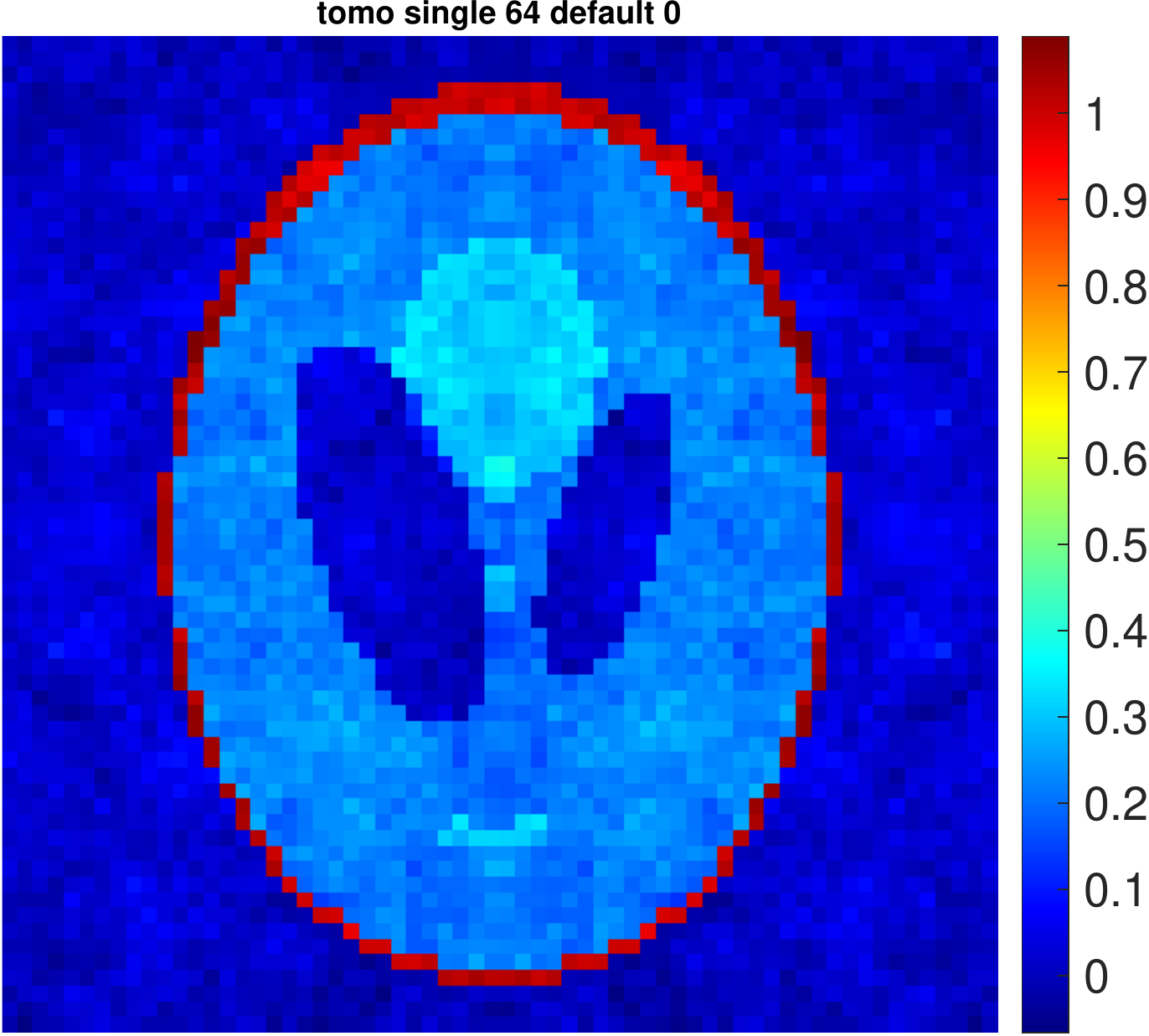}
  \caption{Single precision, size 64, zero noise.}\label{fig:Tcs_single_0noise_64}
  \includegraphics[width=\linewidth]{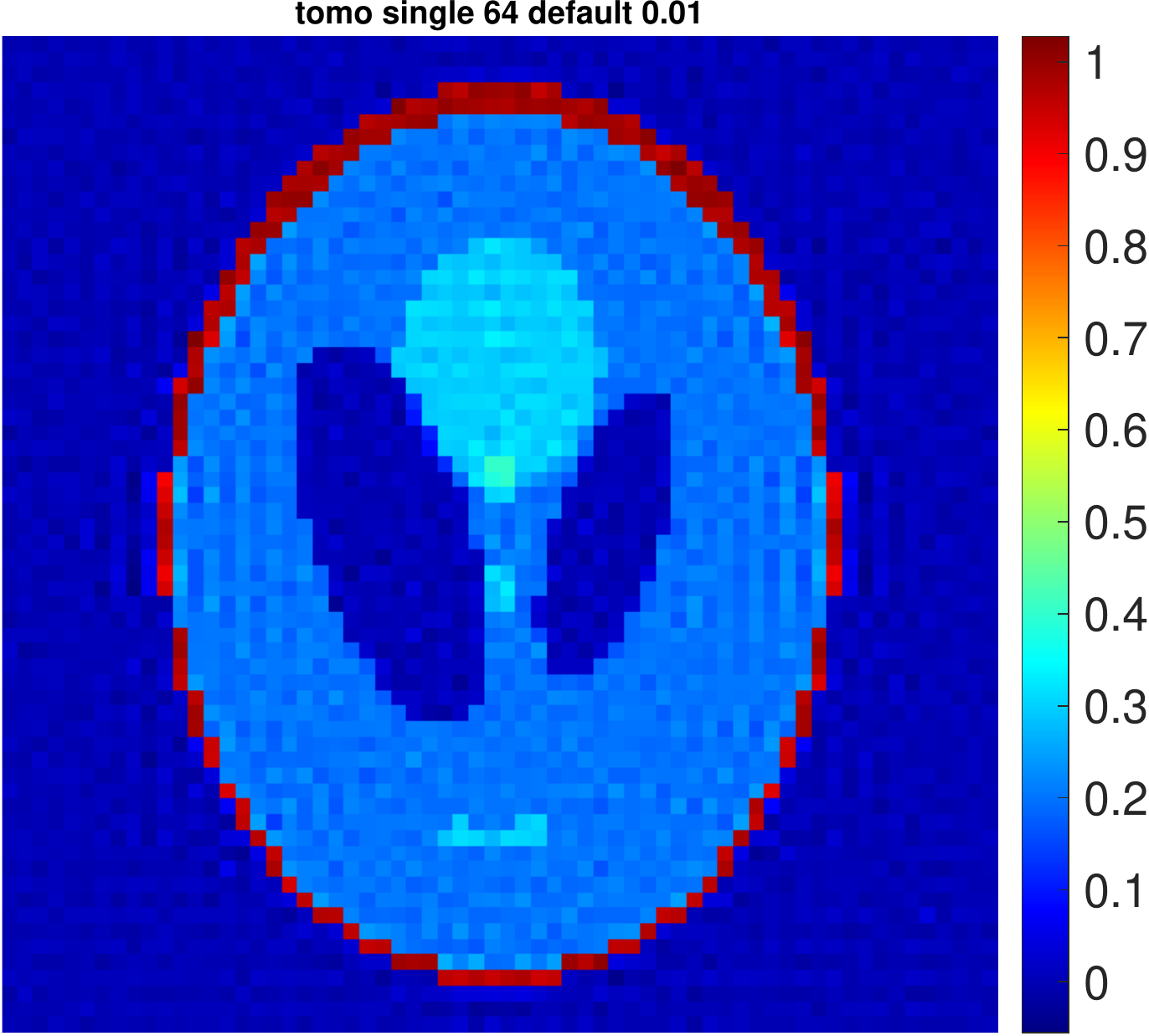}
  \caption{Single precision, size 64, 1\% noise.}\label{fig:Tcs_single_1pnoise_64}
  \includegraphics[width=\linewidth]{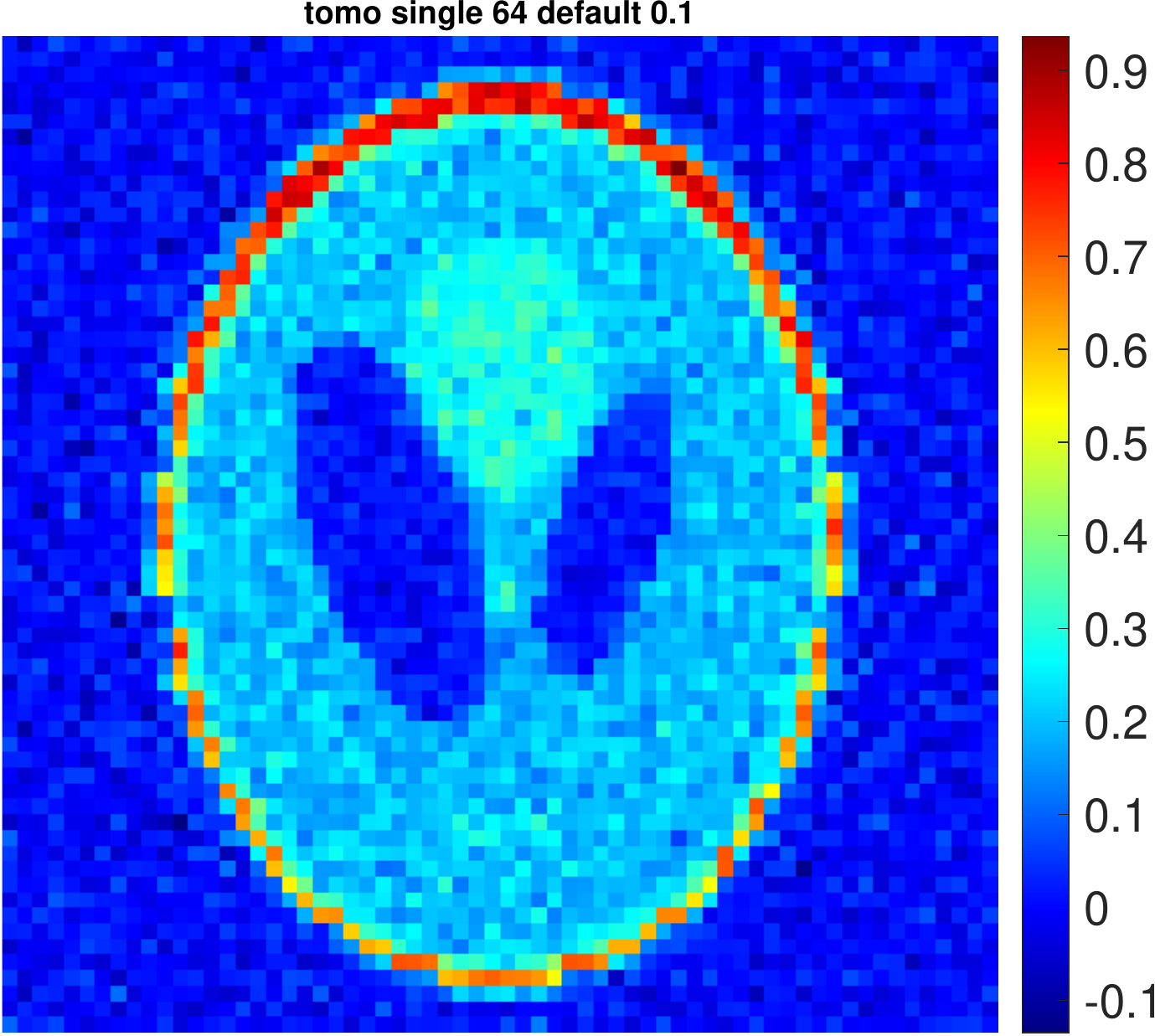}
  \caption{Single precision, size 64, 10\% noise.}\label{fig:Tcs_single_10Pnoise_64}
\endminipage\hfill
\minipage{0.32\textwidth}%
  \includegraphics[width=\linewidth]{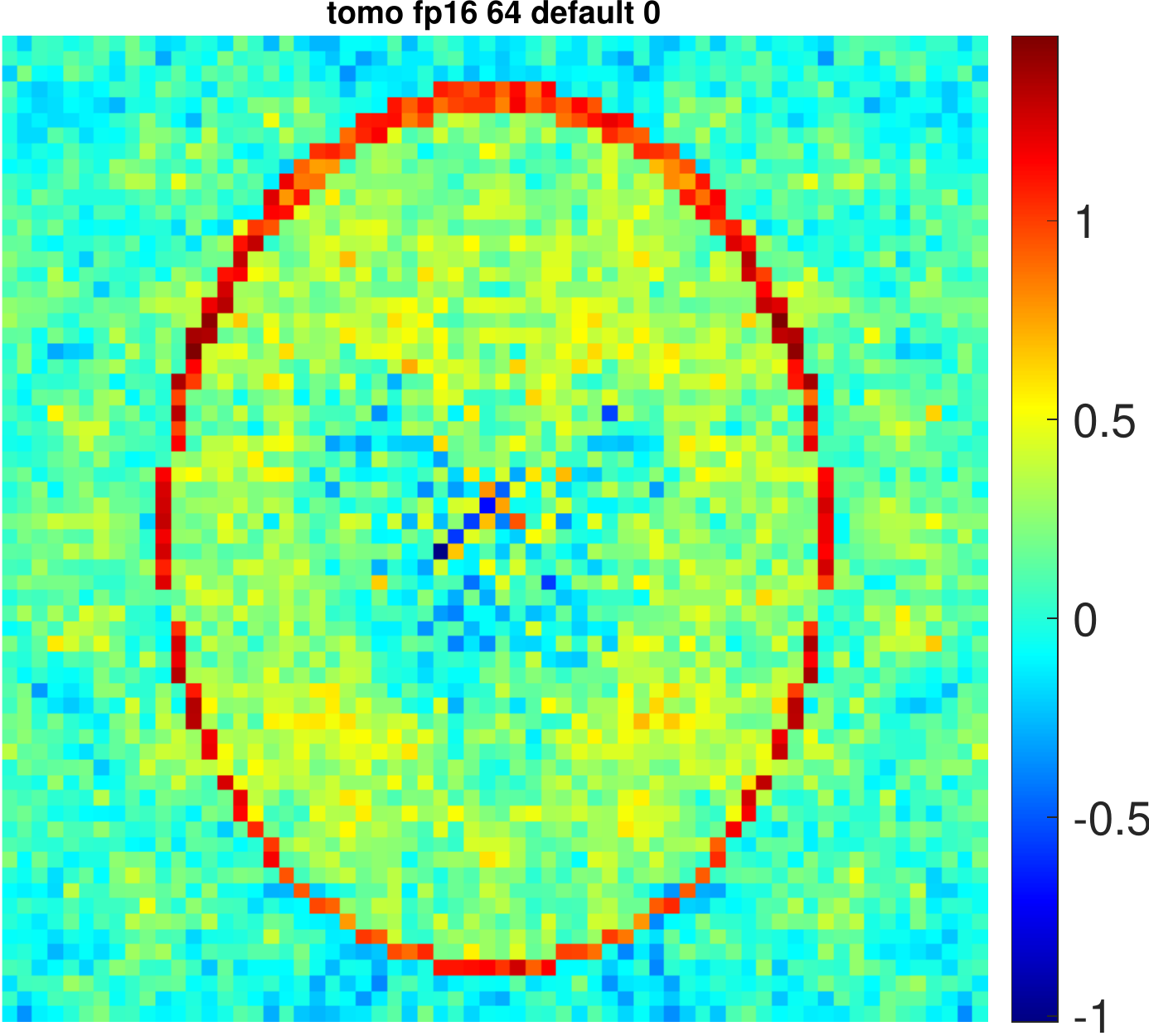}
  \caption{Half precision, size 64, zero noise (after rescaling).}\label{fig:Tcs_half_0noise_64}
  \includegraphics[width=\linewidth]{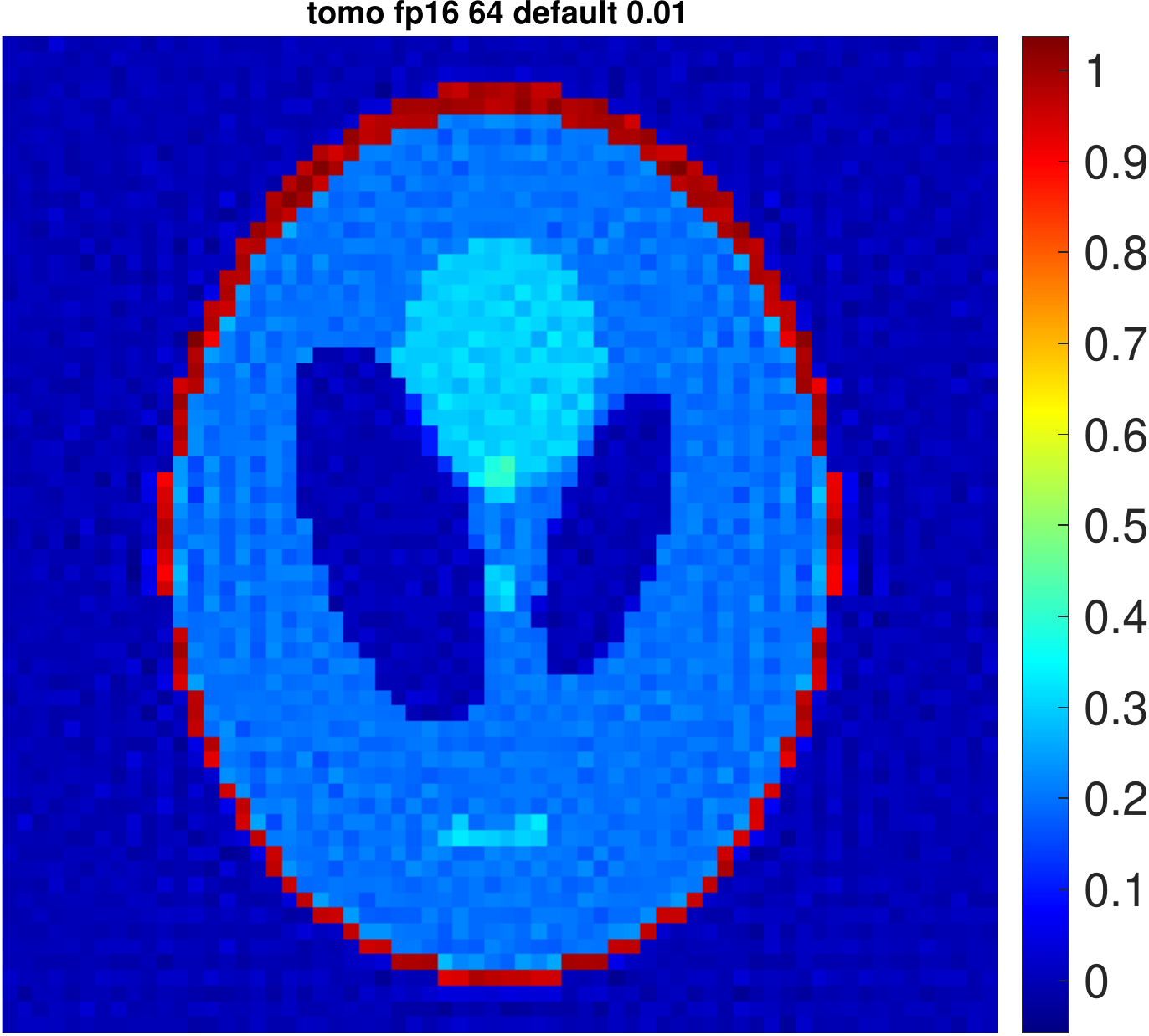}
  \caption{Half precision, size 64, 1\% noise (after rescaling).}\label{fig:Tcs_half_1pnoise_64}
  \includegraphics[width=\linewidth]{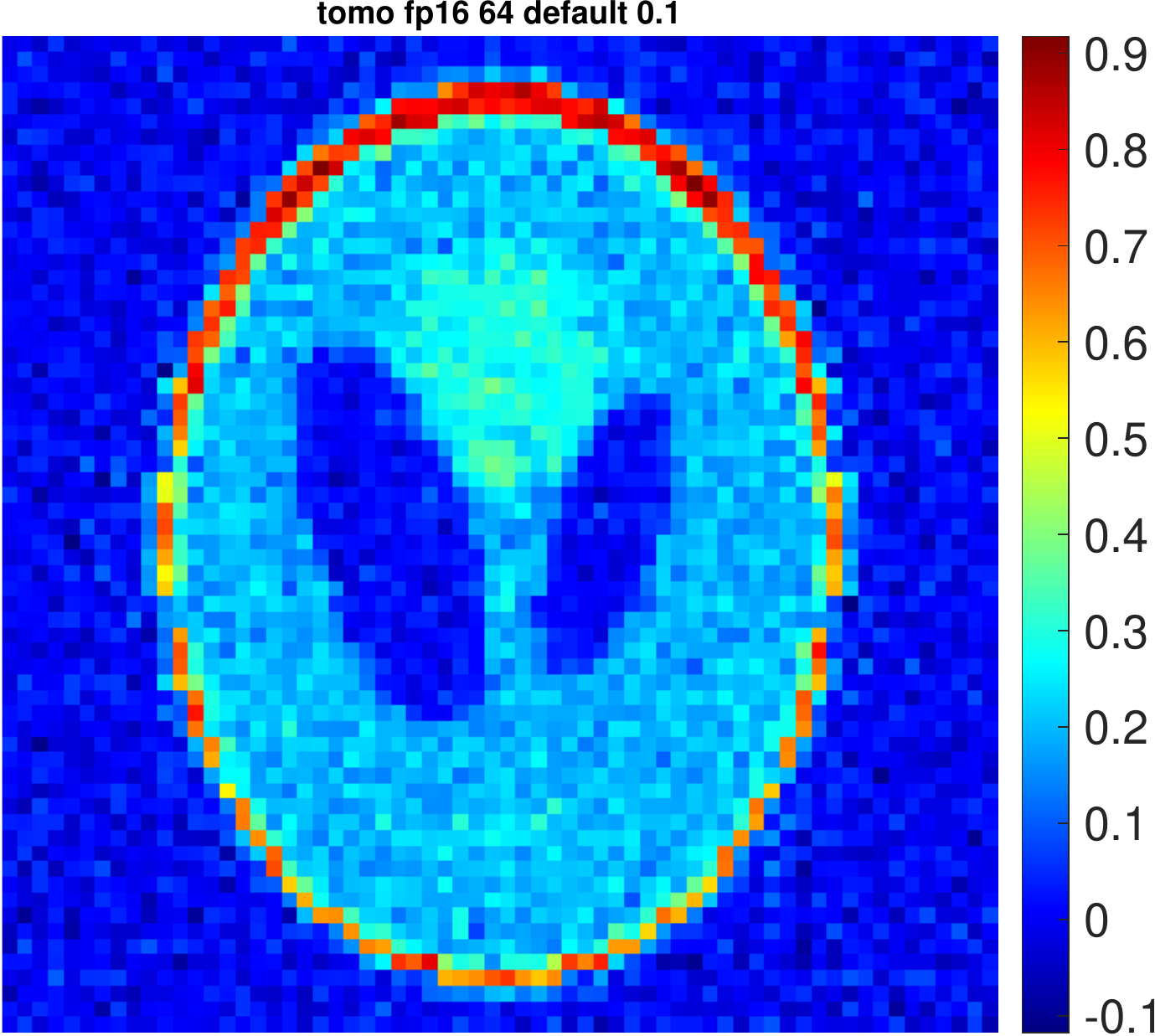}
  \caption{Half precision, size 64, 10\% noise (after rescaling).}\label{fig:Tcs_half_10Pnoise_64}
\endminipage
\end{figure}
\\After rescaling, when there is noise on the right hand side $\bfb$ of the problem, the resulting image at low precision is as valid as images produced at high precision. With $10\%$ noise (Figures \ref{fig:Tcs_double_10pnoise_64}, \ref{fig:Tcs_single_10Pnoise_64}, and \ref{fig:Tcs_half_10Pnoise_64}), the CS algorithm ran $32$ iterations, and produced good reconstructions with clear boundaries and background for all three precision levels.\\
\\However, when there is zero noise, the result at half precision is a poor reconstruction. The background and the object both look blurry and blend together, and the object contents are hardly visible. We believe the dissatisfying image is a combined result of accumulation of round-off errors and under-regularization.\\
\\After a closer look at the result of each iteration, we noticed that the reconstructions in the first few iterations look smooth and improved as the iteration moves on, but at some point they start to become noisy and blurry. Unlike CGLS where the output image at each iteration refines steadily, the results from CS have more oscillations from iteration to iteration, and it has a general trend of becoming noisy as the iteration goes on. Below we show the resulting images at the first iteration (Figure \ref{fig:Tcs_half_0noise_it1_64}), the $5^{th}$ iteration (Figure \ref{fig:Tcs_half_0noise_it5_64}) where the error norm is the smallest, and the $250^{th}$ iteration (Figure \ref{fig:Tcs_half_0noise_it250_64}). At the $250^{th}$ iteration the image is already very noisy as round-off errors accumulated along the way.
\\
\\
Figures \ref{fig:size64Tomo0Noise}, \ref{fig:size64Tomo1Noise} and \ref{fig:size64Tomo10Noise} present the error norms at different precision levels across different noise levels. As expected, the error norms overlap for test problems with noise, implying that the quality of the reconstructions is similar. However, for the noise-free problems, the error norms follow a decreasing trend for double and single precision and a slightly increasing trend for half precision. Besides, the norms oscillate more as noise level decreases, which corresponds with our observation from the resulting images.
\\

\section{Discussion}

\begin{figure}
\begin{center}
\includegraphics[width=7cm]{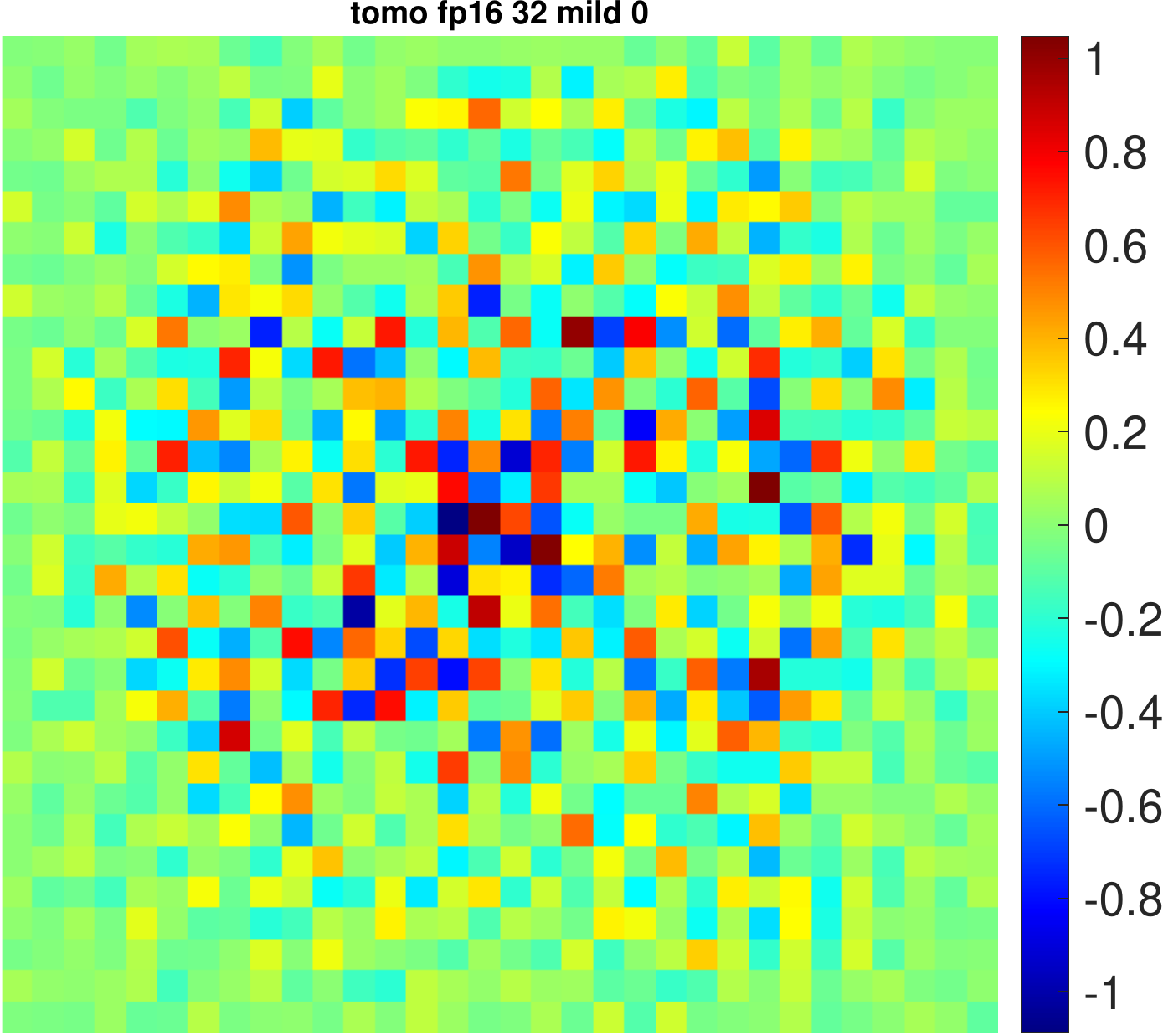}
\caption{Half precision, size 32, zero noise.}\label{fig:size32TomoNoNoise}
\end{center}
\end{figure}

We incorporated CGLS with regularization to compare its performance with that of CS more fairly. The effect of regularization is apparent, as the resulting image is less noisy, especially for low precision. As the noise level increases, the difference between high precision and low precision becomes less significant. Yet still, rescaling is necessary for half precision to avoid overflow for the tomography reconstruction problem. Below we showed the result of CGLS with regularization for the two test problems.
\begin{figure}[!htb]
\minipage{0.32\textwidth}
  \includegraphics[width=\linewidth]{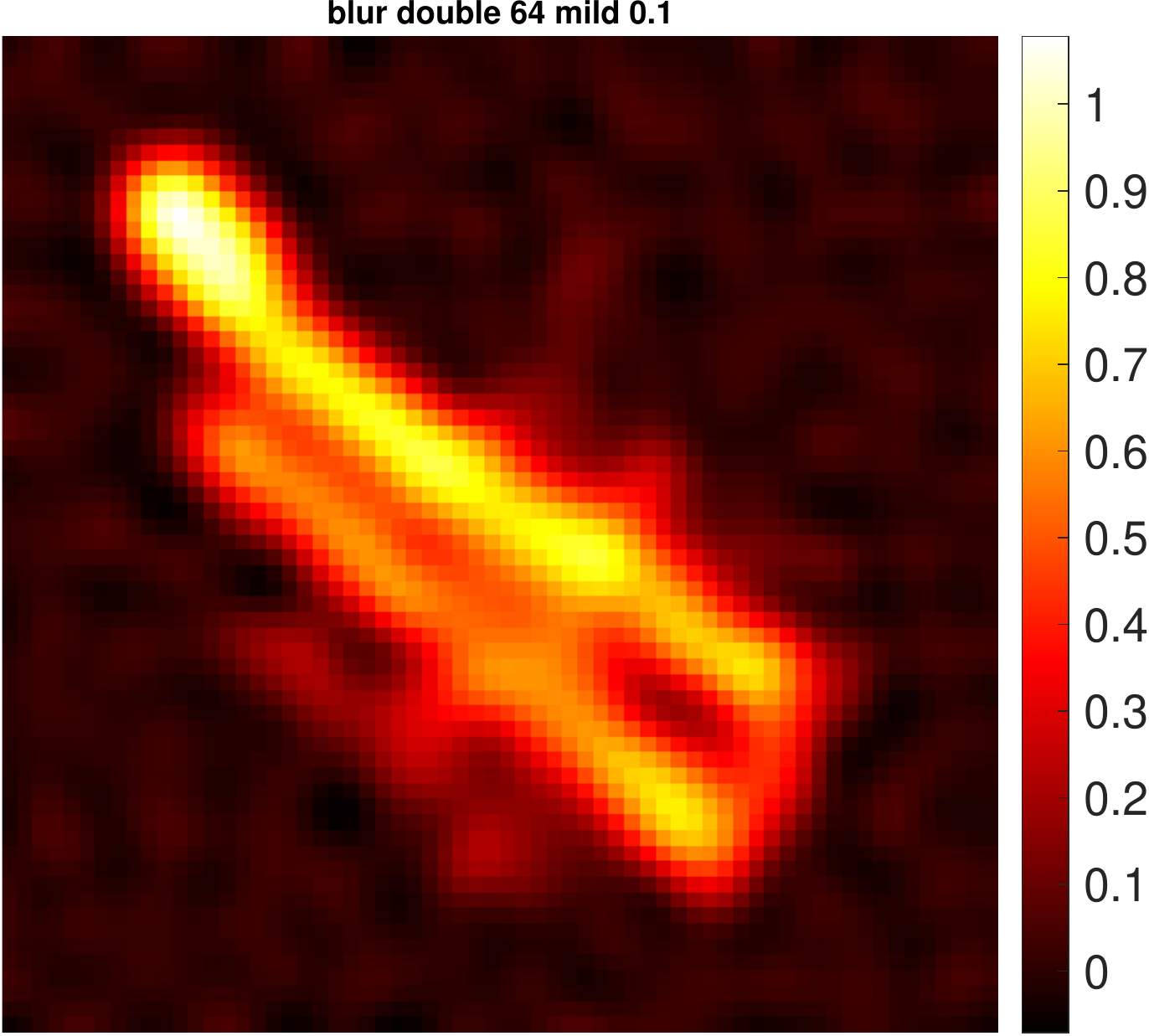}
  \caption{CGLS with regularization, double precision, 10\% noise.}\label{fig:Bcgreg_double_0.1noise_64}
\endminipage\hfill
\minipage{0.32\textwidth}
  \includegraphics[width=\linewidth]{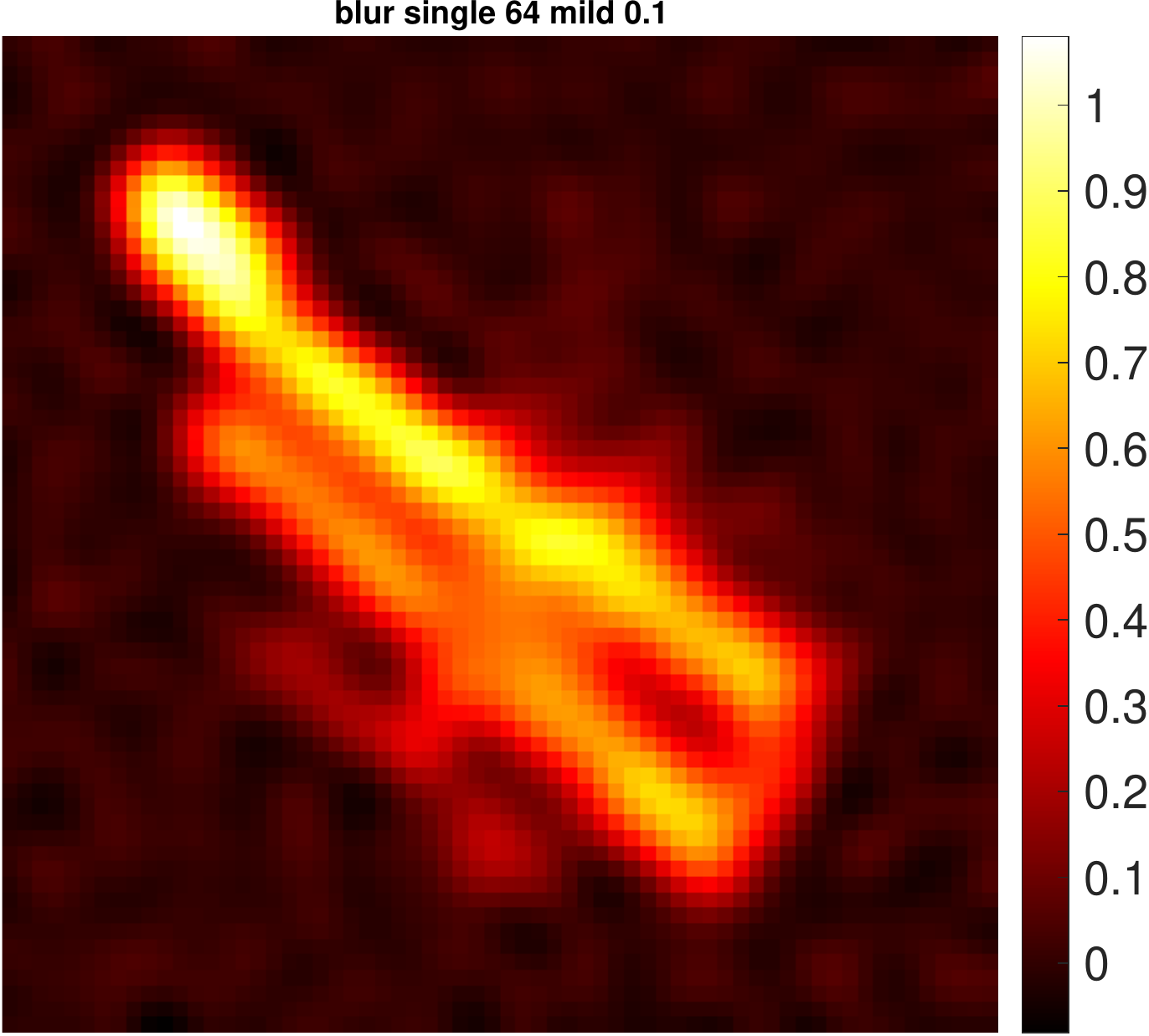}
  \caption{CGLS with regularization, single precision, 10\% noise.}\label{fig:Bcgreg_single_0.1noise_64}
\endminipage\hfill
\minipage{0.32\textwidth}
  \includegraphics[width=\linewidth]{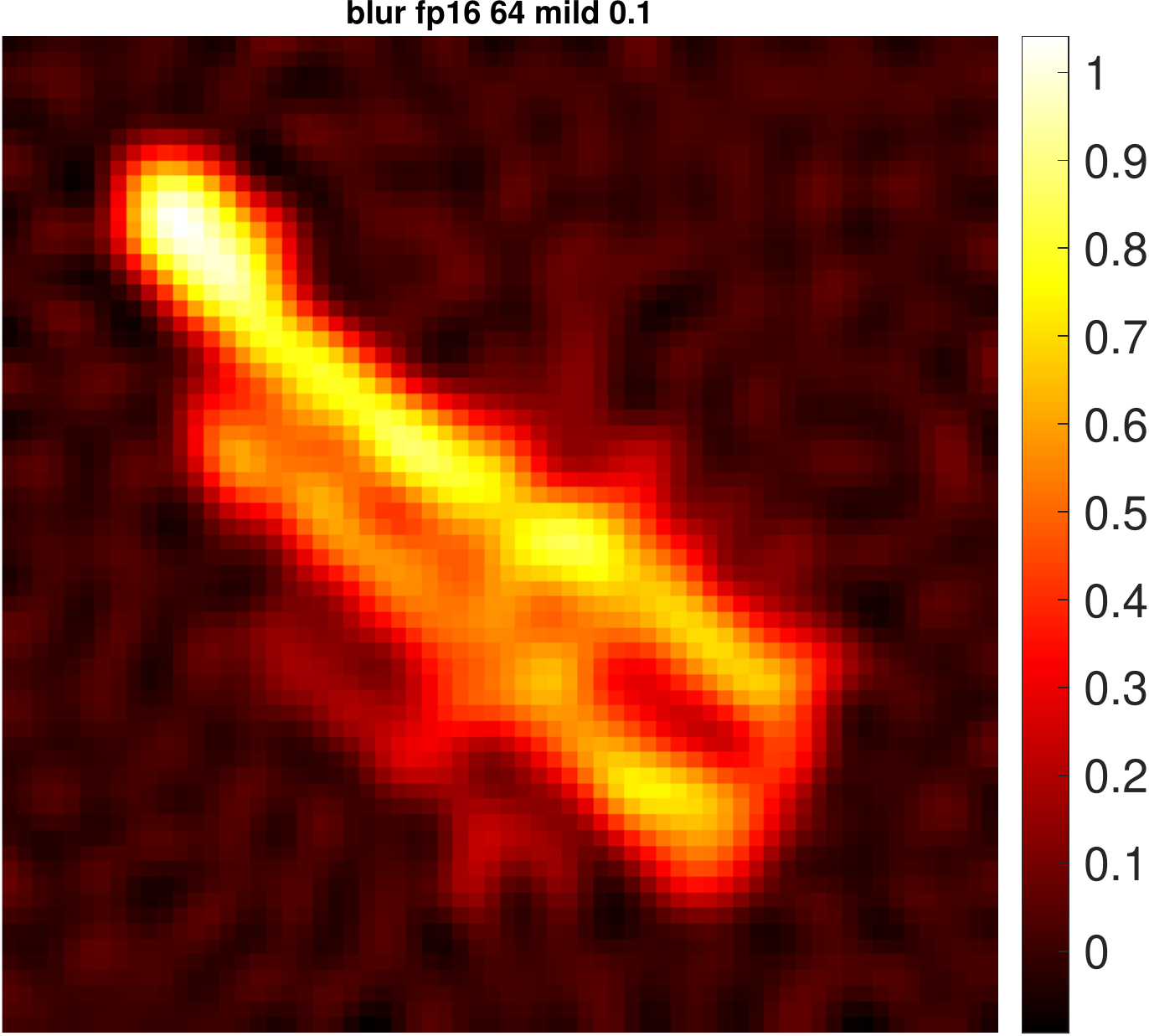}
  \caption{CGLS with regularization, half precision, 10\% noise.}\label{fig:Bcgreg_half_0.1noise_64}
\endminipage
\end{figure}
\begin{figure}[!htb]
\minipage{0.32\textwidth}
  \includegraphics[width=\linewidth]{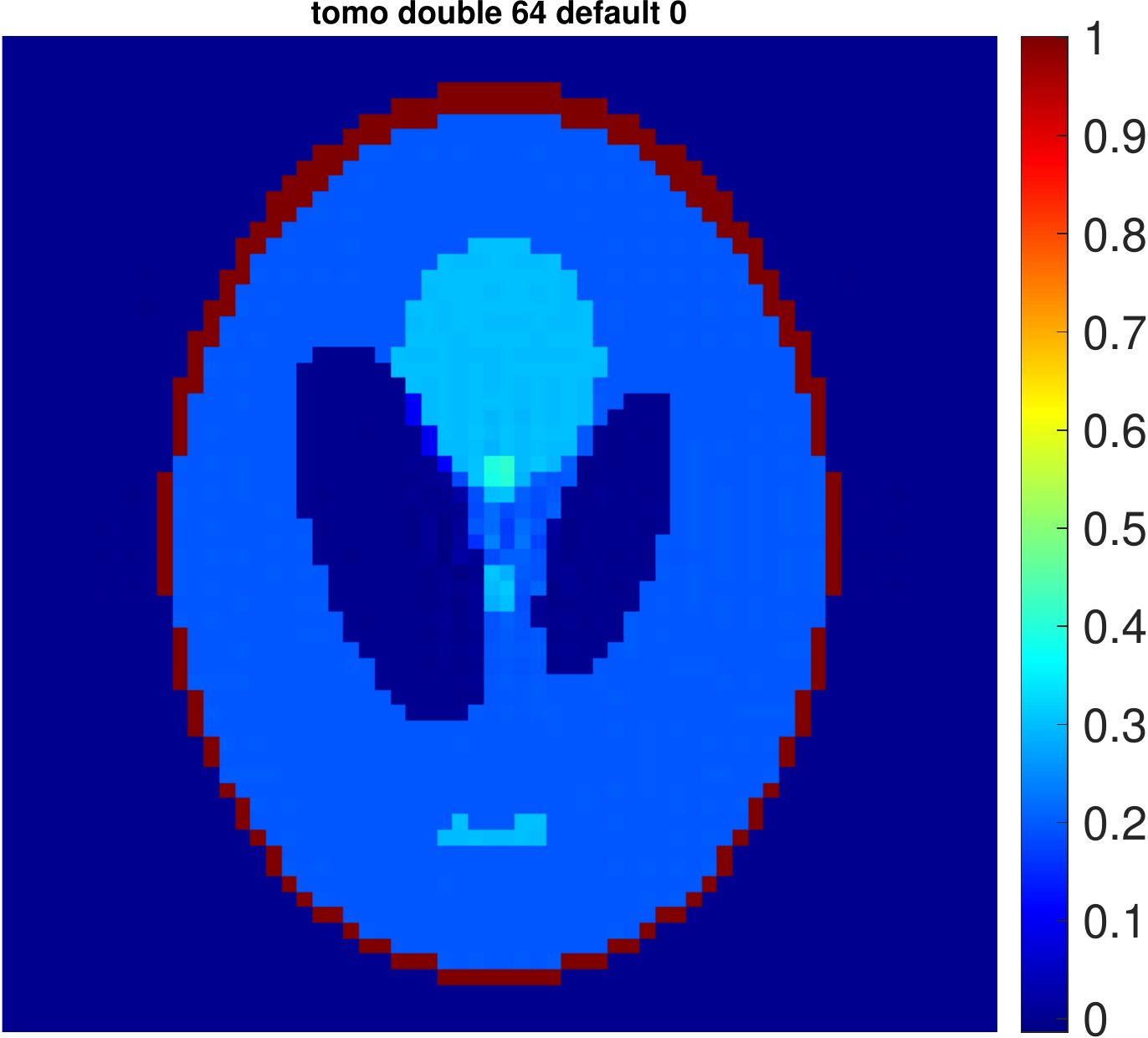}
  \caption{CGLS with regularization, double precision, zero noise.}\label{fig:Tcgreg_double_0noise_64}
\endminipage\hfill
\minipage{0.32\textwidth}
  \includegraphics[width=\linewidth]{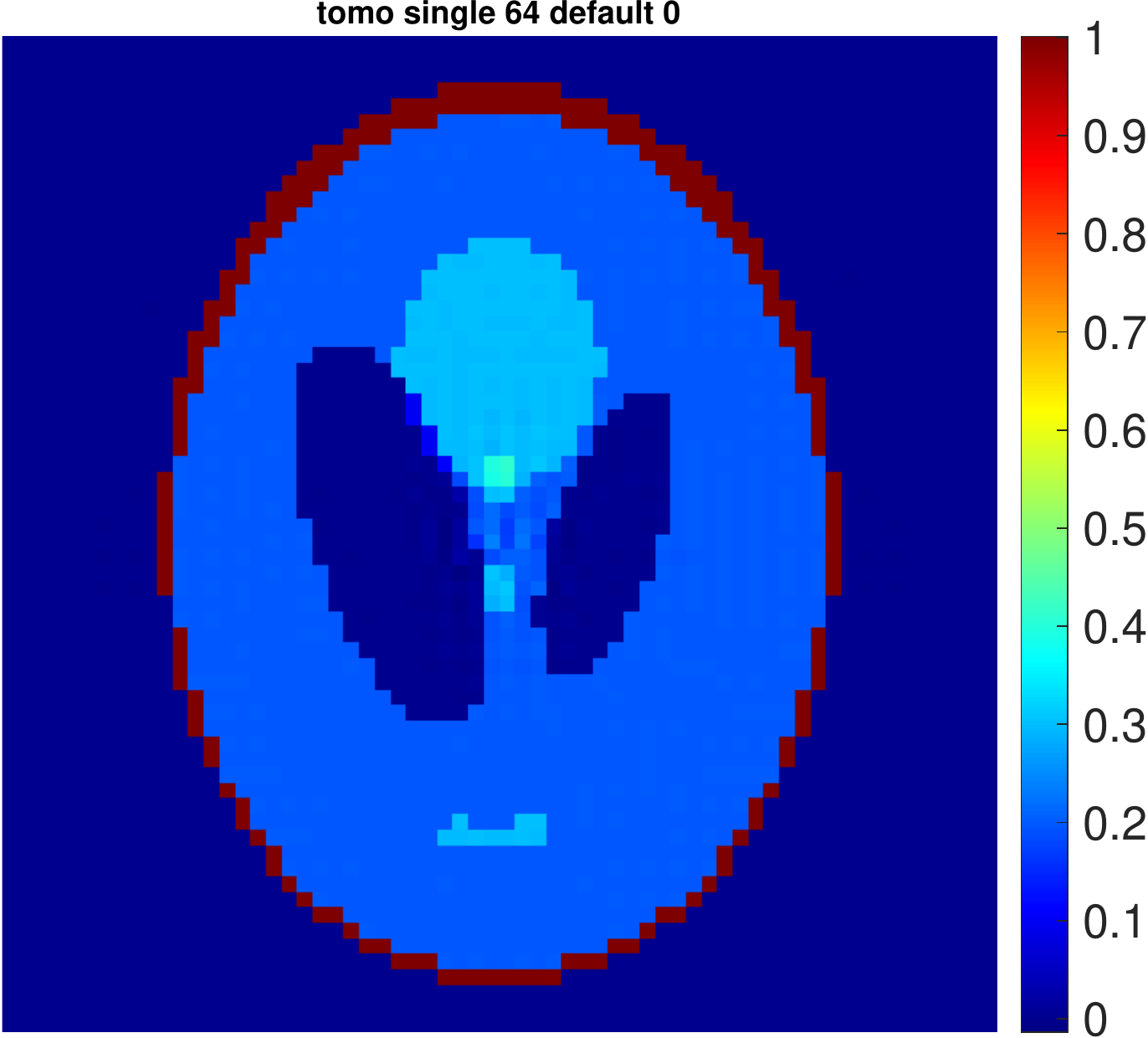}
  \caption{CGLS with regularization, single precision, zero noise.}\label{fig:Tcgreg_single_0noise_64}
\endminipage\hfill
\minipage{0.32\textwidth}%
  \includegraphics[width=\linewidth]{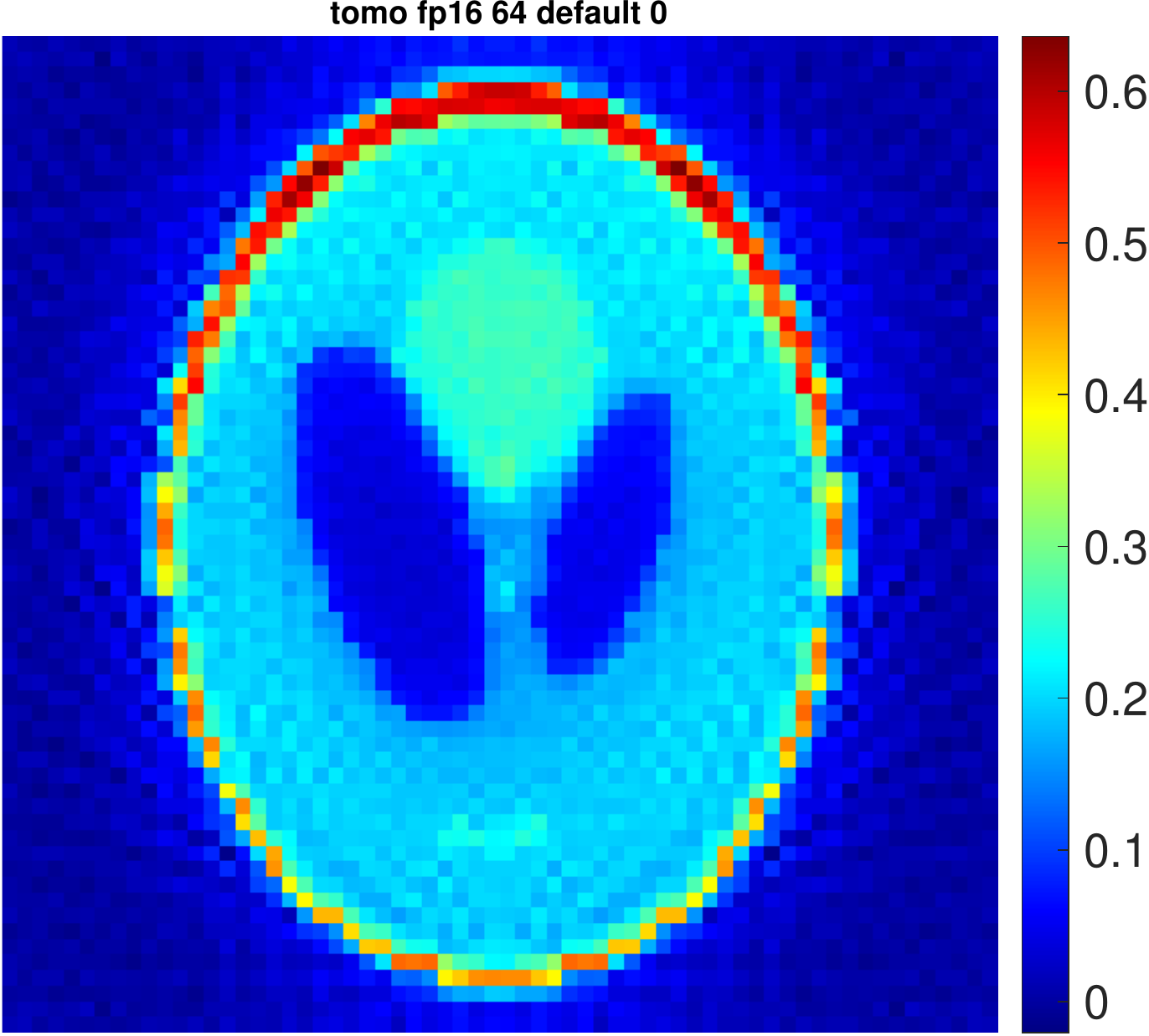}
  \caption{CGLS with regularization, half precision, zero noise (after rescaling).}\label{fig:Tcgreg_half_0noise_64}
\endminipage
\end{figure}
\begin{figure}[!htb]
\minipage{0.32\textwidth}
  \includegraphics[width=\linewidth]{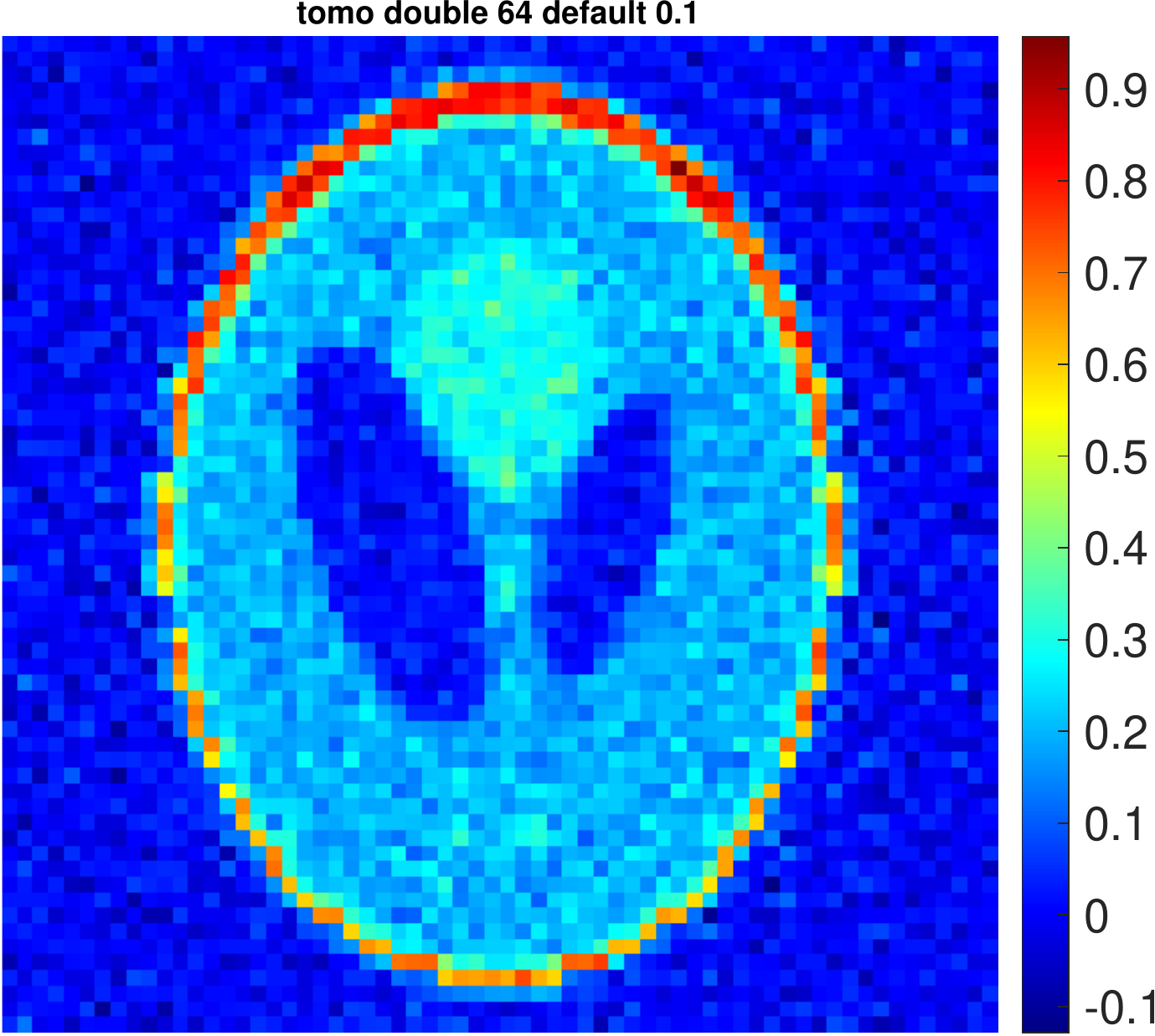}
  \caption{CGLS with regularization, double precision,10\% noise.}\label{fig:Tcgreg_double_10noise_64}
\endminipage\hfill
\minipage{0.32\textwidth}
  \includegraphics[width=\linewidth]{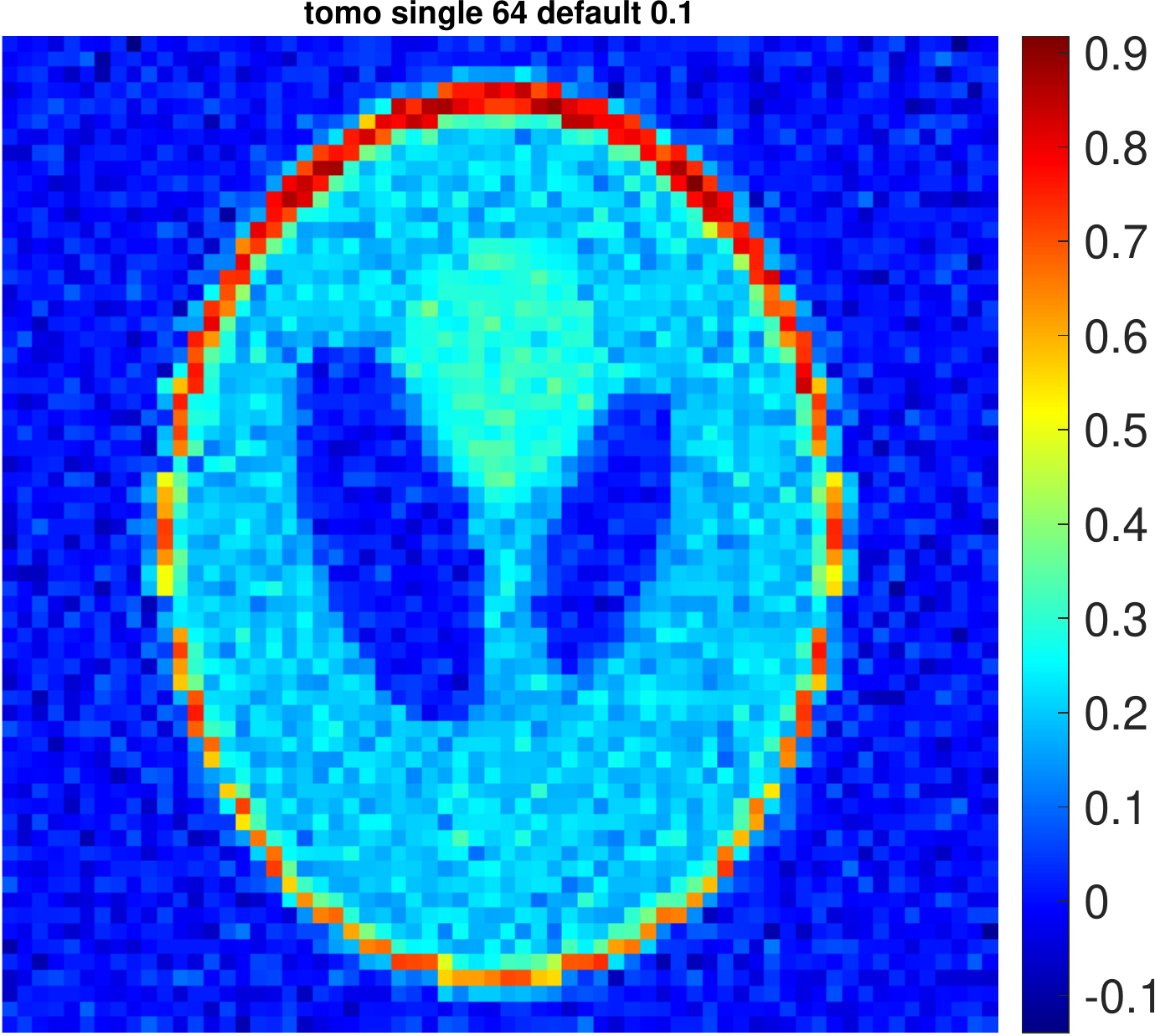}
  \caption{CGLS with regularization, single precision,10\% noise.}\label{fig:Tcgreg_single_10noise_64}
\endminipage\hfill
\minipage{0.32\textwidth}%
  \includegraphics[width=\linewidth]{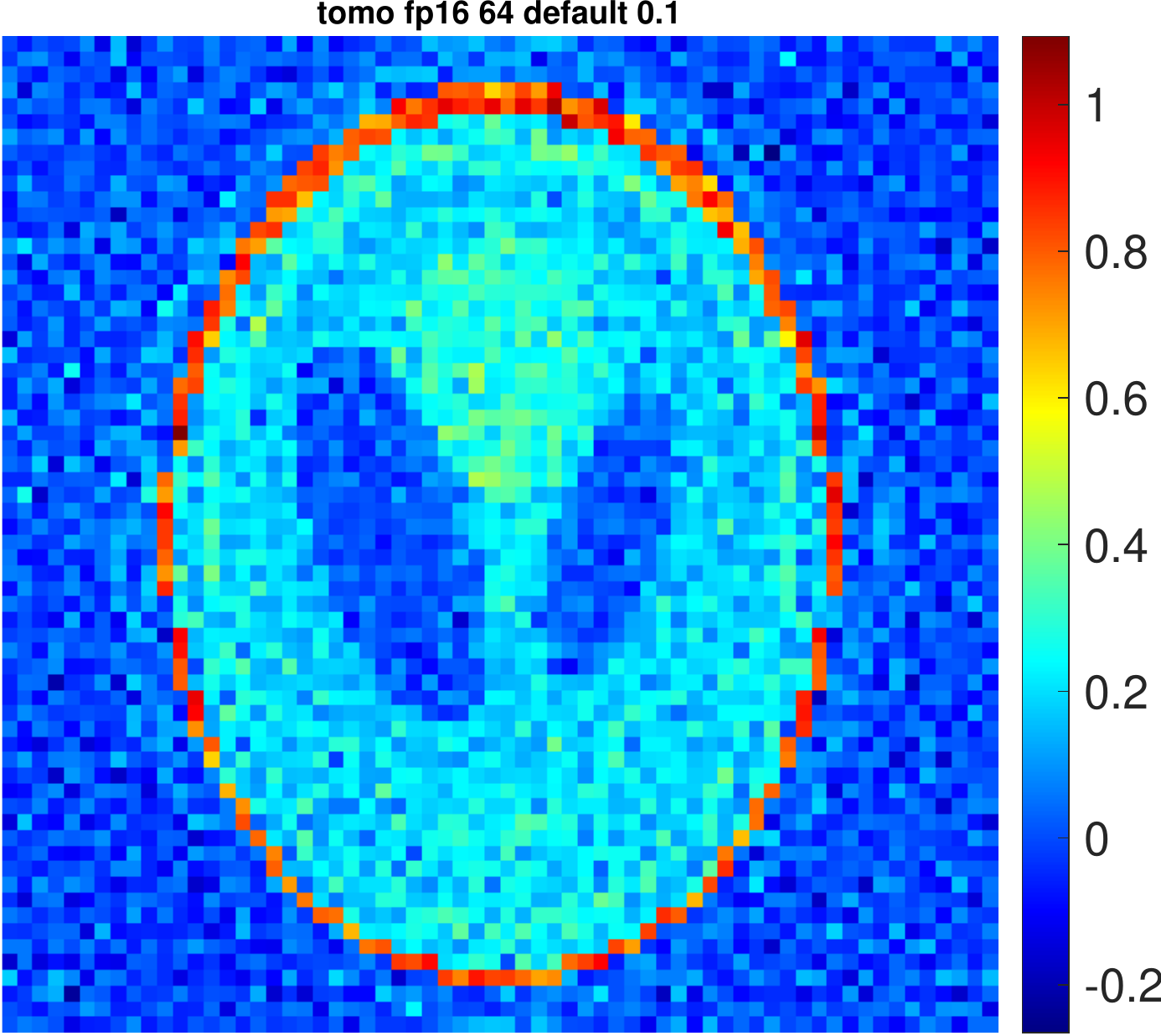}
  \caption{CGLS with regularization, half precision,10\% noise (after rescaling).}\label{fig:Tcgreg_half_10noise_64}
\endminipage
\end{figure}
\\Overall, CGLS is a more stable algorithm than CS and has less oscillations in the solution. The result given by CGLS is steadily improving in the iterations even with the presence of noise. Switching to lower precision does not have a large impact on the performance of CGLS, as long as no overflow or underflow occurs. And if they do occur, we can rescale the problem to delay the occurrence of overflow and NaNs so that a relatively good result can still be obtained. However, we do expect there to be still significant problems with overflow and NaNs for larger problems.\\
\\\noindent The main problem with CGLS at low precision is to find the suitable rescaling parameter, which depends on the matrix $A$ and vector $\bfb$ of each problem. We did not have issues of overflow/underflow at all for the image deblurring problem, while NaNs started to appear at the first iteration for the tomography reconstruction problem. We tried several rescaling parameter before we found a suitable one, and sometimes it is hard to determine whether a parameter is suitable or not if we do not have an idea of what the real image looks like. The purpose of rescaling is to avoid NaNs in the solution, but a solution without NaNs does not imply it has been rescaled properly. For example, Figure \ref{fig:size32TomoNoNoise} is the result of a size-$32$ tomography reconstruction test problem given by CGLS after $100$ iterations. No NaN occurred during the iteration process, but the resulting image is far from the true solution. This is because the overflow leads to an underflow of $\bfx$ to zero immediately in the first iteration. Though it did not result in NaNs, the algorithm did not capture any information about $\bfx$ in that iteration as well. Therefore, the image is not rescaled properly despite the absence of NaNs.
\\
\\\noindent CS, on the other hand, does not necessarily require rescaling at low precision as it avoids the calculation of inner products, which is the main source of overflow. The major advantage of CS is that it is able to run more iterations before the first NaNs occur, so that the resulting images have more chances to get refined and are less likely to be cut off before the optimal point. There are also cases when CGLS overflows at the first iteration and produces anything but NaNs, while CS is able to run several iterations and produce a relatively meaningful image, as for the tomography reconstruction test problem. In these circumstances, if only a rough estimate of the original image is desired, CS can be a more convenient choice as it saves the trouble of finding a proper rescaling parameter.\\
\\\noindent However, CS requires the user to have an estimate of the bounds of $A$'s singular values, which may be hard to estimate. Even with Tikhonov regularization, the resulting bound depends on the regularization parameter $\lambda$, which in turn depends on the noise level and the problem itself. When $\lambda$ is small, the computation becomes more risky, especially at low precision. Generally speaking, when the noise is large, we would need a larger $\lambda$ to filter out the noise. Therefore, it is not surprising that CS has better performance when the noise level is high. However, when the noise is negligible, the more stable CGLS is a better choice.

\section{Conclusion}
In this project, we explored methods for solving inverse problems at low precision. We first modified several MATLAB built-in operations using \textbf{chop} for custom precision levels and applied the blocking method to decrease the error bounds of the simulation. Then we ran the modified CGLS and CS methods on image deblurring and tomography reconstruction tasks to compare their performance across different precision levels. We added Tikhonov regularization to both methods to balance signal and noise. \\
\\After comparing results given by the two algorithms, we concluded that CGLS is more stable than CS. Its output image steadily gets closer to the real image. However, CGLS is likely to suffer from overflow at low precision as it involves the calculation of inner products. One solution to this issue is to rescale the problem, but finding a suitable rescaling parameter requires trial and error. Moreover, the rescaling is unlikely to be effective for extremely large problems when using low precision. In the future, we hope to find a more direct way that finds a suitable rescaling parameter given the matrix $A$, vector $\bfb$, and precision level. Based on intuition gained from experiments, a factor that ``rescales" the values to one digit is often a good choice. We would also look at algorithms for normalizing $A$ as described in the survey by Higham and Mary~\cite{higham2022mixed}.\\
\\The performance of the CS method is less stable and depends more on the noise level of the problem. It has more oscillations than CGLS during the iterating process. However, CS requires no calculation of inner products, and it is therefore less likely to overflow. The method avoids the need of rescaling for some problems; with rescaling, it can refine the result with even more iterations and often ends up with a better resulting image than CGLS. The price paid here is that CS needs some prior knowledge of the range of the matrix's singular values. When the matrix is large and sparse, for low precision levels, we need a regularization parameter large enough to obtain a valid lower bound and avoid the accumulation of round-off errors. For problems with noise on the right hand side $\bfb$, most of the time we would naturally end up with a large enough $\lambda$. Therefore, CS has great performance for noisy problems. However, when the problem is noise-free, the estimated $\lambda$ given by the parameter selection methods are often too small, and CS performs poorly. In the future, we hope to develop better methods for choosing suitable regularization parameters for CS at low precision.\\
\\When examining the error norms of the solution at each iteration, we noticed that sometimes the error norm did not match the so-called ``eyeball norm." An image with clear shape and background may have a higher error norm than a blurry, noisy image. The latter one, though with a smaller error norm, is obviously less informative than the former one. Therefore, more research is needed on ways to take other aspects of the output image into the account of error measurement so that we could find an image that conveys the most information.\\
\\Moreover, we would like to explore more iterative methods other than CGLS and CS and implement them in low precision, as well as mixed precision, which is expected to have both the accuracy of high precision and the gains of speed from low precision.

\bibliographystyle{abbrv}
\newpage
\bibliography{main}

\begin{thebibliography}{10}

\bibitem{bjorck1996numerical}
{\AA}.~Bj{\"o}rck.
\newblock {\em Numerical Methods for Least Squares Problems}.
\newblock SIAM, 1996.

\bibitem{Eps08}
C.~L. Epstein.
\newblock {\em Introduction to the Mathematics of Medical Imaging, Second
  Edition}.
\newblock SIAM, Philadelphia, PA, 2007.

\bibitem{gazzola2019ir}
S.~Gazzola, P.~C. Hansen, and J.~G. Nagy.
\newblock {IR} tools: a {MATLAB} package of iterative regularization methods
  and large-scale test problems.
\newblock {\em Numerical Algorithms}, 81(3):773--811, 2019.

\bibitem{nvidiatechnicalblog_2021}
G.~Gupta.
\newblock Using tensor cores for mixed-precision scientific computing.
\newblock NVIDIA Technical Blog,
  \href{https://developer.nvidia.com/blog/tensor-cores-mixed-precision-scientific-computing/}{https://developer.nvidia.com/blog/tensor-cores-mixed-precision-scientific-computing/},
  Oct-2021.

\bibitem{gutknecht2002chebyshev}
M.~H. Gutknecht and S.~R{\"o}llin.
\newblock The {C}hebyshev iteration revisited.
\newblock {\em Parallel Computing}, 28(2):263--283, 2002.

\bibitem{hansen2010discrete}
P.~C. Hansen.
\newblock {\em Discrete Inverse Problems: Insight and Algorithms}.
\newblock SIAM, 2010.

\bibitem{hestenes1952methods}
M.~R. Hestenes and E.~Stiefel.
\newblock Methods of conjugate gradients for solving.
\newblock {\em Journal of research of the National Bureau of Standards},
  49(6):409, 1952.

\bibitem{higham2002accuracy}
N.~J. Higham.
\newblock {\em Accuracy and Stability of Numerical Algorithms}.
\newblock SIAM, 2002.

\bibitem{higham2022mixed}
N.~J. Higham and T.~Mary.
\newblock Mixed precision algorithms in numerical linear algebra.
\newblock {\em Acta Numerica}, 31:347--414, 2022.

\bibitem{higham2019simulating}
N.~J. Higham and S.~Pranesh.
\newblock Simulating low precision floating-point arithmetic.
\newblock {\em SIAM Journal on Scientific Computing}, 41(5):C585--C602, 2019.

\bibitem{meng2014lsrn}
X.~Meng, M.~A. Saunders, and M.~W. Mahoney.
\newblock {LSRN}: A parallel iterative solver for strongly over-or
  underdetermined systems.
\newblock {\em SIAM Journal on Scientific Computing}, 36(2):C95--C118, 2014.

\bibitem{san2021low}
P.~San~Juan, R.~Rodr{\'\i}guez-S{\'a}nchez, F.~D. Igual, P.~Alonso-Jord{\'a},
  and E.~S. Quintana-Ort{\'\i}.
\newblock Low precision matrix multiplication for efficient deep learning in
  nvidia carmel processors.
\newblock {\em The Journal of Supercomputing}, 77(10):11257--11269, 2021.

\bibitem{trefethen1997numerical}
L.~N. Trefethen and D.~Bau~III.
\newblock {\em Numerical Linear Algebra}, volume~50.
\newblock SIAM, 1997.

\bibitem{wang2015chebyshev}
H.~Wang.
\newblock A {C}hebyshev semi-iterative approach for accelerating projective and
  position-based dynamics.
\newblock {\em ACM Transactions on Graphics (TOG)}, 34(6):1--9, 2015.

\end{thebibliography}
\end{document}